\documentclass[mathpazo]{cicp}

\makeindex

\usepackage{makeidx}         
\usepackage{graphicx}        
\usepackage{amsmath,amssymb}

\usepackage[american]{babel} 

\usepackage[babel]{csquotes} 
\usepackage{dsfont} 
  \newcommand{\IP}{\ensuremath\mathds{P}}
  \newcommand{\IR}{\ensuremath\mathds{R}}
  \newcommand{\IN}{\ensuremath\mathds{N}}
  \newcommand{\IU}{\ensuremath\mathds{U}}
  
  \newcommand{\IH}{\ensuremath\mathds{H}}
  \newcommand{\IQ}{\ensuremath\mathds{Q}}

\usepackage{tikz}
\usepackage{tabularx}             
  \newcolumntype{R}{>{\raggedleft\arraybackslash}X}  \newcolumntype{L}{>{\raggedright\arraybackslash}X}  \newcolumntype{C}{>{\centering\arraybackslash}X}
\usepackage{multirow} 
\usepackage{rotating} 
\usepackage{booktabs} 
\usepackage{bm}    
\usepackage{mathtools}
\usepackage{fixmath}
\usepackage{hyperref}
  \hypersetup{breaklinks=true,colorlinks=true,pdfborder={0 0 0}}
  \hypersetup{linkcolor=black, anchorcolor=black, citecolor=black, filecolor=black, urlcolor=black}
\usepackage{subfigure}
\usepackage{enumitem}
\newcommand{\elem}{{T}}
\newcommand{\Proj}{{\Pi}}
\usepackage[linesnumbered,ruled,vlined]{algorithm2e}
\newcommand*{\setT}{\ensuremath{\mathcal{T}}}
\newcommand*{\setE}{\ensuremath{\mathcal{E}}}
\newcommand*{\setV}{\ensuremath{\mathcal{V}}}

\newcommand*{\vphi}{\varphi}                                     
\renewcommand*{\vec}[1]{{\boldsymbol{#1}}}                       
\newcommand*{\mapEE}{\hat{\vec{\vartheta}}}                      
\DeclareMathAlphabet{\mathbfsf}{\encodingdefault}{\sfdefault}{bx}{n}
\newcommand*{\vecc}[1]{\mathbfsf{#1}}                            
\newcommand*{\grad}{\vec{\nabla}}                                
\renewcommand*{\div}{\vec{\nabla}\cdot}                          
\newcommand*{\dd}{\mathrm{d}}                                    
\newcommand*{\abs}[1]{\ensuremath{\left|#1\right|}}              
\newcommand*{\card}[1]{\ensuremath{\##1}}                        
\newcommand*{\transpose}[1]{{#1}^\mathrm{T}}                     
\newcommand*{\jump}[1]{[[#1]]}
\DeclareMathOperator*{\diag}{diag}                               
        
\newcommand*{\llbrace}{\left\lbrace\!\left\vert}\newcommand*{\rrbrace}{\right\vert\!\right\rbrace}
  \newcommand*{\avg}[1]{\llbrace{#1}\rrbrace}                    
\newcommand{\I}{I} 
\newcommand{\II}{I\!I} 
\newcommand{\III}{I\!I\!I} 
\newcommand{\IV}{I\!V} 
\newcommand{\V}{V} 
\newcommand{\VI}{V\!I} 
\newcommand{\VII}{V\!I\!I} 
\newcommand{\VIII}{V\!I\!I\!I} 
\newcommand{\IX}{I\!X} 
\newcommand{\X}{X} 
\newcommand{\XI}{X\!I} 
\newcommand{\XII}{X\!I\!I} 
\newcommand{\XIII}{X\!I\!I\!I} 
\newcommand{\XIV}{X\!I\!V} 
\newcommand{\XV}{X\!V}
\newcommand{\XVI}{X\!V\!I}
\newcommand{\XVII}{X\!V\!I\!I}
\newcommand{\XVIII}{X\!V\!I\!I\!I}
\newcommand{\XIX}{X\!I\!X} 
\newcommand{\XX}{X\!X} 
\newcommand{\XXI}{X\!X\!I} 
\newcommand{\XXII}{X\!X\!I\!I} 
\newcommand{\XXIII}{X\!X\!I\!I\!I} 
\newcommand{\XXIV}{X\!X\!I\!V}

\newcommand{\tO}{\widetilde{\Omega}}      
\newcommand{\thh}{{\tilde{h}}}            
\newcommand{\tq}{{\tilde{\vec{q}}}}       
\newcommand{\tf}{{\tilde{f}}}             
\newcommand{\tgN}{{\tilde{g}_\mathrm{N}}} 
\newcommand{\tHH}{{\tilde{H}}}            
\newcommand{\tQ}{{\tilde{Q}}}             
\newcommand{\specS}{S_{\!0}}              
\newcommand{\zb}{\zeta_\mathrm{b}}
\newcommand{\Hs}{H_\mathrm{s}}

\newcommand{\Err}[1]{{\ensuremath{\text{\scriptsize Err}\left(#1\right)}}}
\newcommand{\EOC}[1]{{\ensuremath{\text{\scriptsize EOC}\left(#1\right)}}}

\definecolor{faublue}{RGB}{0,56,101}
\definecolor{faulblue}{RGB}{114,159,207}
\definecolor{faullblue}{RGB}{221,229,240}

\definecolor{faugrey}{RGB}{150,150,150}
\definecolor{faulgrey}{RGB}{210,213,215}
\definecolor{faullgrey}{RGB}{235,236,238}

\definecolor{faunat}{RGB}{0,155,119}
\definecolor{faulnat}{RGB}{170,207,189}
\definecolor{faullnat}{RGB}{229,239,234}

\newcommand{\Matlab}{\mbox{MATLAB}}
\newcommand{\Octave}{GNU~\mbox{Octave}}
\newcommand{\MatOct}{\Matlab\,/\,\allowbreak\Octave} 

\usepackage{textcomp} 
\usepackage{listings}
  \definecolor{keywordcolor}{rgb}{0, 0.25, 0.5}
  \definecolor{commentcolor}{rgb}{0.2, 0.5, 0.2}
  \definecolor{stringcolor}{rgb}{0.5, 0.5, 0.2}
\newcommand*{\code}[1]{\mbox{\lstinline[basicstyle=\ttfamily\small]{#1}}}

\newcommand*{\revised}[1]{#1}

\begin{document}

\raggedbottom

\title{FESTUNG: A~\MatOct~toolbox for the discontinuous Galerkin method. Part IV: Generic problem framework and model-coupling interface}

\author[Reuter B et.~al.]{Balthasar~Reuter\affil{1}, Andreas~Rupp\affil{2}\comma\affil{1}, Vadym~Aizinger\affil{3}\comma\affil{1}\comma\corrauth, Florian~Frank\affil{1}, and Peter~Knabner\affil{1}}



\address{ \affilnum{1}\ Friedrich--Alexander University of Erlangen--N\"urnberg, Department of Mathematics, 
Cauerstra{\ss}e~11, 91058~Erlangen, Germany. \\
          \affilnum{2}\ Ruprecht-Karls-Universit\"at Heidelberg, Interdisciplinary Center for Scientific Computing, Im Neuenheimer Feld 205, 69120 Heidelberg, Germany.\\
          \affilnum{3}\ University of Bayreuth, Chair of Scientific Computing, 95447~Bayreuth, Germany. }


\emails{ \texttt{reuter@math.fau.de} (B.~Reuter), \texttt{rupp@math.fau.de} (A.~Rupp), \texttt{vadym.aizinger@uni-bayreuth.de} (V.~Aizinger), \texttt{florian.frank@fau.de} (F.~Frank), \texttt{knabner@math.fau.de} (P.~Knabner) }

\date{Received: \today / Accepted: ...}


\begin{abstract}
This is the fourth installment in our series on implementing the discontinuous Galerkin (DG) method as an~open source \MatOct~toolbox. Similarly to its predecessors, this part presents new features for application developers employing DG~methods and follows our strategy of relying on fully vectorized constructs and supplying a~comprehensive documentation. The specific focus of the current work is the newly added generic problem implementation framework and the highly customizable model-coupling interface for multi-domain and multi-physics simulation tools based on this framework. The functionality of the coupling interface in the FESTUNG toolbox is illustrated using a~two-way coupled free-surface\,/\,groundwater flow system as an~example application.
\end{abstract}




\maketitle

 
\section{Introduction}

The previous papers in the FESTUNG series dealt with the most common differential operators such as the linear diffusion~\cite{FESTUNG1} or advection~\cite{FESTUNG2} as well as with different types of discontinuous Galerkin (DG) discretizations, namely the standard local DG (LDG~\cite{FESTUNG1}) or hybridized DG (HDG~\cite{FESTUNG3}). 
In the time span between the previous installment~\cite{FESTUNG3} and the current work, several applications involving non-linear equations were implemented using our \MatOct~FESTUNG~\cite{FESTUNGGithub} toolbox; those include the two-dimensional shallow-water equations~\cite{HajdukHAR2018} and mean curvature flow~\cite{BungertAF2017,AizingerBF2018}. \revised{In addition, a~general overview paper documenting the current state of development of FESTUNG has been submitted~\cite{FESTUNG-SI} and is currently under review.}
The present study, however, sets a~much more ambitious goal: 
Presenting a~model-coupling interface that allows to create very complex simulation systems consisting of multiple self-contained simulation tools and interacting by exchanging data usually in the form of solution vectors. 
The more specific objectives of this work include developing an~abstract coupling concept capable of supporting very general types of model coupling as well as implementing and testing this concept in the framework of FESTUNG. 
Just as in the previous parts of this series, our performance-optimized, fully vectorized implementation with detailed documentation is freely available as an~open source software.
\par
Very few physical systems are truly isolated and thus can be modeled in a~standalone fashion without the help of some more or less realistic assumptions. 
Often, the setting can be simplified, and action and reaction of other physical systems can be accounted for in a single-physics model via boundary conditions, forcing terms, parametrizations, etc. 
However, many important problems do not lend themselves readily to such simplifications; hence the need for coupled models. 
Examples of widely used multi-physics applications include, in particular, ocean-atmosphere-land systems constituting the staple of climate modeling or fluid-structure interaction models very common in civil engineering and medical applications. 
This multi-physics label can sometimes be somewhat misleading (see a~discussion of this issue in a~very comprehensive overview paper~\cite{Keyes2013}) and is also often used to describe interactions between different mathematical representations, or numerical schemes, or grids, or motion scales within the same `physics'.
\par
For coupled multi-physics problems, one can generally identify three main classes of setups: shared domain/ multiple physics (e.\,g.\ coupled subsurface flow/geomechanics), multiple domain/shared physics (e.\,g.\ regional grid nesting in ocean or atmosphere simulations), and multiple domain/multiple physics. 
The first two permit a~number of performance-relevant simplifications; although our new problem implementation framework accommodates either, here we focus on the third, the most general type of setup, for which a~number of key aspects have to be considered:
\begin{itemize}[nosep]
\item {\bf Modeling issues:} Those include the physical and mathematical consistency questions such as conservation properties, physically meaningful interface conditions, well-posedness of the mathematical problem, etc.

\item {\bf Numerical issues:} Different mathematical models often require specialized numerical methods; in addition, computational meshes corresponding to separate physical systems do not necessarily match at the interface (or coincide in the case of a~shared domain setup). 
This gives rise to a~number of challenges such as controlling interpolation and conservation errors between discrete solutions, matching meshes at the model interface or using, e.\,g., mortars~\cite{Maday1989} to bridge between sub-domain meshes. 
Another host of rather difficult issues arises when one attempts to analyze and to improve the stability and the accuracy of the discretizations for such coupled models. 
These difficulties concern the stability and the convergence of spatial discretizations and---in the case of time-dependent problems---can also affect the temporal convergence that may be degraded by operator splitting techniques very common in coupled applications.

\item {\bf Algorithmic issues:} In addition to differences in the discretization, the solution algorithms of each sub-model may be such that a~well thought-out solution strategy becomes necessary for the coupled model. 
This aspect is particularly important for time-dependent applications, where different PDE system types and resolved spatial/temporal scales may require a~careful handling of the time stepping procedure. 
The algorithmic approaches encompass the one-way (when the data/information only flows from one model to the other) and the two-way coupling paradigms. 
The latter is much more general and includes three main classes of schemes: 
\begin{itemize}[nosep]
\item Fully implicit (monolithic)---all sub-problems together with the interface conditions are treated as one large system solved without operator splitting. This approach incurs high implementational and computational expenses but often produces the most robust and physically consistent schemes. 
It is the method of choice for stationary problems but it is also known to have been used for time-dependent applications---specifically to deal with highly non-linear systems of PDEs. 
It also tends to allow for (significantly) larger time-steps than the other coupling schemes. 
The main difficulties are connected with computing the Jacobi matrix for the coupled system (or solving the non-linear system in a~matrix-free fashion)---a~difficult task if the sub-models are complex and implemented in separate software packages.

\item Iterative (internal) coupling---all sub-systems are solved in each coupling step alternately until a given convergence criterion is satisfied. 
This methodology represents a~compromise between the fully implicit and non-iterative couplings and is widely used in many applications. 
This approach allows to use a~specialized software package for each sub-model and provides a~mechanism to control the coupling (i.e., operator splitting) error. 
An~iterative couping requires, however, particular care with regard to the conservation and stability properties of the discrete scheme; in addition, the operator splitting may cause convergence to a~non-physical solution if the underlying systems have, e.\,g., non-smooth coefficients~\cite{Keyes2013}.

\item Non-iterative (external) coupling works similarly to the iterative coupling but forgoes iterations between the sub-models. 
It is simple and cheap but not sufficiently robust and accurate for some applications.
\end{itemize}

\item {\bf Software issues:} The last but not least difficulty arising when implementing coupled models concerns the software packages used for each sub-model and the interface between them. 
The most advanced and user-friendly features are offered by general-purpose modeling frameworks such as DUNE\footnote{\url{https://www.dune-project.org}}, FEniCS\footnote{\url{https://www.fenicsproject.org}}, or similar software projects.
For software packages not originating from the same family, a~number of specialized couplers exist that support data exchange and asynchronous execution modes, e.\,g. OASIS\footnote{\url{https://portal.enes.org/oasis}}. 
\end{itemize}

Due to their conceptual, algorithmic, and computational complexity, coupled applications have been almost exclusively a~domain of large groups/companies/institutions and tend to be developed for a~specific class of applications with the corresponding requirements in mind. 
Thus the coupled numerical models of ocean/atmosphere/land used for climate studies are usually optimized for massively parallel execution and efficient handling of large volumes of forcing and grid data; they mostly rely on a~loose (external) coupling between single-physics sub-models~\cite{Valcke2012} realized in a~form of a~dedicated piece of software such as OASIS or MCT~\cite{Larson2005}. 
Fluid-structure interaction or groundwater flow-geomechanics models, on the other hand, are often designed from scratch in a monolithic fashion and thus solve a~fully coupled non-linear equation system for all sub-models at once~\cite{RinTGVT2017}.

One of the main motivations for the current work is to make the model coupling and multi-physics (in a~very general sense) capability available to users without access to sophisticated software and high performance computing resources and to enable fast prototyping and testing of multi-physics applications. 
We limit the focus to the algorithmic and software aspects of model coupling; the numerical issues are rather straightforward here due to DG discretizations employed in all sub-models; modeling issues will be presented in a brief form. 
As an~example application to illustrate the coupling mechanism, we consider a~free-surface/groundwater flow problem. 
Free-surface flow is represented by the \emph{2D shallow water equations} in a~vertical slice (2Dv SWE---also known as primitive hydrostatic equations in a~vertical slice), whereas the groundwater is modeled using \emph{Darcy's law} (also in a~two-dimensional vertical slice). 
A~stability analysis of this discretization using the same interface conditions is presented in~\cite{ReuterRAK2019}.

The rest of this article is structured as follows:
The model problem is presented in Section~\ref{sec:model} accompanied by its LDG discretization in Section~\ref{sec:discretization}.
The generic problem framework is the subject of Section~\ref{sec:problem-framework} (with explicit forms of the terms in the linear system given in Appendix~\ref{sec:free-flow:mat-vec}), and an~implementation of the model problem in the context of this framework is demonstrated in Section~\ref{sec:coupling} (more implementation details are provided in Appendix~\ref{sec:assembly}).
Numerical examples can be found in Section~\ref{sec:results}, and 
a~conclusion and outlook make up the remainder of this paper.

\section{Model problem} 
\label{sec:model}

\subsection{Computational domain and mesh}\label{sec:model:mesh}
\begin{figure}[h]
\center
\includegraphics{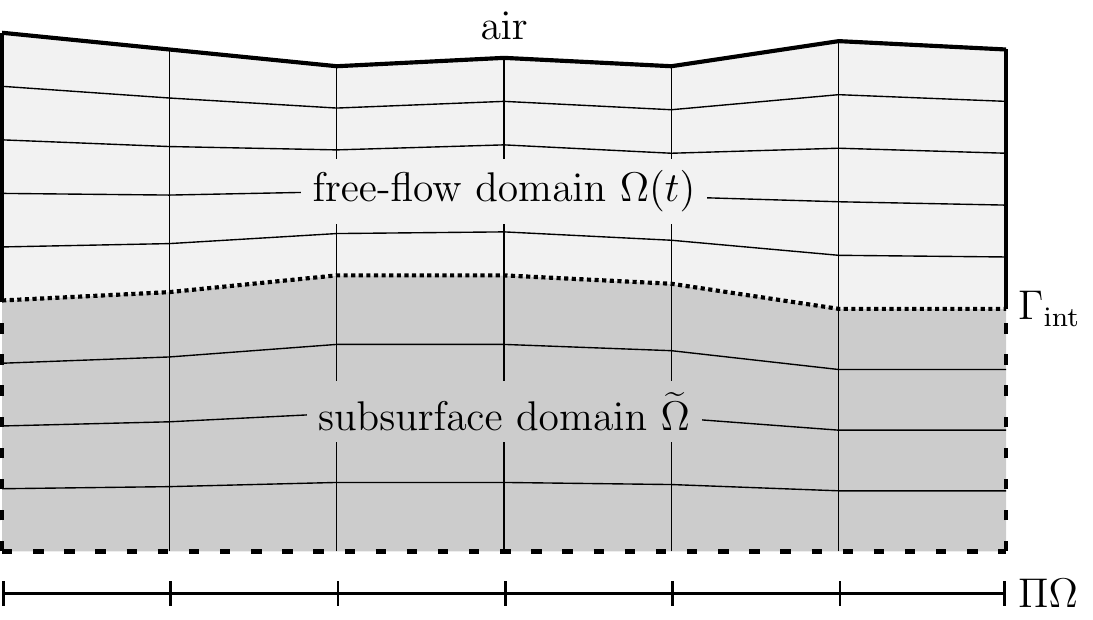}
\caption{Coupled domain $\Omega(t)\uplus \tO$ with interior boundary $\Gamma_\mathrm{int}$, mesh $\mathcal T_{\Delta}$, and projected one-dimensional domain $\Proj\Omega$ and mesh.}\label{fig:Omega}
\end{figure}

Let~$J\coloneqq\,(0,t_\mathrm{end})\,$ be a~finite time interval and~$\Omega(t) \uplus \tO$ the coupled domain consisting of a~free-flow subdomain~$\Omega(t)$ on top of the subsurface subdomain $\tO$ (cf.~Fig.~\ref{fig:Omega}), both of which with compact closure and assumed to be polygonally bounded and Lipschitz.  
Free-flow and subsurface subdomains are separated by an~interior (transition) boundary~$\Gamma_\mathrm{int}$.
Let $\setT_{\Delta}$ be a~regular family of non-overlapping partitions of~$\Omega \in \{\Omega(t), \tO\}$ into~$K$ closed elements~$T$ such that $\displaystyle \Omega=\cup T$.
The discretization of the free-flow problem~\eqref{eq:free:nonmixed} requires a~mesh consisting of trapezoidal elements with strictly vertical parallel edges (cf.~Fig~\ref{fig:Omega}). Under this condition, the two-dimensional elements are aligned in vertical columns, each corresponding to a~one-dimensional element on the $x^1$-axis later used to discretize the water height.
For simplicity, the subsurface problem also uses a~mesh of the same type. 
Denoting by $\Proj$ the standard orthogonal projection operator onto the $x^1$-axis, this geometry and arrangement of mesh elements produces a~non-overlapping partition of $\Proj \Omega$ denoted by $\Proj \setT_{\Delta} $ by simply projecting elements of $\setT_{\Delta}$.

\subsection{Subsurface problem}
\label{sec:model:pm}
Let the boundary of~$\tO$ be subdivided into Dirichlet~$\partial\tO_\mathrm{D}$ and Neumann~$\partial\tO_\mathrm{N}$ parts. 
We consider the \emph{time-dependent Darcy equation}
$\specS\partial_t \thh - \nabla \cdot \big(\vecc{K} \, \nabla \thh\big) = f$
describing water transport through fully saturated porous media, where $\thh$ is generally understood as the hydraulic head.  The constant coefficient~$\specS$ denotes the \emph{specific storativity} of the porous medium and $\vecc{K}$ the \emph{hydraulic conductivity}.
Division by~$\specS$ and setting~$\vecc{\widetilde{D}}\coloneqq \vecc{K}/\specS$, $\tf\coloneqq f/\specS$ yields
\begin{subequations}\label{eq:sub:nonmixed}
\begin{align}
\partial_t \thh(t, \vec{x})  -  \div \Big( \vecc{\widetilde D}(t, \vec{x})\, \grad \thh(t, \vec{x}) \Big) 
&\;=\; \tf(t, \vec{x}) &&\text{in}~J\times\tO~,\label{eq:sub:nonmixed:u}\\
\intertext{where $\vecc{\widetilde D}:J\times\tO\rightarrow \IR^{2\times2}$ (a~uniformly symmetric positive definite matrix) and $\tf:J\times\tO\rightarrow \IR$ are considered time and space dependent. 
Eq.~\eqref{eq:sub:nonmixed:u} is complemented by the following boundary and initial conditions:}
\thh(t,\vec{x})  &\;=\; \thh_\mathrm{D}(t,\vec{x})  &&\text{on}~J\times{\partial\tO}_{\mathrm{D}}\;, &&\label{eq:sub:nonmixed:D}\\
- \vecc{\widetilde D}(t,\vec{x})\, \grad \thh(t,\vec{x})\cdot\vec{\nu}   &\;=\; \tgN(t,\vec{x}) &&\text{on}~J\times{\partial\tO}_\mathrm{N}\;,\label{eq:sub:nonmixed:N}&&\\
\thh(0,\vec{x}) &\;=\; \thh_0 (\vec{x})           &&\text{in}~\tO\,,&&
\end{align}
\end{subequations}
where $\vec{\nu}$ denotes the outward unit normal, and 
$\tilde{h}_0$ is the given initial and $\tilde{h}_\mathrm{D}, \;\tilde{g}_\mathrm{N}$ the boundary data.

\subsection{Free-flow problem}
\label{sec:model:ff}
The boundary of the free-flow domain~$\partial\Omega(t)$ is assumed to consist of top~$\partial\Omega_\mathrm{top}(t)$, bottom~$\partial\Omega_\mathrm{bot}$, and lateral sections with the latter subdivided into land~$\partial\Omega_\mathrm{land}(t)$, open sea $\partial\Omega_\mathrm{os}(t)$, river~$\partial\Omega_\mathrm{riv}(t)$, and radiation~$\partial\Omega_\mathrm{rad}(t)$ parts. All boundaries except $\partial\Omega_\mathrm{bot}$ are regarded as time-dependent, but, for brevity, the time variable is omitted.
Then the 2Dv shallow water equations for velocity $\vec{u}(t,\vec{x}) = \transpose{[u^1, u^2]}$ and total water height $h(t, x^1) = \xi(t, x^1) - \zb(x^1)$ are given as
\begin{subequations}\label{eq:free:nonmixed}
\begin{align}
&\partial_t h(t, x^1) + \partial_{x^1} \bigg( \int_{\zb(x^1)}^{\xi (t, x^1)} u^1(t, \vec{x}) \, \dd x^2 \bigg) = 0 &&\text{in}~J\times \Proj \Omega \;, \label{eq:free:nonmixed:surf}\\
&\begin{aligned}
  \partial_t u^1(t, \vec{x}) + \div \left( u^1(t, \vec{x}) \, \vec{u}(t, \vec{x}) \right) + g\,\partial_{x^1} h(t, \vec{x}) \; & \\
  - \div \left( \vecc{D}(t,\vec{x}) \, \grad u^1(t, \vec{x}) \right) &= f(t, \vec{x}) - g \, \partial_{x^1} \zb (t, x^1)
\end{aligned} &&\text{in}~J\times\Omega(t) \;, \label{eq:free:nonmixed:momentum}\\
&\div \vec{u}(t, \vec{x}) = 0 &&\text{in}~J\times\Omega(t) \;, \label{eq:free:nonmixed:cont}
\end{align}
where $g$~denotes the acceleration due to gravity, $f$~is a~source term, $\xi$ and $\zb$ are the free-surface elevation and the bathymetry with respect to some datum, respectively, and $\vecc{D}$ is a~diffusion tensor satisfying the same conditions as those imposed on $\vecc{\widetilde D}$.
The following boundary and initial conditions are specified for Eqs.~\eqref{eq:free:nonmixed:surf}--\eqref{eq:free:nonmixed:cont}:
\begin{align}
u^1(t,\vec{x}) &\;=\; u^1_\mathrm{D}(t,\vec{x}) &&\text{on}~J\times  ( \partial\Omega_\mathrm{land} \cup \partial\Omega_\mathrm{riv} \cup \partial\Omega_\mathrm{bot} )\,,\label{eq:free:nonmixed:bc:u1}\\
u^2(t,\vec{x}) &\;=\; u^2_\mathrm{D}(t,\vec{x}) &&\text{on}~J\times\partial\Omega_\mathrm{bot}\,,\label{eq:free:nonmixed:bc:u2}\\
- \vecc{D}(t,\vec{x})\, \grad u^1(t,\vec{x})\cdot\vec{\nu} &\;=\; q_\mathrm{D}(t,\vec{x}) &&\text{on}~J\times ( \partial\Omega_\mathrm{top} \cup \partial\Omega_\mathrm{os} \cup \partial\Omega_\mathrm{rad}  )\,,\label{eq:free:nonmixed:bc:q}\\
h(t,x^1) &\;=\; h_\mathrm{D}(t,x^1) &&\text{on}~J\times\Proj  (\partial\Omega_\mathrm{os} \cup \partial\Omega_\mathrm{riv}  )\,,\label{eq:free:nonmixed:bc:h}\\
u^1(0,\vec{x}) &\;=\; u^1_0 (\vec{x}) &&\text{in}~\Omega(0)\,,\\
h(0,x^1) &\;=\; h_0 (x^1) &&\text{in}~\Proj \Omega\,.
\end{align}
\end{subequations}
To summarize the different types of boundary conditions, note the following:
\begin{enumerate}[nosep]
  \item Bottom boundary~$\partial\Omega_\mathrm{bot}$: no-slip ($\vec{u}\cdot\vec{\tau} = 0$ with tangential vector~$\vec{\tau}$) and prescribed normal flow~$\vec{u}\cdot\vec{\nu}$, which translates to~$u^1_\mathrm{D} = u^2_\mathrm{D} = 0$ in uncoupled simulations and interface condition~\eqref{eq:coupling:ff} in coupled simulations;
  \item free surface~$\partial\Omega_\mathrm{top}$ and radiation boundary~$\partial\Omega_\mathrm{rad}$: vanishing normal derivative of the flow, i.\,e., $q_\mathrm{D} = 0$;
  \item land boundary~$\partial\Omega_\mathrm{land}$: no normal flow resulting in $u^1_\mathrm{D} = 0$ due to the strictly vertical boundary;
  \item open sea boundary~$\partial\Omega_\mathrm{os}$: prescribed water height~$h_\mathrm{os}$ (e.\,g., due to tidal forcing), i.\,e., $h_\mathrm{D} = h_\mathrm{os}$ and vanishing normal derivative of the flow, $q_\mathrm{D} = 0$;
\item river boundary~$\partial\Omega_\mathrm{riv}$: prescribed horizontal velocity $u^1_\mathrm{D} = u^1_\mathrm{riv}$ and water height $h_\mathrm{D} = h_\mathrm{riv}$.
\end{enumerate}
Note that different types of boundaries should be pieced together in a~{\em compatible} manner in order to ensure well-posedness of the initial-boundary-value problem.

\subsection{Interface conditions}
\label{sec:model:coupled}
The coupling conditions at the interface boundary have been motivated and described in some detail in our analysis paper~\cite{ReuterRAK2019}. 
They impose the continuity of the normal flux~\eqref{eq:coupling:ff} and the continuity of the dynamic pressure\,/\,head~\eqref{eq:coupling:pm}.
\begin{subequations}\label{eq:coupling}
\vspace{-2mm}
 \begin{align} 
  {\textstyle\frac{1}{\specS}} \vec{u}(t,\vec{x}) \cdot \vec{\nu} & \; = \; \vecc{\widetilde D}(t,\vec{x}) \, \grad \thh(t,\vec{x}) \cdot \tilde{\vec{\nu}} & & \mathrm{on}~J \times \Gamma_\mathrm{int} \,, \label{eq:coupling:ff}\\
  \thh(t,\vec{x}) & \; = \; \xi(t,x^1) + {\textstyle\frac{1}{2\,g}} \left( u^{1} \right)^2 & & \mathrm{on}~J \times \Gamma_\mathrm{int} \,,\label{eq:coupling:pm}
 \end{align}
\end{subequations}
where $\vec{\nu}$ and $\tilde{\vec{\nu}}$ denote the outward unit normals on~$\Gamma_\mathrm{int}$ with respect to~$\Omega(t)$ and~$\tO$, correspondingly.

\section{LDG Discretization}
\label{sec:discretization}

\subsection{Variational formulation of the subsurface system}
\label{sec:variational-sub}

To formulate an LDG~scheme for system~\eqref{eq:sub:nonmixed}, we generally follow~\cite{RuppKnabner2017,RuppKD2018,AizingerRSK2018} and first introduce an~auxiliary vector-valued unknown~$\tq \coloneqq -\grad \thh$ and re-write~\eqref{eq:sub:nonmixed} in mixed form also making the necessary changes to the boundary conditions:
\begin{subequations}\label{eq:sub}
\begin{align}
\partial_t \thh(t,\vec{x})  + \div (\vecc{\widetilde D}(t,\vec{x}) \,\tq(t,\vec{x})) &\;=\; \tf(t,\vec{x})  &&\text{in}~J\times\tO\;,                           \label{eq:sub:u}\\
\tq(t,\vec{x}) + \grad \thh(t,\vec{x})        &\;=\; 0             &&\text{in}~J\times\tO\;,        \label{eq:sub:q}\\
\thh(t,\vec{x})                                                               &\;=\; \thh_\mathrm{D}(t,\vec{x})  &&\text{on}~J\times{\partial\tO}_{\mathrm{D}}\;,    \label{eq:sub:D}\\
\vecc{\widetilde D}(t,\vec{x}) \, \tq(t,\vec{x})\cdot\vec{\nu}                    &\;=\; \tgN(t,\vec{x})  &&\text{on}~J\times{\partial\tO}_\mathrm{N}\;,      \label{eq:sub:N}\\
\thh(0, \vec{x})                                                              &\;=\; \thh_0(\vec{x})             &&\text{in}~\tO\;.                     \label{eq:sub:0}
\end{align}
\end{subequations}

The discontinuous nature of DG approximations permits formulating the variational system on an~element-by-element basis. Denoting by $\vec{\nu}_T$ the unit exterior normal to~$\partial T$, we multiply both sides of equations~\eqref{eq:sub:u},~\eqref{eq:sub:q} with smooth test functions~$\vec{\tilde{y}}: T \rightarrow\IR^2$, $\tilde{w}: T \rightarrow\IR$ , correspondingly, and integrate by parts over $T\in\setT_{\Delta}$ to obtain
\begin{align*}
&\int_{T} \tilde{w}\,\partial_t \thh\,\dd\vec{x} - \int_{T}\grad \tilde{w} \cdot \big( \vecc{\widetilde D}\,\tq\big)\,\dd\vec{x} + \int_{\partial T} \tilde{w}\,\vecc{\widetilde D}\,\tq\cdot\vec{\nu}_T \,\dd\sigma
\;=\; \int_{T} \tilde{w}\,\tf\,\dd\vec{x}\;,
\\
&\int_{T} \vec{\tilde{y}}\cdot\tq\,\dd\vec{x}  - \int_{T} \div\vec{\tilde{y}} \, \thh \,\dd\vec{x} 
+ \int_{\partial T} \vec{\tilde{y}}\cdot\vec{\nu}_T\,\thh\,\dd\sigma
\;=\;0 \;.
\end{align*}
To improve readability, the space and time arguments are omitted whenever no ambiguity is possible.

\subsection{Variational formulation of the free-flow system}
\label{sec:variational-freeflow}

In the elevation equation~\eqref{eq:free:nonmixed:surf} and momentum equations~\eqref{eq:free:nonmixed:momentum}, we combine the advective fluxes in the \emph{primitive numerical fluxes}
\begin{equation*}
C_h(h, u^1) \coloneqq \int_{\zb}^\xi u^1 \, \dd x^2
\qquad\text{and}\qquad
\vec{C}_u(h, \vec{u}) \coloneqq u^1 \vec{u} + \begin{bmatrix} gh \\ 0 \end{bmatrix},
\end{equation*}
and, as in the subsurface problem, write our system in mixed form introducing an~auxiliary unknown~$\vec{q} \coloneqq - \grad u^1$:
\begin{subequations}\label{eq:free:mixed}
 \begin{align}
  \partial_t h + \partial_{x^1} C_h(h,u^1) &= 0 &&\text{in}~J \times \Proj \Omega \,, \label{eq:free:mixed:surf}\\
  \partial_t u^1 + \div \left( \vec{C}_u(h, \vec{u}) + \vecc{D} \, \vec{q} \right) &= f - g \, \partial_{x^1} \zb &&\text{in}~J\times\Omega \,, \label{eq:free:mixed:u}\\
  \vec{q} + \grad u^1 &= 0 &&\text{in}~J\times\Omega \,, \label{eq:free:mixed:q} \\
  \div \vec{u} &= 0 &&\text{in}~J\times\Omega \,, \label{eq:free:mixed:cont}\\
  u^1 &\;=\; u^1_\mathrm{D} &&\text{on}~J\times ( \partial\Omega_\mathrm{land} \cup \partial\Omega_\mathrm{riv} \cup \partial\Omega_\mathrm{bot})\,,\\
  u^2 &\;=\; u^2_\mathrm{D} &&\text{on}~J\times\partial\Omega_\mathrm{bot}\,,\\
  \vecc{D} \vec{q} \cdot \vec{\nu} &\;=\; q_\mathrm{D} &&\text{on}~J\times ( \partial\Omega_\mathrm{top} \cup \partial\Omega_\mathrm{os} \cup \partial\Omega_\mathrm{rad} )\,,\\
  h &\;=\; h_\mathrm{D} &&\text{on}~J\times\Proj ( \partial\Omega_\mathrm{os} \cup \partial\Omega_\mathrm{riv} )\,,\\
  u^1(0,\vec{x}) &\;=\; u^1_0 (\vec{x}) &&\text{in}~\Omega\,,\\
  h(0,x^1) &\;=\; h_0 (x^1) &&\text{in}~\Proj \Omega\,.
 \end{align}
\end{subequations}

Similarly to Sec.~\ref{sec:variational-sub}, we formulate the variational system on an element-by-element basis, multiply both sides of equations~\eqref{eq:free:mixed:u}--\eqref{eq:free:mixed:cont} by smooth test functions $z: T \rightarrow \IR$, $\vec{y}: T \rightarrow \IR^2$, and $w: T \rightarrow \IR$, correspondingly, and integrate by parts over $T \in \setT_{\Delta}$ to obtain
\begin{align*}
&\int_T z \, \partial_t u^1 \dd\vec{x}
- \int_T \grad z \cdot \Big( \vec{C}_u(h, \vec{u}) + \vecc{D} \, \vec{q} \Big) \,\dd \vec{x} 
+ \int_{\partial T} z \, \left(\vec{C}_u(h, \vec{u}) + \vecc{D} \, \vec{q} \right) 
  \cdot \vec{\nu}_T \, \dd\sigma
= \int_T z \, \left( f - g \, \partial_{x^1} \zb \right) \,\dd \vec{x} \;, \\
&\int_T \vec{y} \cdot \vec{q} \, \dd\vec{x} - \int_T \div \vec{y} \, u^1 \dd\vec{x} 
+ \int_{\partial T} \vec{y} \cdot \vec{\nu}_T \, u^1 \, \dd\sigma
= 0 \;, \\
&-\int_T \grad w \cdot \vec{u} \, \dd\vec{x} 
+ \int_{\partial T} w \, \vec{u} \cdot \vec{\nu}_T \, \dd\sigma
= 0 \;.
\end{align*}
For Eq.~\eqref{eq:free:mixed:surf}, we denote by $\overline{\elem} \coloneqq \Proj \elem$ the one-dimensional element corresponding to~$T$ (cf.~Sec.~\ref{sec:model:ff}), multiply both sides by a~smooth test function $\overline{w}: \Proj \Omega \rightarrow \IR$, and integrate by parts.
\begin{equation*}
 \int_{\overline{T}} \overline{w} \, \partial_t h \, \dd x^1 - \int_{\overline{T}} \left( \partial_{x^1} \overline{w} \right)  C_h(h,u^1) \, \dd x^1 + \int_{\partial\overline{T}} \overline{w} \, C_h(h,u^1) \nu \, \dd \sigma 
= 0 \;.
\end{equation*}
To keep the notation uniform, integrals over zero-dimensional domains in~$\partial\overline{T}$ denote evaluation at the respective point.

The computation of the depth-integrated velocity used in $C_h(h,u^1)$ can be performed by summing the integrals over elements~$T_k$ for which $\Proj T_k = \overline{T}$ holds~\cite{DawsonAizinger2005,AizingerPDPN2013}.
Then the depth-integrated horizontal velocity on~$\overline{T}$ is given by
\begin{equation}
\label{eq:depth-int-vel}
\overline{u}^1(x^1) \big|_{\overline{T}} \coloneqq 
\int_{\zb}^\xi u^1(x^1) \, \dd x^2 =
\sum_{k=1}^L \int_{\zeta_{k-1}}^{\zeta_k} u^1(x^1) \big|_{T_k} \dd x^2\,,
\end{equation}
with $\zb = \zeta_0 < \zeta_1 < \dots < \zeta_L = \xi$ and $L$ being the number of elements $T$ for which $\Proj T = \overline{T}$.
For brevity and readability we drop the \enquote{$\big|_T$} in the following.

\subsection{Definitions and preliminaries}

Before describing the DG~scheme for~\eqref{eq:sub},~\eqref{eq:free:mixed},~\eqref{eq:depth-int-vel}, we introduce some notation.
Let~$\Omega \in \{\Omega(t), \tO\}$, then $\setE_\Omega$ denotes the set of interior edges, $\setE_{\partial\Omega}$ the set of boundary edges, and $\setE \coloneqq \setE_\Omega\cup\setE_{\partial\Omega}$ the set of all edges $E$. 
For an interior edge $E\in\setE_\Omega$ shared by elements $T^-$ and $T^+$, we define the one-sided values of a~scalar quantity~$w=w(\vec{x})$ on~$E$ by
\begin{equation*}
w^-(\vec{x})\;\coloneqq\;\lim_{\varepsilon \to 0^{+}} w(\vec{x} - \varepsilon\,\vec{\nu}_{T^-})
\qquad\text{and}\qquad
w^+(\vec{x})\;\coloneqq\;\lim_{\varepsilon \to 0^{+}} w(\vec{x} - \varepsilon\,\vec{\nu}_{T^+})\;,
\end{equation*}
respectively.  
For a~boundary edge~$E\in\setE_{\partial\Omega}$, only the definition on the left is meaningful. 
The \emph{average} and the \emph{jump} of~$w$ on~$E \in \setE_{\partial\Omega}$ are then given by
\begin{equation*}
\avg{w} \;\coloneqq\; (w^- + w^+)/2\qquad  \mbox{and} \qquad \jump{w} \;\coloneqq\; w^- \vec{\nu}_{T^-} + w^+ \vec{\nu}_{T^+} \;=\; (w^- - w^+)\,\vec{\nu}_{T^-}\,,
\end{equation*}
respectively. 
Note that $\jump{w}$ is a~vector-valued quantity.
Further notation is introduced on the first use and is summarized in the \emph{Index of notation}.
\par
For $\IQ_p(T)$, the standard tensor-product polynomial space of degree~$p$ in each variable on~$T\in\setT_{\Delta}$, we define the \emph{broken polynomial space} by
\begin{equation*}
\IQ_p(\mathcal{T}_{\Delta}) \coloneqq \{ w_{\Delta} \in L^2(\Omega) \, : \, \forall \elem \in \mathcal{T}_{\Delta},~w_{\Delta} |_\elem \in \IQ_p(\elem)\}\,.
\end{equation*}
Note that we could also use~$\IP_p(T)$ instead of~$\IQ_p(T)$ with some changes to the computation of the depth-integrated velocity in the free flow problem~\cite{Aizinger2004}.
Moreover, the \emph{$L^2$-projection}~~$\pi : L^2(\mathcal{T}_{\Delta}) \to \IQ_p(\mathcal{T}_{\Delta})$ is defined as
\begin{equation*}
\forall \varphi \in \IQ_p(\mathcal{T}_{\Delta}) ,\quad \int_\Omega \big(\pi (w) - w\big) \, \varphi \; \dd\vec x = 0~.
\end{equation*}
The corresponding component-wise generalization for vector-valued functions is denoted by the same symbol.

\subsection{Local basis representation and transformation rules}
\label{sec:gridbasis}
\begin{figure}[t]
\begin{center}
\includegraphics{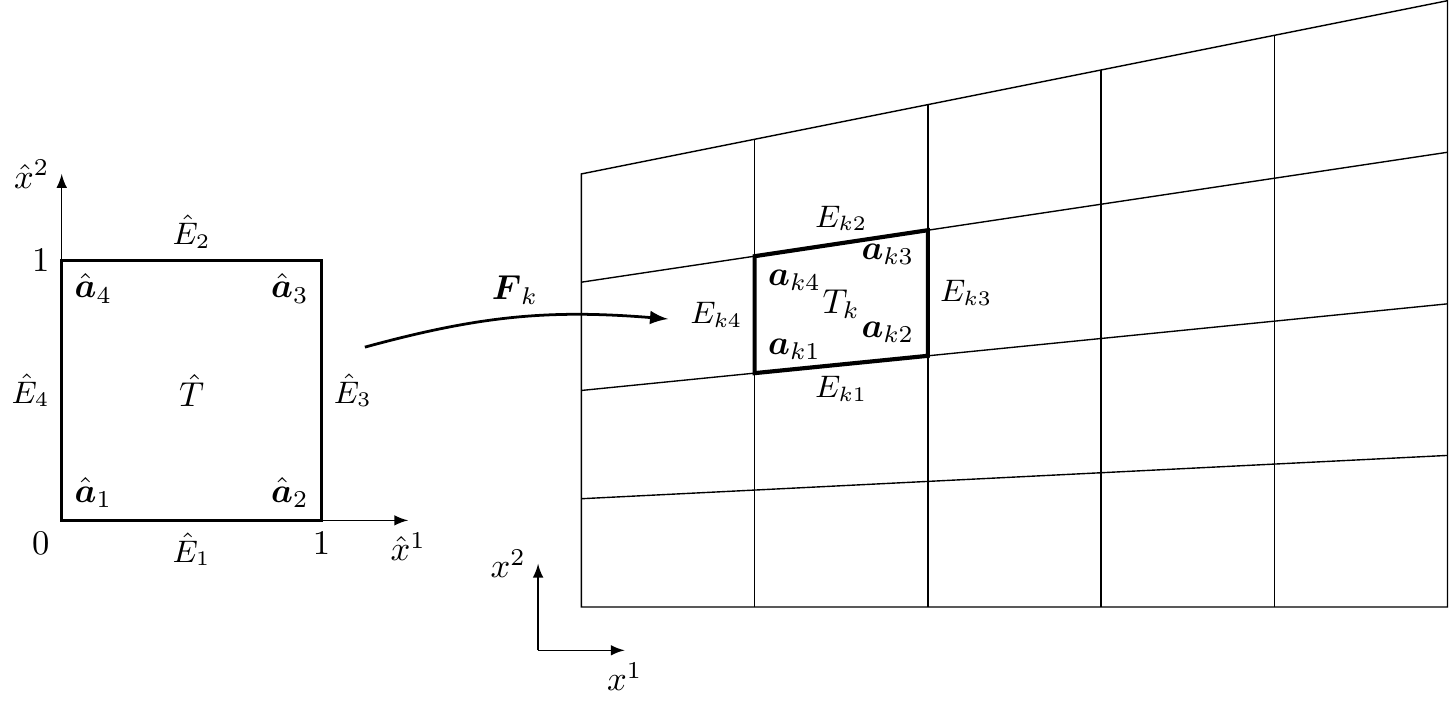}
\end{center}
\caption{Trapezoidal mesh~$\setT_{\Delta}$ with parallel vertical element edges and transformation $\vec{F}_k$ from the~reference square~$\hat{T}$}
\label{fig:setT}
\end{figure}

In the following, we denote by~$T_k$ an element of~$\setT_{\Delta}$ and use a~local numbering scheme to identify its vertices~$\vec{a}_{ki}$ and edges~$E_{kn}$, $i,n \in \{1,2,3,4\}$.
Consequently, an~interior edge belonging to neighboring elements~$T_{k^-}, T_{k^+}$ is identified by $E_{k^-n^-} = E_{k^+n^+}$.
Due to the structure of the mesh (cf.~Fig.~\ref{fig:setT}), the local edge index~$n^+$ is directly deduced from~$n^-$ via
\begin{equation*}
  n^+ = \left\{\begin{matrix}
    2 & \text{if }\;n^- = 1 \\
    1 & \text{if }\;n^- = 2 \\
    4 & \text{if }\;n^- = 3 \\
    3 & \text{if }\;n^- = 4
  \end{matrix}\right\}\,.
\end{equation*}

As in our previous works~\cite{FESTUNG1,FESTUNG2,FESTUNG3}, we employ a~mixture of algebraic and numerical indexing styles: 
for instance, $E_{kn}\in\partial T_k\cap\setE_\Omega$ means \emph{all possible} combinations of a~fixed element index~$k\in\{1,\ldots,K\}$ with local edge indices $n\in\{1,2,3,4\}$ such that $E_{kn}$ lies in~$\partial T_k\cap\setE_\Omega$. 
This implicitly fixes the numerical indices which accordingly can be used to index matrices or arrays.
At some points of the discretization of the free-flow problem, we require a~distinction between lateral edges on the one hand and  top and bottom edges on the other.
For simplicity, we refer to them as `vertical' and `horizontal' edges and mark all sets of edges correspondingly using superscripts `v' or `h' although the latter are not necessarily orthogonal to the direction of gravity (cf.~Fig.~\ref{fig:setT}).
For example, $\setE_\Omega^\mathrm{v}$ and $\setE_\Omega^\mathrm{h}$ are the sets of all vertical and horizontal interior edges, respectively.
The numbering scheme of the mesh ensures that we have local edge index $n\in\{1,2\}$ for horizontal and $n\in\{3,4\}$ for vertical edges.

Following the standard practice, we define our discrete solution $c_{\Delta}(\vec{x}) = \sum_{j = 1}^N C_j \, \vphi_j(\vec{x})$ expressed in some finite basis ~$\{\vphi_i \}_{i=1,\ldots,N}$ of $\IQ_p(\setT_{\Delta})$ using a~reference element. 
For this, we use the \emph{unit reference square}~$\hat{T}$ as shown in Fig.~\ref{fig:setT} and specify for~$T_k\in\setT_{\Delta}$ a~$C^1$-diffeomorphism
\begin{equation*}
\vec{F}_k :\;  \hat{T} \ni \hat{\vec{x}}\mapsto \vec{x}\in T_k\,.
\end{equation*}
We choose a~set of orthonormal with respect to the $L^2$-inner product basis functions~$\{ \hat{\vphi}_j \}_j$ for~$\IQ_p(\hat{\elem})$, define the basis of~$\IQ_p(T_k)$ as~$\hat{\vphi}_j \circ \vec{F}_k^{-1}$, and obtain representation
\begin{equation*}
c_{\Delta}\big|_{\elem_k} = \sum_{j=1}^N C_{kj} \left( \hat{\vphi}_j \circ \vec{F}_k^{-1} \right)\,.
\end{equation*}
With~$\hat{T}$ explicitly defined, the mapping can be expressed in terms of the vertices~$\vec{a}_{ki} = \transpose{[a_{ki}^1,a_{ki}^2]}, i\in\{1,\ldots,4\}$ of~$T_k$:
\begin{subequations}\label{eq:mappings}
\begin{equation} \label{eq:mappings:Fk}
\vec{F}_k(\hat{\vec{x}})\;\coloneqq\; 
\vec{a}_{k1} + (\vec{a}_{k2} - \vec{a}_{k1})\,\hat{x}^1 + \begin{bmatrix}
0 \\
(a_{k4}^2 - a_{k1}^2) + \left((a_{k3}^2-a_{k2}^2) - (a_{k4}^2-a_{k1}^2)\right) \hat{x}^1
\end{bmatrix} \,  \hat{x}^2 \,.
\end{equation}
The Jacobian of the mapping~$\vecc{J}_k(\hat{\vec{x}}) \coloneqq \hat{\grad}\vec{F}_k(\hat{\vec{x}})$ and its determinant are given as
\begin{align} 
\vecc{J}_k(\hat{\vec{x}}) 
&= \begin{bmatrix}
a_{k2}^1 - a_{k1}^1 & 0 \\
a_{k2}^2 - a_{k1}^2 + \left((a_{k3}^2-a_{k2}^2) - (a_{k4}^2-a_{k1}^2)\right) \hat{x}^2 & a_{k4}^2 - a_{k1}^2 + \left((a_{k3}^2-a_{k2}^2) - (a_{k4}^2-a_{k1}^2)\right) \hat{x}^1
\end{bmatrix} \nonumber\\
&= \begin{bmatrix}
a_{k2}^1 - a_{k1}^1 & 0 \\
a_{k2}^2 - a_{k1}^2 & a_{k4}^2 - a_{k1}^2
\end{bmatrix} + \begin{bmatrix}
0 & 0 \\
0 & (a_{k3}^2-a_{k2}^2) - (a_{k4}^2-a_{k1}^2)
\end{bmatrix} \hat{x}^1 \label{eq:mappings:Jacobian}\\
&\qquad+ \begin{bmatrix}
0 & 0 \\
(a_{k3}^2-a_{k2}^2) - (a_{k4}^2-a_{k1}^2) & 0
\end{bmatrix} \hat{x}^2 
\;\eqqcolon\; \vecc{J}_k^1 + \vecc{J}_k^2 \hat{x}^1 + \vecc{J}_k^3 \hat{x}^2 \;,\nonumber\\
\det\left(\vecc{J}_k\right)(\hat{\vec{x}}) 
&= (a_{k2}^1 - a_{k1}^1) \left[ a_{k4}^2-a_{k1}^2 + \left((a_{k3}^2-a_{k2}^2) - (a_{k4}^2-a_{k1}^2)\right) \hat{x}^1 \right] \nonumber\\
&= (a_{k2}^1 - a_{k1}^1) (a_{k4}^2-a_{k1}^2) + (a_{k2}^1 - a_{k1}^1) \left((a_{k3}^2-a_{k2}^2) - (a_{k4}^2-a_{k1}^2)\right) \hat{x}^1 \label{eq:mappings:detJacobian} \\
&\eqqcolon J_k^1 + J_k^2 \hat{x}^1 \;.\nonumber
\end{align}
\end{subequations}
Due to generally non-parallel top and bottom edges of the trapezoidal elements, the entries of the Jacobian matrix $\vecc{J}_k(\hat{\vec{x}})$ are not constants but rather affine-linear functions of~$\hat{\vec{x}}$.
The local node numbering in our mesh preserves the orientation of the reference element, thus $\forall k \in \{1,\dots,K\}$, $\forall \hat{\vec{x}} \in \hat{T}:\det{\vecc{J}_k}(\hat{\vec{x}}) > 0$ holds.
Using the Jacobian we obtain the component-wise rule for the transformation of the gradient for~$m\in\{1,2\}$:
{\allowdisplaybreaks
\begin{align}\nonumber
\partial_{x^m} c(\vec{x}) \;=\; &
\Big[ \vecc{J}_k^{-T}(\hat{\vec{x}}) \Big]_{m,:} \hat{\grad} \hat{c}(\hat{\vec{x}}) \nonumber \\
\;=\; &
\frac{1}{\det(\vecc{J}_k)(\hat{\vec{x}})} \left( \Big[ \vecc{J}_k(\hat{\vec{x}}) \Big]_{3-m,3-m} \partial_{\hat{x}^m} \hat{c}(\hat{\vec{x}}) - \Big[ \vecc{J}_k(\hat{\vec{x}}) \Big]_{3-m,m} \partial_{\hat{x}^{3-m}} \hat{c}(\hat{\vec{x}}) \right) \label{eq:rule:gradient} \\
\;=\; &
\frac{1}{\det(\vecc{J}_k)(\hat{\vec{x}})} \left( 
\Big[ \vecc{J}_k^1 \Big]_{3-m,3-m} \partial_{\hat{x}^m} \hat{c}(\hat{\vec{x}}) - 
\Big[ \vecc{J}_k^1 \Big]_{3-m,m} \partial_{\hat{x}^{3-m}} \hat{c}(\hat{\vec{x}}) + 
\Big[ \vecc{J}_k^2 \Big]_{3-m,3-m} \hat{x}^1 \partial_{\hat{x}^m} \hat{c}(\hat{\vec{x}})  \right. \nonumber\\
&\left. - 
\Big[ \vecc{J}_k^2 \Big]_{3-m,m} \hat{x}^1 \partial_{\hat{x}^{3-m}} \hat{c}(\hat{\vec{x}}) + 
\Big[ \vecc{J}_k^3 \Big]_{3-m,3-m} \hat{x}^2 \partial_{\hat{x}^m} \hat{c}(\hat{\vec{x}}) - 
\Big[ \vecc{J}_k^3 \Big]_{3-m,m} \hat{x}^2 \partial_{\hat{x}^{3-m}} \hat{c}(\hat{\vec{x}}) \right)\,. \nonumber
\end{align}
}

For our choice of~$\IQ_p(T)$ as the space of ansatz functions, we define basis functions~$\hat{\vphi}_i$ on the reference square~$\hat{\elem}$ via tensor products of one-dimensional Legendre polynomials~$\hat{\phi}_m: [0,1] \rightarrow \IR$ as
\begin{equation}\label{eq:basis:tensor-prod}
\hat{\vphi}_i\left(\hat{x}^1, \hat{x}^2\right) = \hat{\phi}_m\left(\hat{x}^1\right)\, \hat{\phi}_n\left(\hat{x}^2\right) \,,
\qquad i \in \{1, \dots, N\} 
\;\text{ and }\; m,n \in \{1, \dots, \overline{N}\} 
\;\text{ with }\; N = \overline{N}^2,
\end{equation}
where $\overline{N} = p+1$ and $i \coloneqq i(m,n) = \left(\max(m,n) - 1\right)^2 + \max(m,n) - m + n$ (see~Table~\ref{tab:basis-index}).
Closed-form expressions for the one-dimensional basis functions on the reference interval $[0,1]$ up to order three are given by:
\begin{align*}
&\hat{\phi}_1(\hat{x})=1\,,\qquad
\hat{\phi}_2(\hat{x})= \sqrt{3}\,(1 - 2 \hat{x})\,,\qquad
\hat{\phi}_3(\hat{x})= \sqrt{5}\,\big((6\hat{x}-6)\hat{x}+1\big)\,,\\
&\hat{\phi}_4(\hat{x})= \sqrt{7}\,\Big(\big((20\hat{x}-30)\hat{x}+12\big)\hat{x}-1\Big)\,.
\end{align*}

\begin{table}[t]
\setlength{\tabcolsep}{0.68mm}
\centering
\begin{tabularx}{\linewidth}{@{}c|c|ccc|ccccc|ccccccc@{}}
\toprule
$p$     & 0       & \multicolumn{3}{c|}{1}      & \multicolumn{5}{c|}{2} & \multicolumn{7}{c}{3}\\
$[m,n]$ & $[1,1]$ & $[2,1]$ & $[2,2]$ & $[1,2]$ & $[3,1]$ & $[3,2]$ & $[3,3]$ & $[2,3]$ & $[1,3]$ & $[4,1]$ & $[4,2]$ & $[4,3]$ & $[4,4]$ & $[3,4]$ & $[2,4]$ & $[1,4]$ \\
$i$     & 1       & 2       & 3       & 4       & 5       & 6       & 7       & 8    & 9   & 10      & 11      & 12      & 13      & 14      & 15      & 16\\
\bottomrule
\end{tabularx}
\caption{Correlation between indices of one-dimensional ($\hat{\phi}_m$, $\hat{\phi}_n$) and two-dimensional basis-functions ($\hat{\vphi}_i$) for polynomial orders $p=0, 1, 2, 3$.}
\label{tab:basis-index}
\end{table}

Any function~$\,c:T_k\rightarrow \IR\,$ implies $\hat{c}:\hat{T}\rightarrow \IR\,$ by $\,\hat{c}=c\circ \vec{F}_k\,$, i.\,e., $\,c(\vec{x}) = \hat{c}(\hat{\vec{x}})\,$, and, in particular, $\vphi_{ki}(\vec{x}) = \hat{\vphi}_{i}(\hat{\vec{x}})$ for all $k$.
For integrals over element~$\elem_k$ and edge~$E_{kn}\subset \elem_k$, we use transformation formulas
\begin{subequations}
\begin{align}
&\int_{T_k}c(\vec{x})\,\dd\vec{x} \;=\; \int_{\hat{T}} ( c\circ\vec{F}_k ) (\hat{\vec{x}})\,\abs{\det{\vecc{J}_k}(\hat{\vec{x}})}\,\dd\hat{\vec{x}}\;=\; \int_{\hat{T}}\hat{c}(\hat{\vec{x}})\,\abs{\det{\vecc{J}_k}(\hat{\vec{x}})}\,\dd\hat{\vec{x}}
\;,\label{eq:trafoRule:T} \\
&\int_{E_{kn}} c(\vec{x})\,\dd\sigma 
\;=\; \frac{\abs{E_{kn}}}{\abs{\hat{E}_n}} \int_{\hat{E}_n} ( c\circ\vec{F}_k ) (\hat{\vec{x}})\,\dd\hat{\sigma}\;=\; \abs{E_{kn}} \int_{\hat{E}_n} \hat{c}(\hat{\vec{x}})\,\dd\hat{\sigma}
\label{eq:trafoRule:E}
\end{align}
and, for integrals over one-dimensional elements~$\overline{T}_{\overline{k}} = [a_{\overline{k}1}^1, a_{\overline{k}2}^1]$, we have
\begin{equation}\label{eq:trafoRule:barT}
\int_{\overline{T}_{\overline{k}}} \overline{c}(x^1) \,\dd x^1 = |\overline{T}| \int_0^1 \hat{\overline{c}}(\hat{x}) \,\dd \hat{x}
\end{equation}
\end{subequations}
with~$\overline{c}:\overline{T}_{\overline{k}} \rightarrow\IR$ implying~$\hat{\overline{c}}:[0,1]\rightarrow\IR$ by~$\hat{\overline{c}} = \overline{c}\circ \overline{F}_{\overline{k}}$  and using the standard linear mapping~$\overline{F}_{\overline{k}} : [0,1] \rightarrow  \overline{T}_{\overline{k}}$.

\subsection{Semi-discrete formulation for the subsurface system}
Setting $\vecc{\widetilde D}_{\Delta}(t,\cdot) \coloneqq \pi (\vecc{\widetilde D}(t,\cdot)), \, \thh_{\Delta}(0) \coloneqq \pi (\thh_0)$, we seek $\left(\tq_{\Delta}(t), \thh_{\Delta}(t)\right)\in  [\IQ_p(\setT_{\Delta})]^2\times \IQ_p(\setT_{\Delta})$ such that the following holds for a.\,e.~$t\in J$, $\forall T^-\in\setT_{\Delta},\, \forall \vec{\tilde{y}}_{\Delta}\in [\IQ_p(\setT_{\Delta})]^2,\, \forall \tilde{w}_{\Delta}\in \IQ_p(\setT_{\Delta})\,$:
\begin{subequations}\label{eq:sub:semidisc}
\begin{align}
& \int_{T^-} \vec{\tilde{y}}_{\Delta}\cdot\tq_{\Delta} \,\dd\vec{x}
\; - \int_{T^-} \div\vec{\tilde{y}}_{\Delta} \,\thh_{\Delta} \,\dd\vec{x}
\; + \int_{\partial T^-} \vec{\tilde{y}}_{\Delta}^-\cdot\vec{\nu}_{T^-}\,\left\lbrace
\begin{array}{cl}
\avg{\thh_{\Delta} }   &\text{on}~ \setE_\Omega\\
\thh_\mathrm{D} &\text{on}~ \setE_{\mathrm{D}}\\
\thh_{\Delta}^-        &\text{on}~ \setE_\mathrm{N}
\end{array}\right\rbrace \,\dd\sigma  = 0\,,
\label{eq:sub:semidisc:q}
\\
& \begin{aligned}
\int_{T^-} \tilde{w}_{\Delta}\,\partial_t \thh_{\Delta} \,\dd\vec{x}
\;- \int_{T^-}\grad \tilde{w}_{\Delta} \cdot\Big(\vecc{\widetilde D}_{\Delta}\,\tq_{\Delta} \Big) \,\dd\vec{x} \qquad\qquad\qquad\qquad\qquad&\\
\;+ \int_{\partial T^-} \tilde{w}_{\Delta}^-\,
\begin{Bmatrix}
\displaystyle \avg{\vecc{\widetilde D}_{\Delta}\,  \tq_{\Delta}}\cdot\vec{\nu}_{T^-} + \frac{\eta}{|E|}  \jump{\thh_{\Delta}}\cdot\vec{\nu}_{T^-}& \text{on}~\setE_\Omega\\
\displaystyle \vecc{\widetilde D}_{\Delta}^-\,  \tq_{\Delta}^-\cdot\vec{\nu}_{T^-}  + \frac{\eta}{|E|} \big( \thh^-_{\Delta} - \thh_\mathrm{D} \big)  & \text{on}~\setE_{\mathrm{D}}\\
\displaystyle \tgN     & \text{on}~\setE_\mathrm{N}
\end{Bmatrix} \,\dd\sigma
&\;=\; \int_{T^-} \tilde{w}_{\Delta}\,\tf_{\Delta} \,\dd\vec{x}\,,
\end{aligned}
\label{eq:sub:semidisc:u}
\end{align}
\end{subequations}
where $\eta>0$ is a~penalty coefficient, and $|E|$ denotes the Lebesgue measure (length) of the interface. 
The penalty terms in~\eqref{eq:sub:semidisc:u} are required  to ensure 
a~full rank of the system in the absence of the time derivative~\cite[Lem.~2.15]{RuppKnabner,RuppKnabnerDawson2018,Riviere2008}.  
For analysis purposes, the above equations are usually summed over all elements~$T\in\setT_{\Delta}$.  
In the implementation that follows, however, it is more convenient to work with element-local equations. 

\subsubsection{System of equations}
%
Testing~\eqref{eq:sub:semidisc:q} with 
$\vec{\tilde{y}}_{\Delta} = \transpose{[\vphi_{ki}, 0]}, \transpose{[0, \vphi_{ki}]}$ and~\eqref{eq:sub:semidisc:u} with $\tilde{w}_{\Delta} = \vphi_{ki}$ for $i\in\{1,\ldots,N\}$ yields a~\emph{time-dependent system of equations}.
The resulting system and the terms involved have the same form as in our first paper in series, we refer to Sections~2.4.2--2.4.4 of~\cite{FESTUNG1} for a~full presentation of the discretization steps and the component-wise definitions of the block-matrices and right-hand side vectors.
The system written in matrix form is given by
\newcommand\scalemath[2]{\scalebox{#1}{\mbox{\ensuremath{\displaystyle #2}}}}
\begin{equation}\label{eq:sub:matsys}
\scalemath{0.80}{%
\underbrace{\begin{bmatrix}
\vecc{0}&\vecc{0}& \vecc{0}
\\
\vecc{0}&\vecc{0}&\vecc{0}
\\
\vecc{0}&\vecc{0}&\tilde{\vecc{M}}
\end{bmatrix}}_{\eqqcolon \;\tilde{\vecc{W}}}
\begin{bmatrix}
\partial_t\vec{\tQ}^1\\ \partial_t \vec{\tQ}^2\\ \partial_t\vec{\tHH}
\end{bmatrix}
+
\underbrace{\begin{bmatrix}
\tilde{\vecc{M}}& \cdot  & -\tilde{\vecc{H}}^1{+}\tilde{\vecc{Q}}^1{+}\tilde{\vecc{Q}}_{\mathrm{N}}^1\\
\cdot & \tilde{\vecc{M}} & -\tilde{\vecc{H}}^2{+}\tilde{\vecc{Q}}^2{+}\tilde{\vecc{Q}}_{\mathrm{N}}^2\\
~-\tilde{\vecc{G}}^1{+}\tilde{\vecc{R}}^1 {+}\tilde{\vecc{R}}^1_\mathrm{D}
& -\tilde{\vecc{G}}^2{+}\tilde{\vecc{R}}^2 {+}\tilde{\vecc{R}}^2_\mathrm{D}
& \eta\,\big( \tilde{\vecc{S}} {+} \tilde{\vecc{S}}_\mathrm{D}\big)
\end{bmatrix}}_{\eqqcolon \;\tilde{\vecc{A}}(t)}
\,
\underbrace{\begin{bmatrix}
\vec{\tQ}^1\\\vec{\tQ}^2\\\vec{\tHH}
\end{bmatrix}}_{\eqqcolon \;\tilde{\vec{Y}}(t)}
~~=~~ 
\underbrace{\begin{bmatrix}
-\tilde{\vec{J}}_{\mathrm{D}}^1\\
-\tilde{\vec{J}}_{\mathrm{D}}^2\\
\eta\, \tilde{\vec{K}}_{\mathrm{D}} - \tilde{\vec{K}}_{\mathrm{N}} {+}  \tilde{\vec{L}}
\end{bmatrix}}_{\eqqcolon \;\tilde{\vec{V}}(t)}
}
\end{equation}
with the representation vectors
\begin{align*}
\vec{\tQ}^m(t) &\;\coloneqq\;
\transpose{\begin{bmatrix}
\tQ_{11}^m(t) & \cdots & \tQ_{1N}^m(t) & \cdots &\cdots & \tQ_{K1}^m(t) & \cdots & \tQ_{KN}^m(t)
\end{bmatrix}}\quad \text{for}~m\in\{1,2\}\;,
\\
\vec{\tHH}(t) &\;\coloneqq\;
\transpose{\begin{bmatrix}
\tHH_{11}(t) & \cdots & \tHH_{1N}(t) & \cdots &\cdots & \tHH_{K1}(t) & \cdots & \tHH_{KN}(t)
\end{bmatrix}}\;.
\end{align*}
Note that blocks~$\tilde{\vecc{G}}^m,\,\tilde{\vecc{R}}^m,\,\tilde{\vecc{R}}^m_\mathrm{D}$ (which contain~$\vecc{\widetilde D}$) and vectors~$\tilde{\vec{J}}^m_\mathrm{D},\,\tilde{\vec{K}}_\mathrm{D},\,\tilde{\vec{K}}_\mathrm{N},\,\tilde{\vec{L}}$ (including boundary data and right-hand side) are time-dependent (time~arguments are suppressed here).

\subsubsection{Time discretization}

System~\eqref{eq:sub:matsys} is equivalent to
\begin{equation} \label{eq:sub:fullsys}
\tilde{\vecc{W}}\partial_t\tilde{\vec{Y}}(t) + \tilde{\vecc{A}}(t)\,\tilde{\vec{Y}}(t)  \;=\; \tilde{\vec{V}}(t)
\end{equation}
with solution $\tilde{\vec{Y}}(t) \in \IR^{3KN}$, right-hand-side vector~$\tilde{\vec{V}}(t) \in \IR^{3KN}$, and matrices~$\tilde{\vecc{A}}(t)$, $\tilde{\vecc{W}}\in\IR^{3KN\times3KN}$ all given in~\eqref{eq:sub:matsys}.
We discretize system~\eqref{eq:sub:fullsys} in time using the implicit Euler method.
Let $0= t^1<t^2<\ldots\leq t_\mathrm{end}$ be a~not necessarily equidistant decomposition of the time interval~$J$, and let $\Delta t^n \coloneqq t^{n+1} -t^n$ denote the time step size. Then
\begin{equation*} 
\left(\tilde{\vecc{W}} + \Delta t^n\,\tilde{\vecc{A}}^{n+1} \right)\,\tilde{\vec{Y}}^{n+1} 
\;=\; \tilde{\vecc{W}}\,\tilde{\vec{Y}}^{n} + \Delta t^n \, \tilde{\vec{V}}^{n+1} 
\end{equation*}
gives the prescription for a~time step, where we abbreviated~$\tilde{\vecc{A}}^n\coloneqq \tilde{\vecc{A}}(t^n)$, etc.

\subsection{Semi-discrete formulation for the free-flow system}
\label{sec:free:semi-disc}

First, we introduce some additional notation: 
As mentioned in Sec.~\ref{sec:gridbasis}, we divide the sets of edges~$\setE_*$ into sets of horizontal and vertical edges denoted by $\setE_*^\mathrm{h}$, $\setE_*^\mathrm{v}$, respectively. Into $\setE_\IH$, $\setE_\IU$, $\setE_\IQ$, we collect boundary edges with prescribed Dirichlet data~$h_\mathrm{D}$, $u^1_\mathrm{D}$, or~$q_\mathrm{D}$ and into $\setE_\mathrm{bot}$ the edges on the bottom boundary of the free-flow domain. 

In order to obtain a semi-discrete formulation we must take into account the discontinuous nature of our broken polynomial ansatz space. 
For integrals on vertical boundaries, we introduce approximations~$\widehat{C}_u(h_{\Delta}^-, \vec{u}_{\Delta}^-, h_{\Delta}^+, \vec{u}_{\Delta}^+)$ and $\widehat{C}_h(h_\mathrm{s}, h_{\Delta}^-, {u_{\Delta}^1}^-, h_{\Delta}^+, {u_{\Delta}^1}^+)$ to the nonlinear boundary fluxes $\vec{C}_u(h_{\Delta}, \vec{u}_{\Delta}) \cdot \vec{\nu}$ and $\vec{u} \cdot \vec{\nu}$, respectively, and use the central numerical fluxes in integrals over horizontal boundaries. 
These fluxes are further detailed in Sec.~\ref{sec:riemann}. 
On external domain boundaries, the fluxes utilize values from the interior for unknowns not specified in the boundary conditions. 

For $\vecc{D}_{\Delta}(t,\cdot) \coloneqq \pi (\vecc{D}(t,\cdot)), \, h_{\Delta}(0) \coloneqq \pi (h_0), \, u^1_{\Delta}(0) \coloneqq \pi (u^1_0)$, we seek $(h_{\Delta}, \vec{u}_{\Delta}, \vec{q}_{\Delta}) \in \IQ_p(\Proj\setT_{\Delta}) \times [\IQ_p(\setT_{\Delta})]^2 \times [\IQ_p(\setT_{\Delta})]^2$ such that~\eqref{eq:free:semidisc} holds for a.\,e.~$t \in J$, $\forall T^- \in \setT_{\Delta}$, $\forall z_{\Delta} \in \IQ_p(\setT_{\Delta})$, $\forall \vec{y}_{\Delta} \in [\IQ_p(\setT_{\Delta})]^2$, $\forall w_{\Delta} \in \IQ_p(\setT_{\Delta})$, $\forall \overline{w}_{\Delta} \in \IQ_p(\Proj\setT_{\Delta})$:
{\allowdisplaybreaks
\begin{subequations}\label{eq:free:semidisc}
\begin{align}
& \int_{T^-} z_{\Delta} \, \partial_t u_{\Delta}^1 \, \dd\vec{x} 
- \int_{T^-} \grad z_{\Delta} \cdot \Big( \vec{C}_u(h_{\Delta}, \vec{u}_{\Delta}) + \vecc{D}_{\Delta} \, \vec{q}_{\Delta} \Big) \,\dd \vec{x} 
\nonumber \\
& \qquad+ \int_{\partial T^-} z_{\Delta}^- \, \begin{Bmatrix}
  \displaystyle \avg{\vecc{D}_{\Delta} \vec{q}_{\Delta}} \cdot \vec{\nu}_{T^-} & \mathrm{on}~\setE_\Omega \\
  q_\mathrm{D} & \mathrm{on}~\setE_\IQ\\
  \displaystyle \vecc{D}_{\Delta} \vec{q}_{\Delta}^- \cdot \vec{\nu}_{T^-} & \text{otherwise}
 \end{Bmatrix}
 \, \dd\sigma 
\nonumber \\
& \qquad+ \int_{\partial T^-} z_{\Delta}^- \, \begin{Bmatrix}
  \displaystyle \avg{ \vec{C}_u(h_{\Delta}, \vec{u}_{\Delta}) } \cdot \vec{\nu}_{T^-} & \mathrm{on}~\setE_\Omega^\mathrm{h} \\
  \widehat{C}_u(h_{\Delta}^-, \vec{u}_{\Delta}^-, h_{\Delta}^+, \vec{u}_{\Delta}^+) & \mathrm{on}~\setE_\Omega^\mathrm{v}\\
  \displaystyle \vec{C}_u(h_\mathrm{bdr}, \vec{u}_\mathrm{bdr}) \cdot \vec{\nu}_{T^-} & \mathrm{on}~\setE_{\partial\Omega}^\mathrm{h} \\
  \widehat{C}_u(h_{\Delta}^-, \vec{u}_{\Delta}^-, h_{\mathrm{bdr}}, \vec{u}_{\mathrm{bdr}}) & \mathrm{on}~\setE_{\partial\Omega}^\mathrm{v}\\
 \end{Bmatrix}\, \dd\sigma
\;=\; \int_{T^-} z_{\Delta} \big( f_{\Delta} - g \, \partial_{x^1} \zb \big) \,\dd \vec{x} \;,
\label{eq:free:semidisc:u} \\
& \int_{T^-} \vec{y}_{\Delta} \cdot \vec{q}_{\Delta} \, \dd\vec{x} 
- \int_{T^-} \div \vec{y}_{\Delta} \, u_{\Delta}^1 \, \dd\vec{x}
+ \int_{\partial T^-} \vec{y}_{\Delta}^- \cdot \vec{\nu}_{T^-} \, \begin{Bmatrix}
  \displaystyle \avg{u_{\Delta}^1} & \mathrm{on}~\setE_\Omega \\
  \displaystyle u^1_{\mathrm{bdr}} & \mathrm{on}~\setE_{\partial\Omega}
\end{Bmatrix} \, \dd \sigma \;=\; 0 \;, 
\label{eq:free:semidisc:q} \\
& - \int_{T^-} \grad w_{\Delta} \cdot \vec{u}_{\Delta} \, \dd\vec{x} 
+ \int_{\partial T^-} w_{\Delta}^- \, \begin{Bmatrix}
  \displaystyle \avg{ u_{\Delta}^1 } \nu_{T^-}^1 + u^2_\uparrow \nu_{T^-}^2 & \mathrm{on}~\setE_\Omega^\mathrm{h} \\
  \displaystyle \widehat{C}_h(\Hs, h_{\Delta}^-, {u_{\Delta}^1}^-, h_{\Delta}^+, {u_{\Delta}^1}^+) & \mathrm{on}~\setE_\Omega^\mathrm{v} \\
  \displaystyle u^1_{\mathrm{bdr}} \nu^1_{T^-} + u^2_\uparrow \nu_{T^-}^2 & \mathrm{on}~\setE_{\partial\Omega}^\mathrm{h}  \\
  \displaystyle \widehat{C}_h(\Hs, h_{\Delta}^-, {u_{\Delta}^1}^-, h_\mathrm{bdr}, u^1_\mathrm{bdr}) & \mathrm{on}~\setE_{\partial\Omega}^\mathrm{v}
\end{Bmatrix} \, \dd\sigma \;=\; 0 \;, 
\label{eq:free:semidisc:cont} \\
& \int_{\overline{T}^-} \overline{w}_{\Delta} \, \partial_t h_{\Delta} \, \dd x^1 
- \int_{\overline{T}^-} \partial_{x^1} \overline{w}_{\Delta} \frac{\overline{u}_{\Delta}^1 \, h_{\Delta}}{\Hs} \, \dd x^1 \nonumber \\
& \qquad + \displaystyle\sum_{\Proj\partial T^- = \partial\overline{T}^-} \int_{\partial T^-} \overline{w}_{\Delta}^- 
  \begin{Bmatrix}
    \displaystyle \widehat{C}_h(\Hs, h_{\Delta}^-, {u^1_{\Delta}}^-, h_{\Delta}^+, {u^1_{\Delta}}^+) & \mathrm{on}~\setE_\Omega^\mathrm{v} \\
    \displaystyle \widehat{C}_h(\Hs, h_{\Delta}^-, {u^1_{\Delta}}^-, h_\mathrm{bdr}, u^1_\mathrm{bdr}) & \mathrm{on}~\setE_{\partial\Omega}^\mathrm{v} 
  \end{Bmatrix} \, \dd\sigma 
  \;=\; 0 
\label{eq:free:semidisc:surf}
\end{align}
\end{subequations}
}
with boundary values $h_{\mathrm{bdr}}$, $\vec{u}_{\mathrm{bdr}} \coloneqq \transpose{[u_\mathrm{bdr}^1, u_\mathrm{bdr}^2]}$, $u^2_\uparrow$ defined as given Dirichlet data, where available, or interior values otherwise:
\begin{align*}
&h_{\mathrm{bdr}} \coloneqq \begin{Bmatrix}
  h_\mathrm{D} & \mathrm{on}~\setE_\IH \\
  h_{\Delta}^- & \text{otherwise}
\end{Bmatrix}\,,&\quad
&u^1_{\mathrm{bdr}} \coloneqq \begin{Bmatrix}
  u^1_\mathrm{D} & \mathrm{on}~\setE_\IU \\
  u^{1-}_{\Delta} & \text{otherwise}
\end{Bmatrix}\,,\quad\\
&u^2_{\mathrm{bdr}} \coloneqq \begin{Bmatrix}
  u^2_\mathrm{D} & \mathrm{on}~\setE_{\mathrm{bot}} \\ 
  u^{2-}_{\Delta} & \text{otherwise}
\end{Bmatrix}\,,\quad&
&u^2_\uparrow \coloneqq \begin{Bmatrix}
  u^2_\mathrm{D} & \mathrm{on}~\setE_{\mathrm{bot}} \\ 
  u^2_{\Delta}\big|_{T_\mathrm{below}} & \text{otherwise}
\end{Bmatrix} \,.
\end{align*}

By $u^2_\uparrow$ we denote taking the value from the element below. 
This is due to the fact that calculating $u^2$ can be interpreted as solving the ordinary differential equation $\partial_{x^1} u^1 + \partial_{x^2} u^2 = 0$ for given $u^1$ and initial condition $u^2 = u^2_D$ at the bottom of the free-flow domain.

\begin{figure}[t]
\centering
\includegraphics{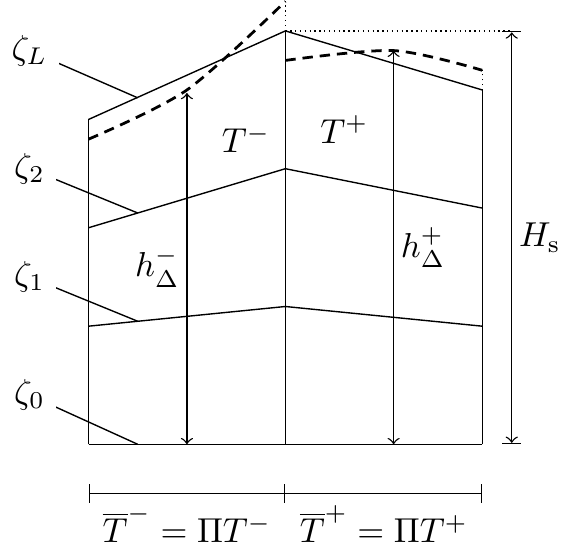}
\caption{Approximation of the free-surface geometry for discontinuous water height~$h_\Delta$ with mesh height~$\Hs$ and~$\zeta_k$, $k\in\{0,\ldots,L\}$ the $x^2$-coordinates on the horizontal element boundaries.}
\label{fig:mesh_smoothing}
\end{figure}

In free-surface equation~\eqref{eq:free:semidisc:surf}, we modified the advective flux in the following way by first combining the integrals into an integration over the lateral boundaries of the corresponding two-dimensional elements and then equivalently representing via multiplication and division by strictly positive water height~$h$ for use in the discrete flux computation detailed in Sec.~\ref{sec:riemann}:
\begin{equation*}
\int_{\partial\overline{T}} \overline{w}\,
  \left( \int_{\zb}^{\xi}u^1\dd x^2 \right) \nu_{\overline{T}} \dd x^1 
= \sum_{\Proj \partial T = \partial \overline{T}} \int_{\partial T} \overline{w}\,
  u^1 \nu^1_T \, \dd \sigma
= \sum_{\Proj \partial T = \partial \overline{T}} \int_{\partial T} \overline{w}\,
  \frac{u^1 h \,\nu^1_T}{h} \, \dd \sigma \,,
\end{equation*}
which becomes after discretization
\begin{equation*}
\sum_{\Proj \partial T^- = \partial \overline{T}^-} \int_{\partial T^-} 
  \overline{w}_{\Delta}^-\,\widehat{C}_h(h_\mathrm{s}, h_{\Delta}^-,u_{\Delta}^{1-},h_{\Delta}^+,u_{\Delta}^{1+}) \, \dd \sigma \,.
\end{equation*}
Here, $h_\mathrm{s}$ is an approximation to~$h$ and chosen to be the mesh height~$\Hs$ (see Fig.~\ref{fig:mesh_smoothing}) and, for consistency, we apply the same to the element integral.
Further details about this transformation are given in~\cite{Aizinger2004}.

For two-dimensional unknowns~$\vec{u}_{\Delta}, \vec{q}_{\Delta}$, we use the same local basis representation as for the subsurface problem (see Sec.~\ref{sec:gridbasis}), and, for the one-dimensional water height~$h_{\Delta}$, we utilize a~one-dimensional local basis representation:
\begin{equation*}
\vec{u}_{\Delta} \big|_{T_k} = \sum_{j=1}^N \begin{bmatrix} U_{kj}^1 \\ U_{kj}^2 \end{bmatrix} \, \vphi_j \,, \qquad
\vec{q}_{\Delta} \big|_{T_k} = \sum_{j=1}^N \begin{bmatrix} Q_{kj}^1 \\ Q_{kj}^2 \end{bmatrix} \, \vphi_j \,, \qquad
h_{\Delta} \big|_{\overline{T}_{\overline{k}}} = \sum_{m=1}^{\overline{N}} H_{\overline{k}m} \, \phi_m \,.
\end{equation*}
For the computation of the depth-integrated velocity (cf.~Sec.~\ref{sec:variational-freeflow}), we exploit the tensor-product structure of our two-dimensional basis functions (cf.~Eq.~\eqref{eq:basis:tensor-prod}) and the properties of Legendre polynomials to streamline the integration of $x^2$-dependent one-dimensional basis functions. We rewrite the discrete version of Eq.~\eqref{eq:depth-int-vel} as follows~\cite{ConroyKubatko2016}:
\begin{align*}
\overline{u}_{\Delta}^1(x^1) &= \sum_{k=1}^L \int_{\zeta_{k-1}}^{\zeta_k} u_{\Delta}^1(x^1,x^2)\, \dd x^2
= \sum_{k=1}^L \sum_{j=1}^N U^1_{kj} \int_{\zeta_{k-1}}^{\zeta_k} \vphi_{kj}(x^1,x^2)\, \dd x^2\\
&= \sum_{k=1}^L \sum_{j=1}^N U^1_{kj} \, \phi_{km}(x^1) \int_{\zeta_{k-1}}^{\zeta_k} \phi_{kn} (x^1,x^2) \, \dd x^2\\
&= \sum_{k=1}^L \sum_{j=1}^N U^1_{kj}\, \phi_{km}(x^1) \, \left(\zeta_k(x^1) - \zeta_{k-1}(x^1)\right) \smash{\overbrace{\int_0^1 \hat{\phi}_n(\hat{x}^2) \, \dd s}^{= \delta_{1n}}}\\
&=  \sum_{k=1}^L \sum_{m=1}^{\overline{N}} U^1_{kj(m)} \,\left(\zeta_k(x^1) - \zeta_{k-1}(x^1)\right)\, \phi_{km}(x^1)\,,
\end{align*}
with $j(m) = (m-1)^2 + 1$ since $\int_0^1 \hat{\phi}_1 \dd s = 1$. 
This means, we can represent the averaged velocity $\overline{u}_{\Delta}^1$ using a~one-dimensional local basis representation
\begin{equation}\label{eq:averagedvelocity}
\overline{u}^1_{\Delta} \big|_{\overline{T}_{\overline{k}}} = \sum_{m=1}^{\overline{N}} \overline{U}_{\overline{k}m} \, \phi_{\overline{k}m} 
\qquad\text{with}\qquad
\overline{U}_{\overline{k}m} = \sum_{k=1}^L U_{kj(m)} \,\left(\zeta_k(x^1) - \zeta_{k-1}(x^1)\right) \,.
\end{equation}
Note that~$\zeta_k(x^1) - \zeta_{k-1}(x^1)$ is the height of element~$T_k$ and thus dependent on the~$x^1$-coordinate.
For details on how we compute the depth-integrated velocity in our implementation and resolve the~$x^1$-dependency, see~\ref{sec:depth-int-vel}.

\subsubsection{System of equations}
Testing~\eqref{eq:free:semidisc:u} with $z_{\Delta} = \vphi_{ki}$, \eqref{eq:free:semidisc:q}~with $\vec{y}_{\Delta} = \transpose{[\vphi_{ki}, 0]}, \transpose{[0, \vphi_{ki}]}$, and~\eqref{eq:free:semidisc:cont} with $w_{\Delta} = \vphi_{ki}$ for $i\in\{1,\ldots,N\}$ yields a~\emph{time-dependent system of equations} whose contribution from~$T_k$ (identified with $T_{k^-}$ in boundary integrals) reads

{\allowdisplaybreaks
\begin{subequations}\label{eq:free:sys}
\begin{align}
 & \underbrace{\sum_{j=1}^N \partial_t U^1_{kj} \int_{T_k} \vphi_{ki} \, \vphi_{kj} \, \dd\vec{x}}_{\I~(\vecc{M})} - \underbrace{\sum_{m=1}^2 \sum_{j=1}^N U^m_{kj} \sum_{l=1}^N U^1_{kl} \int_{T_k} \partial_{x^m} \vphi_{ki} \, \vphi_{kl} \, \vphi_{kj} \, \dd \vec{x}}_{\II~(\vecc{E^m})} \;- \nonumber \\
&\qquad g \, \underbrace{\sum_{j=1}^{\overline{N}} H_{\overline{k}j} 
  \int_{T_k} \partial_{x^1} \vphi_{ki} \, \phi_{\overline{k}j} \, \dd\vec{x}}_{\III~(\check{\vecc{H}})}
- \underbrace{\sum_{r=1}^2 \sum_{m=1}^2 \sum_{l=1}^N D_{kl}^{rm} \sum_{j=1}^N Q^m_{kj} \int_{T_k} \partial_{x^r} \vphi_{ki} \, \vphi_{kl} \, \vphi_{kj} \,\dd \vec{x}}_{\IV~(\vecc{G}^m)} \;+ \nonumber \\
&\qquad \underbrace{\sum_{E\in\setE_\Omega\cap\partial T_{k^-}} \int_E \vphi_{k^-i}\, 
    \frac{1}{2} \sum_{r=1}^2 \nu_{k^-}^r \sum_{m=1}^2 \left( 
      \sum_{l=1}^N D_{k^-l}^{rm} \, \vphi_{k^-l} \sum_{j=1}^N Q_{k^-j}^m \, \vphi_{k^-j} +
      \sum_{l=1}^N D_{k^+l}^{rm} \, \vphi_{k^+l} \sum_{j=1}^N Q_{k^+j}^m \, \vphi_{k^+j}
    \right) \dd \sigma}_{\V~(\vecc{R}^m)} \; + \nonumber\\
&\qquad \underbrace{\sum_{E\in\left(\setE_{\partial\Omega}\setminus\setE_\IQ\right)\cap\partial T_{k^-}} \int_E \vphi_{k^-i}\,\sum_{r=1}^2 \nu_{k^-}^r \sum_{m=1}^2 \sum_{l=1}^N D_{k^-l}^{rm} \, \vphi_{k^-l} \sum_{j=1}^N Q_{k^-j}^m \, \vphi_{k^-j} \,\dd\sigma}_{\VI~(\vecc{R}_\mathrm{bdr}^m)} \; + \;
  \underbrace{\sum_{E\in\setE_\IQ\cap\partial T_{k^-}} \int_E \vphi_{k^-i}\, q_D\, \dd \sigma}_{\VII~(\vec{J}_q)} \; + \nonumber\\
&\qquad \underbrace{\sum_{E\in\setE_\Omega\cap\partial T_{k^-}} \int_E \vphi_{k^-i} \, \begin{Bmatrix}
  \displaystyle \avg{ \vec{C}_u(h_{\Delta}, \vec{u}_{\Delta}) } \cdot \vec{\nu}_{k^-} & \mathrm{on}~\setE_\Omega^\mathrm{h} \\
  \widehat{C}_u(h_{\Delta}^-, \vec{u}_{\Delta}^-, h_{\Delta}^+, \vec{u}_{\Delta}^+) & \mathrm{on}~\setE_\Omega^\mathrm{v}
  \end{Bmatrix} \,\dd\sigma}_{\VIII~(\vecc{P}^m, \check{\vecc{Q}}, \vec{K}_u)} \; + \;
  \underbrace{\sum_{E\in\setE_{\partial\Omega}^\mathrm{h}\cap\partial T_{k^-}} \int_E \vphi_{k^-i} \, 
    \vec{C}_u(h_{\Delta}^-, \vec{u}_\mathrm{bdr}) \cdot \vec{\nu}_{k^-} \, \dd\sigma}_{\IX~(\vecc{P}^{m}_\mathrm{bdr}, \check{\vecc{Q}}_\mathrm{bdr}, \vec{J}^m_{\mathrm{bot}})} \nonumber \\
&\qquad \underbrace{\sum_{E\in\setE_{\partial\Omega}^\mathrm{v}\cap\partial T_{k^-}} \int_E \vphi_{k^-i} \, 
    \widehat{C}_u(h_{\Delta}^-, \vec{u}_{\Delta}^-, h_\mathrm{bdr}, \vec{u}_\mathrm{bdr}) \, \dd\sigma}_{\X~(\vecc{P}^m, \vecc{P}^m_\mathrm{bdr}, \check{\vecc{Q}}_\mathrm{bdr}, \vec{K}_{u}, \vec{J}_{h}, \vec{J}_{uu})}
\;=\; \underbrace{\int_{T_k} \vphi_{ki} \left( f_{\Delta} - g \, \partial_{x^1} \zb \right) \,\dd \vec{x}}_{\XI~(\vec{L}_u, \vec{L}_{\zb})} \;,
\label{eq:free:sys:u} \\
&\underbrace{\sum_{j=1}^N Q^m_{kj} \int_{T_k} \vphi_{ki} \, \vphi_{kj} \, \dd\vec{x}}_{\XII~(\vecc{M})} -
    \underbrace{\sum_{j=1}^N U^1_{kj} \int_{T_k} \partial_{x^m} \vphi_{ki} \, \vphi_{kj} \, \dd\vec{x}}_{\XIII~(\vecc{H}^m)} \;+\;
  \underbrace{\sum_{E\in\setE_{\partial\Omega}\cap\partial T_{k^-}} \int_E \vphi_{k^-i} \, \nu_{k^-}^m \, 
    u^1_{\mathrm{bdr}} \dd \sigma}_{\XIV~(\vecc{Q}^m, \vecc{Q}^m_\mathrm{bdr}, \vec{J}_u^m)} \;+ \nonumber \\
&\qquad \underbrace{\sum_{E\in\setE_\Omega\cap\partial T_{k^-}} \int_E \vphi_{k^-i} \, \nu_{k^-}^m \frac{1}{2}\left( 
    \sum_{j=1}^N U^1_{k^-j} \, \vphi_{k^-j} +
    \sum_{j=1}^N U^1_{k^+j} \, \vphi_{k^+j}
  \right) \,\dd\sigma}_{\XV~(\vecc{Q}^m)}
 = 0 \qquad \qquad \mathrm{for}~m\in\{1,2\}\;, 
\label{eq:free:sys:q} \\
-&\underbrace{\sum_{m=1}^2 \sum_{j=1}^N U^m_{kj} \int_{T_k} \partial_{x^m} \vphi_{ki} \, \vphi_{kj} \, \dd\vec{x}}_{\XVI~(\vecc{H}^m)} \; + 
  \underbrace{\sum_{E\in\setE_\Omega^\mathrm{v}\cap\partial T_{k^-}} \int_E \vphi_{k^-i} \, 
    \widehat{C}_h(\Hs, h_{\Delta}^-, {u_{\Delta}^1}^-, h_{\Delta}^+, {u_{\Delta}^1}^+)  \, \dd\sigma}_{\XVII~(\check{\vecc{P}}, \vec{K}_{h})} \; + \nonumber \\
&\qquad \underbrace{\sum_{E\in\setE_\Omega^\mathrm{h}\cap\partial T_{k^-}} \int_E \vphi_{k^-i} \,\frac{1}{2} \left(
    \sum_{j=1}^N U^1_{k^-j} \, \vphi_{k^-j} +
    \sum_{j=1}^N U^1_{k^+j} \, \vphi_{k^+j}
  \right) \, \nu_{k^-}^1 + \left(
    \sum_{j=1}^N U^2_{k^\uparrow j} \, \vphi_{k^\uparrow j}
  \right) \, \nu_{k^-}^2 \,\dd\sigma }_{\XVIII~(\vecc{Q}_\mathrm{avg}, \vecc{Q}_\mathrm{up})} \; +\nonumber \\
&\qquad \underbrace{\sum_{E\in\setE_{\partial\Omega}^\mathrm{h}\cap\partial T_{k^-}} \int_E \vphi_{k^-i} \, \left(
    u_\mathrm{bdr}^1 \, \nu_{k^-}^1 + u_{\uparrow}^2 \, \nu_{k^-}^2
  \right) \, \dd\sigma}_{\XIX~(\vecc{Q}_\mathrm{up}, \vecc{Q}_\mathrm{bdr}^m, \vec{J}_{u,\mathrm{bot}}^m)} \; + \;
  \underbrace{\sum_{E\in\setE_{\partial\Omega}^\mathrm{v}\cap\partial T_{k^-}}  \int_E \vphi_{k^-i} \, 
    \widehat{C}_h(\Hs, h_{\Delta}^-, {u_{\Delta}^1}^-, h_\mathrm{bdr}, u_\mathrm{bdr}^1)  \, \dd\sigma}_{\XX~(\check{\vecc{P}}, \check{\vecc{P}}_\mathrm{bdr}, \vec{K}_{h}, \vec{J}_{uh}, \check{\vec{J}}_u, \check{\vec{J}}_{h})}
= 0 \;.
\label{eq:free:sys:cont} 
\end{align}
Furthermore, we test the equation for the water height~\eqref{eq:free:semidisc:surf} with $\overline{w}_{\Delta} = \phi_{\overline{k}i}$ for $i \in \{1,\dots,\overline{N}\}$ and obtain another system of equations with a~contribution from $\overline{T}_{\overline{k}} = \Proj T_k$ given as
\begin{align}
& \underbrace{\sum_{j=1}^{\overline{N}} \partial_t H_{\overline{k}j} \int_{\overline{T}_{\overline{k}}} \phi_{\overline{k}i} \, \phi_{\overline{k}j} \, \dd x^1}_{\XXI~(\overline{\vecc{M}})} -
\underbrace{ \sum_{j=1}^{\overline{N}} H_{\overline{k}j} 
 \int_{\overline{T}_{\overline{k}}} \frac{1}{\Hs}  \partial_{x^1} \phi_{\overline{k}i} \, \left( \sum_{l=1}^{\overline{N}} \overline{U}_{\overline{k}l} \, \phi_{\overline{k}l} \right) \, \phi_{\overline{k}j} \, \dd x^1}_{\XXII~(\overline{\vecc{G}})}
\;+ \nonumber \\
&\qquad \underbrace{\sum_{E \in \setE_\Omega^\mathrm{v} \cap \Proj^{-1}\partial \overline{T}_{\overline{k}^-}} \int_E \phi_{\overline{k}^-i} \, 
    \widehat{C}_h(\Hs, h_{\Delta}^-, {u_{\Delta}^1}^-, h_{\Delta}^+, {u_{\Delta}^1}^+) \,\dd\sigma}_{\XXIII~(\overline{\vecc{P}},\overline{\vec{K}}_{h})} 
    \; + \nonumber\\
&\qquad     
  \underbrace{\sum_{E\in\setE_{\partial\Omega}^\mathrm{v} \cap \Proj^{-1}\partial \overline{T}_{\overline{k}^-}} \int_E \phi_{\overline{k}^-i} \, 
  \widehat{C}_h(\Hs, h_{\Delta}^-, {u_{\Delta}^1}^-, h_\mathrm{bdr}, u^1_\mathrm{bdr}) \, \dd\sigma}_{\XXIV~(\overline{\vecc{P}}, \overline{\vecc{P}}_\mathrm{bdr},\overline{\vec{K}}_{h}, \overline{\vec{J}}_{h}, \overline{\vec{J}}_u, \overline{\vec{J}}_{uh})}
~=~ 0 \;. 
\label{eq:free:sys:surf}  
\end{align}
\end{subequations}
}

Here, we deviate from our usual notation with cases inside edge integrals (used, e.\,g., in system~\eqref{eq:free:semidisc}) and use various sums over sets of edges of an~element that possess certain properties instead. 
This makes the presentation more compact and allows to refer directly to the relevant terms later on.
Note that these sets of edges can possibly be empty: for example, an interior element~$T_{k^-}$ with no edges on any domain boundaries results in an~empty set~$\setE_{\partial\Omega} \cap \partial T_{k^-} = \emptyset$, in which case the associated terms in system~$\eqref{eq:free:sys}$ drop out.

To ease relating the terms in system~\eqref{eq:free:sys} to the respective matrices and vectors that are presented explicitly in Appendix~\ref{sec:free-flow:mat-vec} and assembled in Appendix~\ref{sec:assembly}, we denote the matrices and vectors to which they contribute below each term.
Written in matrix form, this gives
{\allowdisplaybreaks
\begin{subequations}\label{eq:free:matsys}
\begin{align}
\vecc{M} \partial_t \vec{U}^1 =&
  \;\vec{L}_u - \vec{L}_{\zb} + \sum_{m=1}^2 \left( \vecc{G}^m - \vecc{R}^m - \vecc{R}^m_\mathrm{bdr} \right) \vec{Q}^m +
  g \left( \check{\vecc{H}} - \check{\vecc{Q}} - \check{\vecc{Q}}_\mathrm{bdr} \right) \vec{H} \nonumber \\
&  \;+ \sum_{m=1}^2 \left(\vecc{E}^m - \vecc{P}^m - \vecc{P}^m_\mathrm{bdr}\right) \vec{U}^m - \vec{K}_u - \vec{J}   \,, \\
\vecc{M} \vec{Q}^m =&
  \left( \vecc{H}^m - \vecc{Q}^m - \vecc{Q}^m_\mathrm{bdr} \right) \vec{U}^1
  -\vec{J}_{u}^m \qquad \mathrm{for}~m\in\{1,2\}\,, \\
\left( \vecc{H}^2 - \vecc{Q}_\mathrm{up} \right) \vec{U}^2 =&
  \left( -\vecc{H}^1 + \vecc{Q}_\mathrm{avg} + \vecc{Q}^1_\mathrm{bdr} + \check{\vecc{P}} + \check{\vecc{P}}_\mathrm{bdr} \right) \vec{U}^1 \nonumber \\
& \; +\vec{K}_{h} + \vec{J}_u^1 + \vec{J}_u^2 + 
  \frac{1}{2}\, \left(\check{\vec{J}}_{u} + \check{\vec{J}}_{h} + \vec{J}_{uh} \right) \,, \\
\overline{\vecc{M}} \partial_t \vec{H} =&
  \left(\overline{\vecc{G}} - \overline{\vecc{P}} - \overline{\vecc{P}}_\mathrm{bdr} \right) \vec{H} 
  - \overline{\vec{K}}_{h} -\frac{1}{2}\left( \overline{\vec{J}}_{h} + \overline{\vec{J}}_{u} + \overline{\vec{J}}_{uh} \right) \,,
\end{align}  
\end{subequations}
with 
$
\vec{J} \coloneqq \vec{J}_\mathrm{bot}^1 + \vec{J}_\mathrm{bot}^2 + \vec{J}_q + \frac{1}{2} (g \vec{J}_{h} + \vec{J}_{uu})
$
and representation vectors
\begin{align*}
\vec{U}^m(t) &\coloneqq \transpose{\begin{bmatrix} 
  U_{11}^m(t) & \cdots & U_{1N}^m(t) & \cdots & \cdots & U_{K1}^m(t) & \cdots & U_{KN}^m(t)
\end{bmatrix}} \,\in\IR^{KN} \quad \text{for } m\in\{1,2\}\,,\\
\vec{Q}^m(t) &\coloneqq \transpose{\begin{bmatrix} 
  Q_{11}^m(t) & \cdots & Q_{1N}^m(t) & \cdots & \cdots & Q_{K1}^m(t) & \cdots & Q_{KN}^m(t)
\end{bmatrix}}  \,\in\IR^{KN}  \quad \text{for } m\in\{1,2\}\,,\\
\vec{H}(t) &\coloneqq \transpose{\begin{bmatrix} 
  H_{11}(t) & \cdots & H_{1\overline{N}}(t) & \cdots & \cdots & H_{\overline{K}1}(t) & \cdots & H_{\overline{K}\overline{N}}(t)
\end{bmatrix}} \,\in\IR^{\overline{K}\overline{N}} \,.
\end{align*}}

Matrices with letters E to H correspond to element integrals, matrices with letters P to R to edge integrals, vectors with letter J represent contributions from Dirichlet boundary data, vectors with letter K stem from the jump term in the Lax--Friedrichs Riemann solver (see the next section), and vectors with letter L are contributed by right-hand side functions.
Furthermore, we use an~overline ($\overline{\,\cdot\,}$) to indicate matrices and vectors that originate from the one-dimensional free-surface equation~\eqref{eq:free:sys:surf} and mark with a~check ($\check{\,\cdot\,}$) all matrices in the two-dimensional equations~\eqref{eq:free:sys:u}--\eqref{eq:free:sys:cont} that concern the one-dimensional water height.
Occasional factors~$\frac{1}{2}$ are due to averaging in the Riemann solver.

\subsubsection{Approximation of non-linear fluxes}
\label{sec:riemann}
System~\eqref{eq:free:sys} contains both linear and non-linear fluxes.
The linear fluxes in terms~$\V-\VII$ and $\XIV-\XV$ stem from the diffusion operator in the momentum equation and can be approximated by central fluxes.
Non-linear ones appear in the form of the primitive fluxes in terms $\VIII-\X$, $\XVII-\XX$, and $\XXIII-\XXIV$.
These are treated in a~different way on vertical and horizontal edges, since the discrete free-surface elevation $H$ is discontinuous over the vertical edges but not over the horizontal ones.
On a~horizontal edge, we rely on simple averaging of the fluxes from both sides of the edge.
On vertical edges, an~approximation to the primitive fluxes is carried out with the help of a~Riemann solver.

We use the \textit{Lax--Friedrichs Riemann solver}, as it is one of the simplest Riemann-solvers that guarantees stability of our methods. It approximates a flux $C(c)$ for given primary variables $c^-,c^+$ as
\begin{equation}
\label{eq:lax-friedrichs}
\widehat{C}(c^-, c^+) = \avg{C(c)}  + \frac{1}{2} \abs{\hat{\lambda}} \jump{c} \cdot \vec{\nu} = \frac{1}{2} \Big( C(c^-) + C(c^+) \Big) + \frac{1}{2} \abs{\hat{\lambda}} \big( c^- - c^+ \big) \;.
\end{equation}
With this, the flux in terms $\VIII$ and $\X$ of the momentum equation~\eqref{eq:free:sys:u} is approximated by
\begin{equation*}
\widehat{C}_u(h_{\Delta}^-, \vec{u}_{\Delta}^-, h_{\Delta}^+, \vec{u}_{\Delta}^+) = \frac{1}{2} \Big( \vec{C}_u(h_{\Delta}^-, \vec{u}_{\Delta}^-) + \vec{C}_u(h_{\Delta}^+, \vec{u}_{\Delta}^+) \Big) \cdot \vec{\nu} + \frac{1}{2} \abs{\hat{\lambda}} \big( {u^1_{\Delta}}^- - {u^1_{\Delta}}^+ \big) \;,
\end{equation*}
and in terms $\XVII,\XX$ of the continuity equation~\eqref{eq:free:sys:cont} and $\XXIII,\XXIV$ of the free-surface equation~\eqref{eq:free:sys:surf} reads
\begin{equation*}
\widehat{C}_h(\Hs, h_{\Delta}^-, {u_{\Delta}^1}^-, h_{\Delta}^+, {u_{\Delta}^1}^+) = \frac{1}{2 \Hs} \Big( C_h(h_{\Delta}^-, {u_{\Delta}^1}^-) + C_h(h_{\Delta}^+, {u_{\Delta}^1}^+) \Big) \, \nu^1 + \frac{1}{2\Hs} \abs{\hat{\lambda}} \big( h_{\Delta}^- - h_{\Delta}^+ \big) \;.
\end{equation*}
Here, $\abs{\hat{\lambda}}$ is the largest (in absolute value) eigenvalue of the Jacobian of the primitive numerical fluxes
\begin{equation*}
\vec{C}(h_{\Delta}, \vec{u}_{\Delta}) = \begin{bmatrix} 
  C_h(h_{\Delta}, u_{\Delta}^1) \\ 
  \vec{C}_u(h_{\Delta}, \vec{u}_{\Delta}) \cdot \vec{\nu}
\end{bmatrix} \;,
\quad\text{where}\quad
\vec{C}'(\avg{h_{\Delta}}, \avg{\vec{u}_{\Delta}}) = \begin{bmatrix}
\avg{u_{\Delta}^1} \nu^1 & \avg{h_{\Delta}} \nu^1 \\
g \nu^1 & 2 \avg{u_{\Delta}^1} \nu^1
\end{bmatrix}
\end{equation*}
is the Jacobian w.r.t. variables $\avg{h_{\Delta}}$, $\avg{u_{\Delta}^1}$.
With $\vec{\nu} = \pm [1,0]$ on vertical edges, we obtain
\begin{align*}
\abs{\hat{\lambda}} &= 
\frac{3}{2} \abs{\avg{u_{\Delta}^1} \nu^1} + \frac{1}{2} \sqrt{ 9\left( \avg{u_{\Delta}^1}  \nu^1\right)^2 - 8\left( \avg{u_{\Delta}^1}  \nu^1\right)^2 + 4 g \avg{h_{\Delta}} \left( \nu^1 \right)^2 } \\
&= \frac{3}{2} \abs{\avg{u_{\Delta}^1}} + \frac{1}{2} \sqrt{ \left( \avg{u_{\Delta}^1} \right)^2 + 4 g \avg{h_{\Delta}}}\;.
\end{align*}
In semi-discrete form, the boundary flux $\widehat{C}_u(h_{\Delta}^-,\vec{u}_{\Delta}^-,h_{\Delta}^+,\vec{u}_{\Delta}^+)$ across vertical edges 
reads as
\begin{align*}
&\widehat{C}_u(h_{\Delta}^-,\vec{u}_{\Delta}^-,h_{\Delta}^+,\vec{u}_{\Delta}^+) \;=\; \frac{1}{2} \,\nu_{k^-}^1\, \left(
  \sum_{j=1}^N U^1_{k^-j} \, \vphi_{k^-j} \sum_{l=1}^{N} U^1_{k^-l} \, \vphi_{k^-l} +
  \sum_{j=1}^N U^1_{k^+j} \, \vphi_{k^+j} \sum_{l=1}^{N} U^1_{k^+l} \, \vphi_{k^+l}
\right)  \\
&\qquad\;+\;\frac{1}{2} \,\nu_{k^-}^1 \,g\, \left(
  \sum_{j=1}^{\overline{N}} H_{\overline{k}^-j} \, \phi_{\overline{k}^-j} + 
  \sum_{j=1}^{\overline{N}} H_{\overline{k}^+j} \, \phi_{\overline{k}^+j}
\right) \;+\; \frac{1}{2}\, \abs{\hat{\lambda}}\, \left(
  \sum_{j=1}^N U^1_{k^-j} \, \vphi_{k^-j} -
  \sum_{j=1}^N U^1_{k^+j} \, \vphi_{k^+j}
\right)\;.
\end{align*}
The boundary flux $\widehat{C}_h(\Hs, h_{\Delta}^-, {u_{\Delta}^1}^-, h_{\Delta}^+, {u_{\Delta}^1}^+)$ in equations~\eqref{eq:free:sys:cont},~\eqref{eq:free:sys:surf} is an approximation to the nonlinear boundary flux $\vec{u} \cdot \vec{\nu}$ across the vertical edges and reads in the semi-discrete form on $\partial T_{k^-} \cap \partial T_{k^+}$ as
\begin{align*}
\widehat{C}_h(\Hs, h_{\Delta}^-, {u_{\Delta}^1}^-, h_{\Delta}^+, {u_{\Delta}^1}^+) \;=\; &\frac{1}{2 \Hs} \,\nu_{k^-}^1 \,\left(
  \sum_{j=1}^N U^1_{k^-j} \, \vphi_{k^-j} \sum_{l=1}^{\overline{N}} H_{\overline{k}^-l} \, \phi_{\overline{k}^-l} +
  \sum_{j=1}^N U^1_{k^+j} \, \vphi_{k^+j} \sum_{l=1}^{\overline{N}} H_{\overline{k}^+l} \, \phi_{\overline{k}^+l}
\right) \;+ \\
&\qquad\frac{1}{2\Hs}\, \abs{\hat{\lambda}}\, \left(
  \sum_{j=1}^{\overline{N}} H_{\overline{k}^-j} \, \phi_{\overline{k}^-j} - 
  \sum_{j=1}^{\overline{N}} H_{\overline{k}^+j} \, \phi_{\overline{k}^+j}
\right)\;.
\end{align*}
On edges at the domain boundary, we employ Dirichlet data instead of the values from the neighboring element~$T_{k^+}$, where available, or use values from the interior of element~$T_{k^-}$.
In the latter case, the jump terms drop out, and the numerical fluxes become simply the one-sided primitive fluxes~$\vec{C}_u(h_{\Delta}^-, \vec{u}_{\Delta}^-)$ and~$C_h(h_{\Delta}^-, {u_{\Delta}^1}^-)$.

\subsubsection{Time discretization}

System~\eqref{eq:free:matsys} can be rewritten in the following way:
\begin{align*}
\vecc{M} \, \partial_t \vec{U}^1(t) =& \; \vec{S}_u(t) + \vecc{A}_{u,h}\, \vec{H}(t) + \sum_{m=1}^2 \left( \vecc{A}^m_{u,u}(t) \,\vec{U}^m(t) + \vecc{A}^m_{u,q}(t) \,\vec{Q}^m(t) \right)\,,\\
\vecc{M} \, \vec{Q}^m(t) =& \; \vec{S}_q^m(t) + \vecc{A}_{q}^m \vec{U}^1(t) \qquad\qquad \text{ for } m\in\{1,2\} \,,\\
\vecc{A}_{w,w} \, \vec{U}^2(t) =& \; \vec{S}_w(t) + \vecc{A}_{w,u}(t) \, \vec{U}^1(t) \,,\\
\overline{\vecc{M}} \, \partial_t \vec{H}(t) =& \; \vec{S}_{h}(t) + \vecc{A}_{h}(t) \, \vec{H}(t)\,,
\end{align*}
with
\begin{equation*}
\begin{array}{ll}
\vec{S}_u(t) \coloneqq\, \vec{L}_u(t) - \vec{L}_{\zb} - \vec{K}_u(t) - \vec{J}(t)\,, 
&\vecc{A}_{u,u}^m (t) \coloneqq\, \vecc{E}^m(t) - \vecc{P}^m(t) - \vecc{P}^m_\mathrm{bdr}(t)\,, \\
\vecc{A}_{u,h}\coloneqq\, \check{\vecc{H}} - \check{\vecc{Q}} - \check{\vecc{Q}}_\mathrm{bdr} \,, 
&\vecc{A}_{u,q}^m (t) \coloneqq\, \vecc{G}^m(t) - \vecc{R}^m(t) - \vecc{R}^m_\mathrm{bdr}(t)\,, \\
\vec{S}_q^m(t) \coloneqq\, -\vec{J}_u^m(t)\,, 
&\vecc{A}_q^m         \coloneqq\, \vecc{H}^m - \vecc{Q}^m - \vecc{Q}^m_\mathrm{bdr}\,,  \\
\vecc{A}_{w,u}(t) \;\coloneqq\, -\vecc{H}^1 + \vecc{Q}_\mathrm{avg} + \vecc{Q}^1_\mathrm{bdr} + \check{\vecc{P}}(t) + \check{\vecc{P}}_\mathrm{bdr}(t) \,,
&\vecc{A}_{w,w}\coloneqq\,  \vecc{H}^2 - \vecc{Q}_\mathrm{up} \,,\\
\multicolumn{2}{l}{\vec{S}_w(t) \coloneqq\, \vec{K}_{h}(t) + \vec{J}_u^1(t) + \vec{J}_u^2(t) + \frac{1}{2}\, \left(\check{\vec{J}}_{u}(t) + \check{\vec{J}}_{h}(t) + \vec{J}_{uh}(t) \right)\,,} \\
\vec{S}_{h}(t) \coloneqq\, - \overline{\vec{K}}_{h} -\frac{1}{2}\left( \overline{\vec{J}}_{h} + \overline{\vec{J}}_{u} + \overline{\vec{J}}_{uh} \right) \,,
&\vecc{A}_{h}(t) \coloneqq\, \overline{\vecc{G}} - \overline{\vecc{P}} - \overline{\vecc{P}}_\mathrm{bdr} \,,
\end{array}
\end{equation*}
for $m\in\{1,2\}$.
Here, we indicate matrices and vectors that have an~explicit time-depen\-dency; however, note that changes to the geometry of the domain (e.\,g., due to movement of the free surface) make re-assembly of all matrices necessary.

We discretize this system in time using the explicit Euler method.
In each time step, we first solve for diagnostic variables~$\vec{q}_{\Delta},u^2_{\Delta}$ at the previous time level~$t^n$ (using~$u^1_{\Delta}$ and~$h_{\Delta}$ from the previous time level) and use those to update water height~$h_{\Delta}$ and horizontal velocity~$u^1_{\Delta}$ at the new time level, i.\,e., a time step implements the following scheme:

Let $0 = t^1 < t^2 < \cdots \leq t_\mathrm{end}$ be a~not necessarily equidistant decomposition of time interval~$J$, and let~$\Delta t^n \coloneqq t^{n+1} - t^n$ denote the time step size.
The update scheme is given by
\begin{subequations}\label{eq:free:timedepsys}
\begin{align}
\vec{Q}^m(t^n) =& \; \vecc{M}^{-1} \left( \vec{S}_q^m(t^n) + \vecc{A}_{q}^m \vec{U}^1(t^n) \right) \qquad\qquad \text{ for } m\in\{1,2\} \,,\\
\vec{U}^2(t^n) =& \; \vecc{A}_{w,w}^{-1} \left( \vec{S}_w(t^n) + \vecc{A}_{w,u}(t^n) \, \vec{U}^1(t^n) \right) \,,\\
\vec{U}^1(t^{n+1}) =& \; \vec{U}^1(t^n) + \Delta t^{n+1} \, \vecc{M}^{-1} \Big[ \vec{S}_u(t^n) + \vecc{A}_{u,h}\, \vec{H}(t^n) \nonumber \\
& \qquad\qquad + \sum_{m=1}^2 \left( \vecc{A}^m_{u,u}(t^n) \,\vec{U}^m(t^n) + \vecc{A}^m_{u,q}(t^n) \,\vec{Q}^m(t^n) \right) \Big]\,,\\
\vec{H}(t^{n+1}) =& \; \vec{H}(t^n) + \Delta t^{n+1} \, \overline{\vecc{M}}^{-1} \Bigg[ \vec{S}_{h}(t^n) + \vecc{A}_{h}(t^n) \, \vec{H}(t^n)\Bigg] \,.
\end{align}
\end{subequations}

In each time step, the $x^2$-coordinates of the surface nodes of the mesh are adjusted according to follow the movement of the free surface as expressed by the height~$h_{\Delta}$ at that time level.
That way, the shape of the top-most elements changes over the course of the simulation and requires to re-assemble stationary matrices in every time step.

Note that FESTUNG has higher-order time integration methods built-in~\cite{FESTUNG2}, which could be used instead of the first-order Euler method.
However, this introduces either a~consistency error due to delayed mesh updates, or the mesh adaptation would have to be applied in each substep of the Runge-Kutta method.
If a small consistency error is acceptable, e.\,g., in cases where the free surface does not change much per time step, computation time can be reduced by executing multiple time steps before adapting the surface nodes.
In our implementation, however, we found the assembly of stationary matrices to be affordable compared to other matrix terms and thus adapt the free surface in every time step.

A~limitation of our current implementation is the fact that surface elements are not allowed to \enquote{dry up}, i.\,e., the movement of the free surface towards the bottom is not allowed to exceed the vertical extent of the top-most elements.

\section{Generic problem framework}
\label{sec:problem-framework}

\begin{figure}[!tb]
\centering
\includegraphics[width=\textwidth]{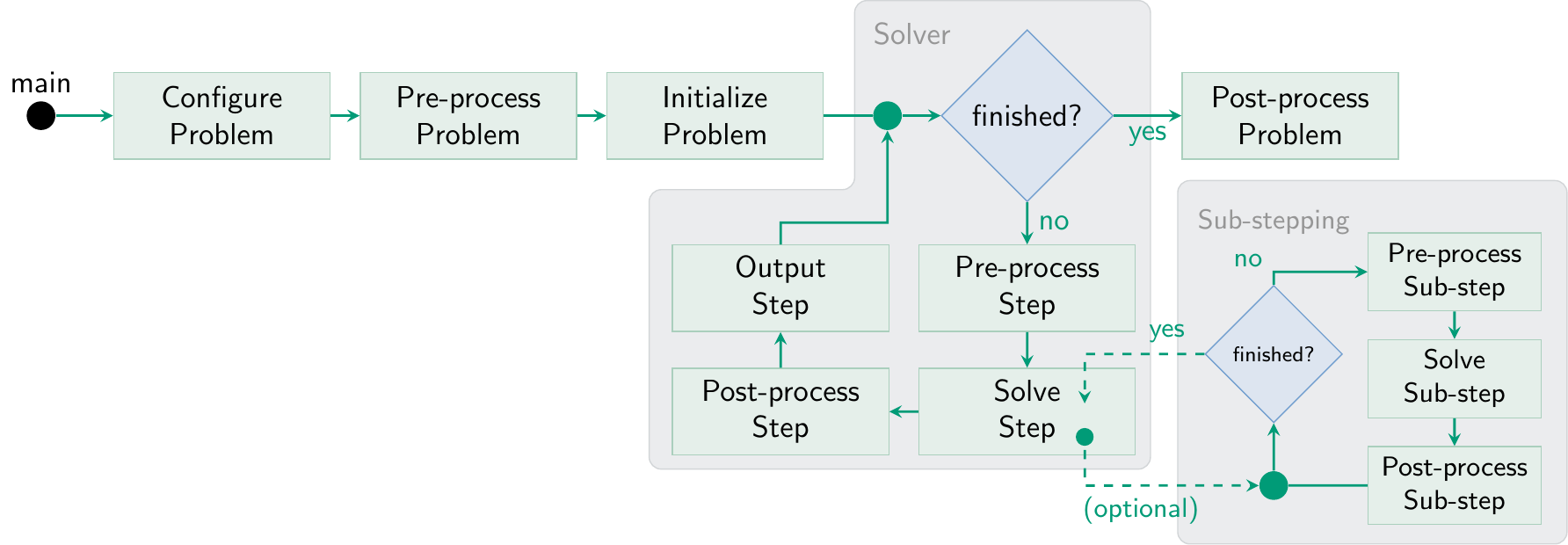}
\caption{The generic solver formulation, including the optional sub-stepping in the solver phase.}
\label{fig:solver}
\end{figure}

The core of our model coupling framework is the generic solver formulation that we first applied in our previous work in series~\cite{FESTUNG3} and properly introduce in this version of FESTUNG.
It is based on the observation that almost every solver for a PDE problem can be subdivided into three major steps:
\begin{enumerate}[nosep]
 \item A~\emph{setup phase}, which defines problem parameters or reads them from a~configuration file, allocates grid data structures, initializes solution vectors, etc.;
 \item an~(iterative) \emph{solver phase}, e.\,g., to build and solve a~linear system (possibly repeatedly for many time steps), or apply an~iterative method for non-linear problems;
 \item a~\emph{post-processing phase} to evaluate errors, write the computed solution to a~file for visualization, etc.
\end{enumerate}
To further structure the setup phase, we split it into a~\emph{configuration} step (definition of problem parameters and array sizes), a~\emph{problem pre-processing} step (assembly of static data structures), and an~\emph{initialization} step (projection of initial data).

The solver phase takes care of the actual work and is designed to be passed through repeatedly, controlled by a~single parameter that indicates whether another iteration of this phase is to be executed.
This way, the number of iterations does not have to be known {\em a~priori}, thus adaptive time stepping and stationary problems can be elegantly implemented with the latter simply marked as finished after one iteration.
Again, this phase is split into sub-steps: a~\emph{pre-processing} step, a~\emph{solver} step, a~\emph{post-processing} step, and an~\emph{output} step.
Any of the stages can be empty if not required by the solution algorithm.
For solvers with nested iterations (e.\,g., a~Newton's method in each time step or multiple stages of a~Runge-Kutta method) we provide an~optional sub-stepping functionality that can be executed, e.\,g., as part of the solver step.
It does, again, carry out pre-processing, solving, and post-processing steps.
The resulting solver structure is outlined in Fig.~\ref{fig:solver}.

Each of these steps is implemented as a~separate \MatOct-function, all of which are stored together in a~sub-folder.
A~simulation run is driven by a~generic \code{main}-Function that is given the name of the folder and executes all steps in the defined order.
Data needed inbetween steps is stored in a~\code{struct} that is passed to each step function together with the current iteration number.

Both, free-flow problem and subsurface problem are implemented individually in this fashion.
This allows to verify, debug, and re-use them as needed.
The coupled problem itself is implemented as a~solver of the same structure that executes the steps of both sub-problems at the appropriate times and maintains separate \code{struct} instances for each sub-problem.
The coupled solver, for example, takes care of executing multiple time steps of one problem during each time step of the other (see Sec.~\ref{sec:coupling}) using the sub-stepping functionality depicted in Fig.~\ref{fig:solver}.
The resulting program flow is shown in Fig.~\ref{fig:solver:coupled}.
All relevant data for the coupling at the interface is passed from one subproblem to the other before and after this sub-stepping by updating each problem's \code{struct} data with the interface vectors described in Sec.~\ref{sec:coupling}.
This way, the number of changes to the sub-problems is minimal and (almost) all coupling logic is kept separately in the coupled solver improving the readability and maintainability.
Particular care has to be exercised because functions with the same name (i.e., the steps for each individual sub-problem and the coupled solver) are contained in different folders -- making it impossible to have all relevant folders in the search path at the same time.
To overcome the overhead involved with changes to the path-variable, function handles to each of the required steps and functions from the folders are stored in the beginning and used in the iterative solver phase.
We found this requirement to be crucial in eliminating all significant overhead of the new framework.

When comparing the runtimes of the solvers implemented in our first papers in series~\cite{FESTUNG1,FESTUNG2} before and after migration to the new framework, we experienced that this added less than~$1\,\mathrm{ms}$ of additional runtime per time step on a~regular desktop computer (Intel Core i7-4790, MATLAB 2017b) that can be attributed to the additional overhead of calling functions for each step.
Compared to the computational cost of the numerical scheme itself, this was in the order of 0.1--1$\%$ of the total runtime for these rather simple models.
For more complicated problems, where the simulation time per time step increases further, this becomes even less significant.

\section{Coupling of subsurface and free-flow problems}
\label{sec:coupling}

In this section, we introduce the necessary changes to the discretizations of free-flow and subsurface problems when coupling both models using the interface conditions given in Sec.~\ref{sec:model:coupled}.
To distinguish between variables, values, etc. belonging to the subsurface and the free-flow problems, we once again indicate the subsurface entities  by a~tilde `$\,\tilde \cdot\,$'. 
The computational mesh of the coupled domain
(see~Fig.~\ref{fig:Omega}) is for simplicity restricted to have matching elements in horizontal direction between the subsurface and free-flow parts.
Vertical boundaries are straight vertical lines due to the requirements for the discretization of the free-flow problem (see Sec.~\ref{sec:gridbasis}).
Horizontal boundaries are allowed to be piecewise linear.
In the following, we denote the set of edges in the subsurface problem on the interior boundary as~$\widetilde{\setE}_\mathrm{int} \coloneqq \{ E \in \setE_{\partial\tO} \mid E \subset \Gamma_\mathrm{int} \}$ and, correspondingly, in the free-flow problem as~$\setE_\mathrm{int} \coloneqq \{ E \in \setE_{\partial\Omega(t)} \mid E \subset \Gamma_\mathrm{int} \}$.

A challenge for coupled simulations is the difference in time scales.
At the surface, the water velocity is in the range of meters per second, while the subsurface flow velocity is in the range of decimeters per day.
Consequently, the free-flow problem's time step must be significantly smaller than the one for the subsurface problem.
Moreover, we must take into account the fact that the free-flow problem is discretized in time using an~explicit Euler method, whereas the subsurface problem relies on an~implicit Euler method, and that the coupling must be mass conservative.

Since $\Delta t / \Delta \tilde t \le 1$, the flux from~$\Omega(t)$ to~$\tO$ has to be time-averaged to preserve the conservation in the interior boundary condition~\eqref{eq:coupling:pm}. 
We use a~so-called \emph{non-simple} or \emph{complex} boundary, which memorizes the flux from $\tO$ into $\Omega(t)$ across the interior boundary, and vice versa, for all sub-steps.
Instead of time-averaging, one could also use water height and horizontal velocity of the latest time step to re-implement some sort of standard implicit Euler scheme for the subsurface problem.
For consistency, we require $\Delta \tilde t = n_\mathrm{substep} \, \Delta t$ with $n_\mathrm{substep} \in \IN$, which allows us to approximate the time average using a~summed trapezoidal rule:

\begin{equation}\label{eq:couple:time-avg}
\begin{aligned}
\overline{\thh_\mathrm{D}} = &
\frac{1}{\Delta \tilde t} \int_{t^n}^{t^{n+1}} \left(
  h_{\Delta} (t) + \frac{\left(u^{1}(t) \right)^2}{2}  \right) \dd t \\
\approx &
\frac{\Delta t}{2 \Delta \tilde t} \sum_{i=1}^{n_\mathrm{substep}} \left(
  h_{\Delta}(t^{n,i-1}) + \frac{\left( u_{\Delta}^{1}(t^{n,i-1}) \right)^2}{2} +
  h_{\Delta}(t^{n,i}) + \frac{\left(u_{\Delta}^{1}(t^{n,i}) \right)^2}{2} 
  \right)\,,
\end{aligned}
\end{equation}
where~$\Delta \tilde t = t^{n+1} - t^n$, and~$t^{n,i}\coloneqq t^n + i \, \Delta t$.
The flux over the interior boundary for interface condition~\eqref{eq:coupling:ff} is always taken from the latest time step of the subsurface problem rescaled with the corresponding time step size. 
This ensures conservation of mass with a~time lag of one time step $\Delta \tilde t$, i.e., if the subsurface problem looses mass, this mass is gained by the free-flow problem in the following $n_\mathrm{substep}$ time steps of $\Delta t$.

\subsection{Changes to the subsurface problem}
\label{sec:coupling:pm}

Interior boundary condition~\eqref{eq:coupling:pm} enters Eq.~\eqref{eq:sub:semidisc:u} as a replacement for Dirichlet boundary contributions for the hydraulic head on edge integrals.
Thus, it produces contributions equivalent to~$\tilde{\vec{J}}^m_\mathrm{D}$ and~$\tilde{\vec{K}}_\mathrm{D}$.
To make the coupling as transparent as possible, we introduce vectors~$\tilde{\vec{J}}^m_\mathrm{int}, \tilde{\vec{K}}_\mathrm{int}$ at the same places in system~\eqref{eq:sub:matsys}, which are by default set to zero.
In the case of a~coupled simulation, the coupled solver fills them in each time step of the subsurface problem with the updated values, which are computed as
\begin{align*}
[\tilde{\vec{J}}^m_{\mathrm{int}}]_{(k-1)N+i} &\;\coloneqq\; 
\sum_{E_{kn}\in\partial T_k\cap\widetilde{\setE}_{\mathrm{int}}} 
\nu_{kn}^m \int_{E_{kn}} \vphi_{ki}\,\overline{\thh_\mathrm{D}}\,\dd\sigma\;,
\\
[\tilde{\vec{K}}_\mathrm{int}]_{(k-1)N+i}  
&\;\coloneqq\;   \sum_{E_{kn}\in\partial T_k\cap\widetilde{\setE}_{\mathrm{int}}} \frac{1}{\abs{E_{kn}}} \int_{E_{kn}} \vphi_{ki} \,\overline{\thh_\mathrm{D}}\,\dd\sigma\;.
\end{align*}
In our implementation, we apply a~quadrature rule (cf.~Sec.~\ref{sec:quadrature}) to evaluate the entries of~$\tilde{\vec{J}}^m_\mathrm{int}$, $\tilde{\vec{K}}_\mathrm{int}$ and determine the time-averaged value for~$\overline{\thh_\mathrm{D}}$ in Eq.~\eqref{eq:couple:time-avg} directly as per-quadrature-point values.

\subsection{Changes to the free-flow problem}
\label{sec:coupling:ff}

The interior boundary condition for the free-flow problem~\eqref{eq:coupling:ff} is of Dirichlet type for the velocity~$\vec{u}$ and enters flux equation~\eqref{eq:free:sys:q} as Dirichlet boundary data for horizontal velocity~$u^1$ and momentum and continuity equations~\eqref{eq:free:sys:u},~\eqref{eq:free:sys:cont} as Dirichlet boundary data for horizontal and vertical velocities.
It replaces the contributions to vectors~$\vec{J}^m_u$ (term~$\XIV$), $\vec{J}_\mathrm{bot}^m$ (term~$\IX$), and $\vec{J}^m_{u,\mathrm{bot}}$ on interior boundary edges.
As in the subsurface problem, we introduce vectors~$\vec{J}^m_{u,\mathrm{int}}$, $\vec{J}_{w,\mathrm{int}}$, and~$\vec{J}^m_{uu,\mathrm{int}}$, which are initialized to be zero and filled by the coupled solver to add the relevant contributions in coupled simulations, with entries due to coupling condition~\eqref{eq:coupling:ff} given by
\begin{align*}
\left[\vec{J}^m_{u,\mathrm{int}}\right]_{(k-1)N+i} \coloneqq &
  \sum_{E_{kn}\in\partial T_k \cap \setE_\mathrm{int}} \nu_{kn}^m 
  \int_{E_{kn}} \vphi_{ki} \, u_\mathrm{D}^{1} \,\dd\sigma 
  = 
  \sum_{E_{kn}\in\partial T_k \cap \setE_\mathrm{int}} \nu_{kn}^m 
  \int_{E_{kn}} \vphi_{ki} \, \left( \left[\vecc{\widetilde D}\right]_{1,:} \cdot \tq_{\Delta} \right) \,\dd\sigma \,,\\
\left[\vec{J}_{w,\mathrm{int}}\right]_{(k-1)N+i} \coloneqq&
  \sum_{E_{kn}\in\partial T_k \cap \setE_\mathrm{int}} \nu_{kn}^2 
  \int_{E_{kn}} \vphi_{ki} \, u_\mathrm{D}^{2} \,\dd\sigma 
  = 
  \sum_{E_{kn}\in\partial T_k \cap \setE_\mathrm{int}} \nu_{kn}^2 
  \int_{E_{kn}} \vphi_{ki} \, \left( \left[\vecc{\widetilde D}\right]_{2,:} \cdot \tq_{\Delta} \right) \,\dd\sigma \,,\\
\left[\vec{J}^m_{uu,\mathrm{int}}\right]_{(k-1)N+i} \coloneqq &
  \sum_{E_{kn}\in\partial T_k \cap \setE_\mathrm{int}} \nu_{kn}^m 
  \int_{E_{kn}} \vphi_{ki} \, u_\mathrm{D}^{1} \, u_\mathrm{D}^{m} \,\dd\sigma \\
  = &
  \sum_{E_{kn}\in\partial T_k \cap \setE_\mathrm{int}} \nu_{kn}^m 
  \int_{E_{kn}} \vphi_{ki} \, \left( \left[\vecc{\widetilde D}\right]_{1,:} \cdot \tq_{\Delta} \right) \,\left( \left[\vecc{\widetilde D}\right]_{m,:} \cdot \tq_{\Delta} \right) \,\dd\sigma \,.
\end{align*}
In our implementation, we apply a~quadrature rule (cf.~\ref{sec:quadrature}) to directly integrate the entries of these vectors.

\subsection{Coupled algorithm}

\begin{figure}
\centering
\includegraphics[width=\textwidth]{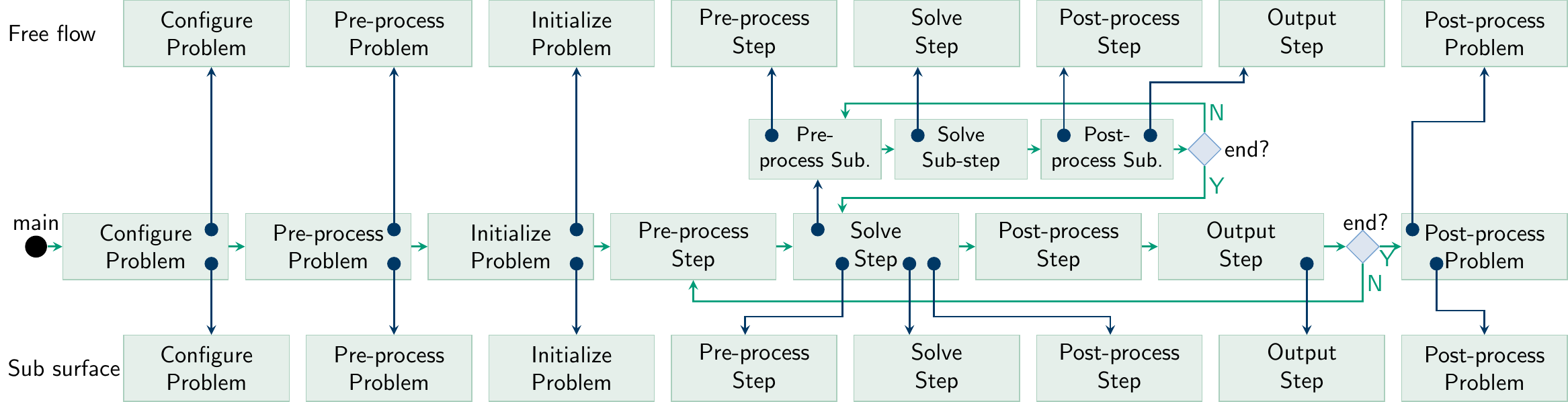}
\caption{Program flow (from left to right) of the coupled free-flow and subsurface solver. Arrows in {\color{faunat}green} indicate program flow advanced by the main driver routine, arrows in {\color{faublue}blue} (with circles at the anchor) represent calls to the respective subroutine.}
\label{fig:solver:coupled}
\end{figure}

The coupled solver is implemented as an instance of the generic problem framework described in Section~\ref{sec:problem-framework}, i.e., it provides all steps depicted in Figure~\ref{fig:solver} as \MatOct\ routines.
In each of these routines, it calls the relevant step routines of the free-flow and subsurface problem and assembles coupling terms where necessary.
The resulting program flow is shown in Figure~\ref{fig:solver:coupled}: 
The initialization phase consists of the respective initialization steps of each sub-problem and of the computation of grid data structures that allow to match the mesh entities on the top of the subsurface and bottom of the free-flow grids.
Each iteration of the time stepping loop begins with resetting the memory variable of the time-averaged boundary condition for the hydraulic head (Eq.~\eqref{eq:couple:time-avg}) and assembling the contributions to the free-flow problem (cf.~Sec.~\ref{sec:coupling:ff}).
The majority of work is done in the solver step, which first executes multiple time steps of the free-flow problem (implemented as sub-steps) before computing the contributions of the coupling interface to the subsurface problem (see Sec.~\ref{sec:coupling:pm}) and carrying out its matching time step.
That way, the free-flow problem gets advanced in time first (in an explicit manner) before bringing the subsurface problem to the same time level in a~single implicit time step.

\section{Numerical results}
\label{sec:results}

The performance of our implementation is demonstrated using two types of problem: 
Analytical convergence tests to verify the solver and a~more realistic setup.

\subsection{Analytical convergence tests}

\begin{table}[!ht]\revised{
\scriptsize
\setlength{\tabcolsep}{2pt}
\begin{tabular}{@{}ccccccccccccccc@{}}
\toprule
$p$ & $j$ & $\Err{h}$ & $\EOC{h}$ & $\Err{u^1}$ & $\EOC{u^1}$ & $\Err{u^2}$ & $\EOC{u^2}$ & $\Err{\thh}$ & $\EOC{\thh}$ & $\Err{\tilde{q}^1}$ & $\EOC{\tilde{q}^1}$ & $\Err{\tilde{q}^2}$ & $\EOC{\tilde{q}^2}$ \\
\toprule
\multirow{5}{*}{0}
& 0 & 7.15e--01 &  --- &  8.89e--01 &  --- &  2.00e--01 &  --- &  1.50e+01  &  --- &  7.78e--03 &  --- &  2.07e+00  & \phantom{--} --- \\
& 1 & 3.59e--01 & 0.99 &  6.76e--01 & 0.40 &  1.74e--01 & 0.20 &  8.29e+00  & 0.86 &  5.27e--03 & 0.56 &  1.11e+00  & \phantom{--}0.90 \\
& 2 & 1.79e--01 & 1.00 &  3.47e--01 & 0.96 &  9.13e--02 & 0.93 &  4.21e+00  & 0.98 &  2.44e--03 & 1.11 &  5.30e--01 & \phantom{--}1.07 \\
& 3 & 8.97e--02 & 1.00 &  1.76e--01 & 0.98 &  5.69e--02 & 0.68 &  2.11e+00  & 1.00 &  1.15e--03 & 1.09 &  2.58e--01 & \phantom{--}1.04 \\
& 4 & 4.49e--02 & 1.00 &  8.81e--02 & 1.00 &  3.83e--02 & 0.57 &  1.06e+00  & 1.00 &  5.63e--04 & 1.03 &  1.27e--01 & \phantom{--}1.02 \\
\midrule
\multirow{5}{*}{1}
& 0 & 1.06e--02 &  --- &  7.01e--01 &  --- &  9.44e--02 &  --- &  5.49e+00  &  --- &  4.17e--03 &  --- &  3.61e--01 & \phantom{--} --- \\
& 1 & 3.11e--03 & 1.76 &  1.76e--01 & 1.99 &  5.79e--02 & 0.71 &  1.32e+00  & 2.05 &  2.86e--03 & 0.54 &  6.02e--01 &           --0.74 \\
& 2 & 7.92e--04 & 1.97 &  4.89e--02 & 1.85 &  3.65e--02 & 0.67 &  3.30e--01 & 2.00 &  2.00e--03 & 0.51 &  4.19e--01 & \phantom{--}0.52 \\
& 3 & 1.99e--04 & 1.99 &  1.24e--02 & 1.98 &  2.05e--02 & 0.83 &  8.24e--02 & 2.00 &  1.14e--03 & 0.82 &  2.34e--01 & \phantom{--}0.84 \\
& 4 & 4.99e--05 & 2.00 &  3.11e--03 & 1.99 &  1.06e--02 & 0.95 &  2.06e--02 & 2.00 &  5.92e--04 & 0.94 &  1.22e--01 & \phantom{--}0.94 \\
\midrule
\multirow{5}{*}{2}
& 0 & 3.87e--03 &  --- &  2.46e--01 &  --- &  8.76e--02 &  --- &  6.33e--01 &  --- &  3.03e--03 &  --- &  9.41e--02 & \phantom{--} --- \\
& 1 & 5.07e--04 & 2.93 &  4.39e--02 & 2.49 &  2.56e--02 & 1.77 &  9.19e--02 & 2.78 &  6.72e--04 & 2.17 &  3.71e--02 & \phantom{--}1.34 \\
& 2 & 6.56e--05 & 2.95 &  5.42e--03 & 3.02 &  3.24e--03 & 2.98 &  1.17e--02 & 2.97 &  8.91e--05 & 2.92 &  4.91e--03 & \phantom{--}2.92 \\
& 3 & 8.26e--06 & 2.99 &  6.81e--04 & 2.99 &  4.13e--04 & 2.97 &  1.48e--03 & 2.99 &  1.10e--05 & 3.02 &  6.21e--04 & \phantom{--}2.98 \\
& 4 & 1.03e--06 & 3.00 &  8.52e--05 & 3.00 &  5.98e--05 & 2.79 &  1.85e--04 & 3.00 &  1.37e--06 & 3.01 &  7.83e--05 & \phantom{--}2.99 \\
\midrule
\multirow{5}{*}{3}
& 0 & 9.06e--04 &  --- &  9.68e--02 &  --- &  2.99e--02 &  --- &  1.21e--01 &  --- &  3.54e--04 &  --- &  8.00e--03 & \phantom{--} --- \\
& 1 & 6.65e--05 & 3.77 &  5.77e--03 & 4.07 &  3.84e--03 & 2.96 &  7.24e--03 & 4.06 &  1.18e--04 & 1.58 &  8.47e--03 &           --0.08 \\
& 2 & 4.23e--06 & 3.98 &  4.02e--04 & 3.84 &  7.28e--04 & 2.40 &  4.51e--04 & 4.01 &  2.36e--05 & 2.32 &  1.43e--03 & \phantom{--}2.57 \\
& 3 & 2.66e--07 & 3.99 &  2.55e--05 & 3.98 &  1.04e--04 & 2.81 &  2.81e--05 & 4.00 &  3.40e--06 & 2.80 &  1.97e--04 & \phantom{--}2.86 \\
& 4 & 1.74e--08 & 3.93 &  1.60e--06 & 3.99 &  1.34e--05 & 2.96 &  1.76e--06 & 4.00 &  4.44e--07 & 2.93 &  2.56e--05 & \phantom{--}2.95 \\
\midrule
\multirow{3}{*}{4}
& 0 & 2.02e--04 &  --- &  1.99e--02 &  --- &  1.32e--02 &  --- &  8.27e--03 &  --- &  3.08e--04 &  --- &  1.24e--03 & \phantom{--} --- \\
& 1 & 6.46e--06 & 4.96 &  8.82e--04 & 4.50 &  1.08e--03 & 3.62 &  3.00e--04 & 4.78 &  1.83e--05 & 4.08 &  2.68e--04 & \phantom{--}2.21 \\
& 2 & 2.08e--07 & 4.96 &  2.69e--05 & 5.03 &  3.61e--05 & 4.90 &  9.59e--06 & 4.97 &  6.49e--07 & 4.81 &  9.39e--06 & \phantom{--}4.83 \\
\bottomrule
\end{tabular}
\caption{$L^2(\Omega)$ discretization errors and estimated orders of convergence for the sub-problems. On the $j$th refinement level, we used $2^{j+1} \times 2^j$ elements and time step~$\Delta \tilde{t} = 4\cdot10^{-5} \cdot 2^{-p \, (j+1)}$ for the subsurface problem and~$\Delta t = 4\cdot10^{-6} \cdot 2^{-p \, (j+1)}$ for the free-flow problem.}}
\label{tab:conv:ffpm}
\end{table}

\begin{table}[!ht]\revised{
\scriptsize
\setlength{\tabcolsep}{2pt}
\begin{tabular}{@{}ccccccccccccccc@{}}
\toprule
$p$ & $j$ & $\Err{h}$ & $\EOC{h}$ & $\Err{u^1}$ & $\EOC{u^1}$ & $\Err{u^2}$ & $\EOC{u^2}$ & $\Err{\thh}$ & $\EOC{\thh}$ & $\Err{\tilde{q}^1}$ & $\EOC{\tilde{q}^1}$ & $\Err{\tilde{q}^2}$ & $\EOC{\tilde{q}^2}$ \\
\toprule
\multirow{5}{*}{0}
& 0 & 7.15e--01 &  --- &  8.89e--01 &  --- &  2.01e--01 &  --- &  1.50e+01  &  --- &  7.78e--03 &  --- &  2.07e+00  & \phantom{--} --- \\
& 1 & 3.59e--01 & 0.99 &  6.76e--01 & 0.40 &  1.74e--01 & 0.21 &  8.29e+00  & 0.86 &  5.27e--03 & 0.56 &  1.11e+00  & \phantom{--}0.90 \\
& 2 & 1.79e--01 & 1.00 &  3.47e--01 & 0.96 &  9.12e--02 & 0.93 &  4.21e+00  & 0.98 &  2.44e--03 & 1.11 &  5.30e--01 & \phantom{--}1.07 \\
& 3 & 8.97e--02 & 1.00 &  1.76e--01 & 0.98 &  5.69e--02 & 0.68 &  2.11e+00  & 1.00 &  1.15e--03 & 1.09 &  2.58e--01 & \phantom{--}1.04 \\
& 4 & 4.49e--02 & 1.00 &  8.81e--02 & 1.00 &  3.83e--02 & 0.57 &  1.06e+00  & 1.00 &  5.63e--04 & 1.03 &  1.27e--01 & \phantom{--}1.02 \\
\midrule
\multirow{5}{*}{1}
& 0 & 1.06e--02 &  --- &  7.01e--01 &  --- &  9.46e--02 &  --- &  5.49e+00  &  --- &  4.17e--03 &  --- &  3.61e--01 & \phantom{--} --- \\
& 1 & 3.11e--03 & 1.76 &  1.76e--01 & 1.99 &  5.78e--02 & 0.71 &  1.32e+00  & 2.05 &  2.86e--03 & 0.54 &  6.01e--01 &           --0.74 \\
& 2 & 7.92e--04 & 1.97 &  4.89e--02 & 1.85 &  3.65e--02 & 0.67 &  3.30e--01 & 2.00 &  2.00e--03 & 0.51 &  4.19e--01 & \phantom{--}0.52 \\
& 3 & 1.99e--04 & 1.99 &  1.24e--02 & 1.98 &  2.05e--02 & 0.83 &  8.24e--02 & 2.00 &  1.14e--03 & 0.82 &  2.34e--01 & \phantom{--}0.84 \\
& 4 & 4.99e--05 & 2.00 &  3.11e--03 & 1.99 &  1.06e--02 & 0.95 &  2.06e--02 & 2.00 &  5.92e--04 & 0.94 &  1.22e--01 & \phantom{--}0.94 \\
\midrule
\multirow{5}{*}{2}
& 0 & 3.87e--03 &  --- &  2.46e--01 &  --- &  8.76e--02 &  --- &  6.33e--01 &  --- &  3.03e--03 &  --- &  9.42e--02 & \phantom{--} --- \\
& 1 & 5.07e--04 & 2.93 &  4.39e--02 & 2.49 &  2.56e--02 & 1.77 &  9.19e--02 & 2.78 &  6.72e--04 & 2.17 &  3.72e--02 & \phantom{--}1.34 \\
& 2 & 6.56e--05 & 2.95 &  5.42e--03 & 3.02 &  3.24e--03 & 2.98 &  1.17e--02 & 2.97 &  8.91e--05 & 2.92 &  4.91e--03 & \phantom{--}2.92 \\
& 3 & 8.26e--06 & 2.99 &  6.81e--04 & 2.99 &  4.17e--04 & 2.96 &  1.48e--03 & 2.99 &  1.10e--05 & 3.02 &  6.21e--04 & \phantom{--}2.98 \\
& 4 & 1.03e--06 & 3.00 &  8.52e--05 & 3.00 &  7.16e--05 & 2.54 &  1.85e--04 & 3.00 &  1.37e--06 & 3.01 &  7.83e--05 & \phantom{--}2.99 \\
\midrule
\multirow{5}{*}{3}
& 0 & 9.06e--04 &  --- &  9.68e--02 &  --- &  2.99e--02 &  --- &  1.21e--01 &  --- &  3.55e--04 &  --- &  9.13e--03 & \phantom{--} --- \\
& 1 & 6.65e--05 & 3.77 &  5.77e--03 & 4.07 &  3.84e--03 & 2.96 &  7.24e--03 & 4.06 &  1.18e--04 & 1.59 &  8.55e--03 & \phantom{--}0.10 \\
& 2 & 4.23e--06 & 3.98 &  4.02e--04 & 3.84 &  7.29e--04 & 2.40 &  4.51e--04 & 4.01 &  2.36e--05 & 2.32 &  1.43e--03 & \phantom{--}2.58 \\
& 3 & 2.66e--07 & 3.99 &  2.55e--05 & 3.98 &  1.10e--04 & 2.73 &  2.81e--05 & 4.00 &  3.40e--06 & 2.80 &  1.97e--04 & \phantom{--}2.86 \\
& 4 & 1.74e--08 & 3.93 &  1.60e--06 & 3.99 &  3.92e--05 & 1.49 &  1.76e--06 & 4.00 &  4.44e--07 & 2.93 &  2.56e--05 & 2.94 \\
\midrule
\multirow{3}{*}{4}
& 0 & 2.02e--04 &  --- &  1.99e--02 &  --- &  1.32e--02 &  --- &  8.27e--03 &  --- &  3.10e--04 &  --- &  4.94e--03 & \phantom{--} --- \\
& 1 & 6.46e--06 & 4.96 &  8.82e--04 & 4.50 &  1.08e--03 & 3.62 &  3.00e--04 & 4.78 &  1.85e--05 & 4.07 &  5.70e--04 & \phantom{--}3.11 \\
& 2 & 2.08e--07 & 4.96 &  2.69e--05 & 5.03 &  5.05e--05 & 4.42 &  9.59e--06 & 4.97 &  7.29e--07 & 4.66 &  6.67e--05 & \phantom{--}3.10 \\
\bottomrule
\end{tabular}
\caption{$L^2(\Omega)$ discretization errors and estimated orders of convergence for the fully coupled model. On the $j$th refinement level, we used $2^{j+1} \times 2^j$ elements and time step~$\Delta \tilde{t} = 4\cdot10^{-5} \cdot 2^{-p \, (j+1)}$ for the subsurface problem and~$\Delta t = 4\cdot10^{-6} \cdot 2^{-p \, (j+1)}$ for the free-flow problem.}}
\label{tab:conv:coupled}
\end{table}

We choose computational domain~\revised{$\Omega(t)\cup\tO \subset\IR^2$ with~$\Omega(t) \coloneqq (0,100) {\,[\mathrm{m}]}\times(\zb,\xi(t)) {\,[\mathrm{m}]}$, $\tO \coloneqq (0,100) {\,[\mathrm{m}]}\times(-20,\zb) {\,[\mathrm{m}]}$, time interval~$J=(0,0.0002)$} , and a~sloped interface between free-flow and subsurface problem~$\zb(x^1) \coloneqq 0.005 x^1$, which gives us a~constant normal vector~$\vec{\nu} = \pm 1/\sqrt{1+0.005^2} \,\transpose{[-0.005, 1]}$.
For a~chosen free-surface elevation~$\xi$ and horizontal velocity~$u^{1}$ that fulfills the no-slip boundary condition at the bottom of the free-flow problem, one can derive matching analytical functions for~$\thh$ (using interface condition~\eqref{eq:coupling:pm}) and $u^{2}$ (using continuity equation~\eqref{eq:free:nonmixed:cont} and interface condition~\eqref{eq:coupling:ff}).
As mentioned before, here we use a~non-homogeneous boundary condition~\eqref{eq:free:nonmixed:bc:q} at the free surface, which gives us more freedom in our choice for~$u^{1}$, resulting in the following analytical solutions
\begin{align*}
\xi(t,x^1) &\coloneqq\; 
  \revised{5 + 0.003\,\sin(0.08 x^1 + t)} \,,\\
u^{1}(t,\vec{x}) &\coloneqq\; 
  y(t,x^1) \left( \cos(0.1 x^2) - \cos(0.1 \zb(x^1)) \right)\,,\\
u^{2}(t,\vec{x}) &\coloneqq\; 
  v(t,\vec{x}) + \varepsilon(t,x^1) \,,\\
\thh(t,\vec{x}) &\coloneqq\;
  \revised{\xi(t,x^1) + \left(\sin(0.1 x^2) - \sin(0.1 \zb(x^1))\right)}
\end{align*}
with diffusion coefficients~\revised{$\vecc{D} \coloneqq 0.001 \,\vecc{I}$}, $\vecc{\widetilde D} \coloneqq 0.01\, \vecc{I}$, and~$v(t,\vec{x})$ chosen such that~$\partial_{x^1} u^{1} + \partial_{x^2} v = 0$ in~$\Omega(t)$,
\begin{equation*}
v(t,\vec{x}) \,\coloneqq\; 
  -\partial_{x^1} y(t,x^1) \left( \frac{1}{0.1} \sin(0.1 x^2) - \cos(0.1\zb(x^1)) \, x^2 \right) - 
  0.1 \cdot 0.005 \cdot y(t,x^1) \, \sin(0.1 \zb(x^1)) \, x^2 \,,
\end{equation*}
and $\varepsilon(t,x^1)$ shifts~$u^{2}$ to fulfill coupling condition~\eqref{eq:coupling:ff}, i.e., 
\revised{
\begin{align*}
\varepsilon(t,x^1) \,\coloneqq\; &
-0.01 \cdot \partial_x \xi(t,x^1) -0.01 \cdot 0.1 \cdot (0.005^2 + 1) \cdot \cos(0.1 \zb(x^1)) \\
&+\partial_x y(t,x^1) \cdot \left(\frac{1}{0.1} \cdot \sin(0.1 \zb(x^1)) - \cos(0.1 \zb(x^1)) \cdot \zb(x^1)\right)\\
&+0.1 \cdot 0.005 \cdot y(t,x^1) \cdot \sin(0.1\zb(x^1)) \cdot \zb(x^1)
\end{align*}}
Function~$y(t,x^1)$ is used to increase the spatial variability in $x^1$-direction and to introduce a~time dependency.
Here, we use
\begin{equation*}
\revised{y(t,x^1)\,\coloneqq\; \sin(0.1 x^1 + t) \,.}
\end{equation*}
We prescribe Dirichlet boundary conditions for all unknowns and derive boundary data, right-hand side functions, and initial data from the analytical solutions.
Using this setup, we compute the solution for a~sequence of increasingly finer meshes with element sizes~$\Delta x_j$ and evaluate $L^2$-errors and estimated orders of convergence for any function~$c_{\Delta}$ as
\begin{equation*}
\mathrm{Err}(c) \,\coloneqq\; \|c_{\Delta x_{j-1}} - c\|_{L^2(\Omega)}\,,\qquad
\mathrm{EOC}(c) \,\coloneqq\; \ln \left(\frac{\|c_{\Delta x_{j-1}} - c\|_{L^2(\Omega)}}{\|c_{\Delta x_{j}} - c\|_{L^2(\Omega)}} \right)\Bigg/ \ln \left(\frac{\Delta x_{j-1}}{\Delta x_j}\right)
\end{equation*}
to compare these values to the analytically predicted ones.
We do this twice: first, for each sub-problem individually using analytically derived Dirichlet boundary conditions on the interior boundary (see Table~\ref{tab:conv:ffpm}), then for the fully coupled problem (see Table~\ref{tab:conv:coupled}). 

\subsection{Realistic example}
\label{sec:showcase}

\begin{figure}[!ht]
\centering
\includegraphics[width=\textwidth]{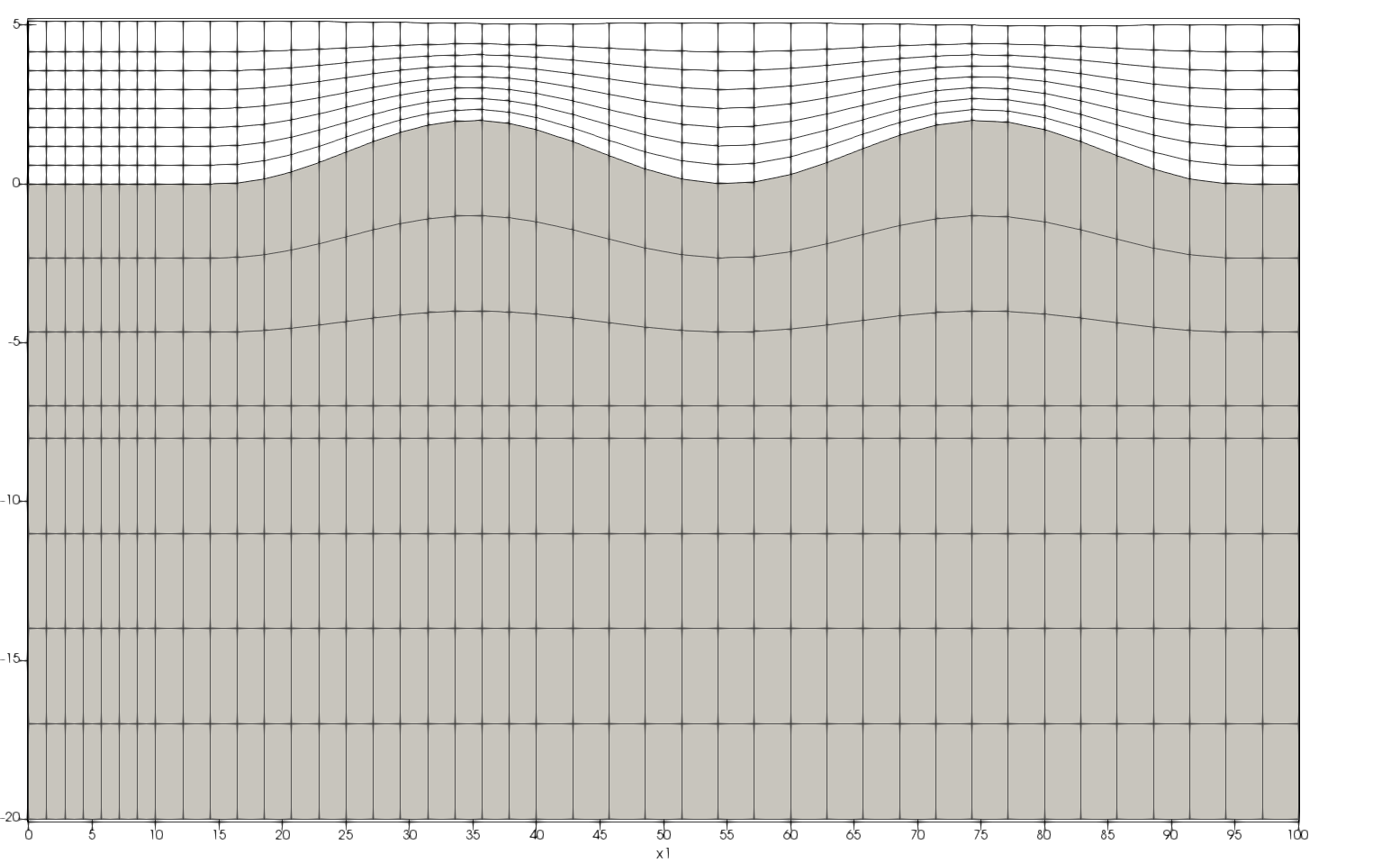}
\caption{Computational mesh at final state for the fully coupled example (cf.~Sec.~\ref{sec:showcase}), with free flow domain (white background) and subsurface domain (gray background). Illustration is compressed in the horizontal direction.}
\label{fig:showcase:grid}
\end{figure}

\begin{figure}[!ht]
\includegraphics[width=\textwidth]{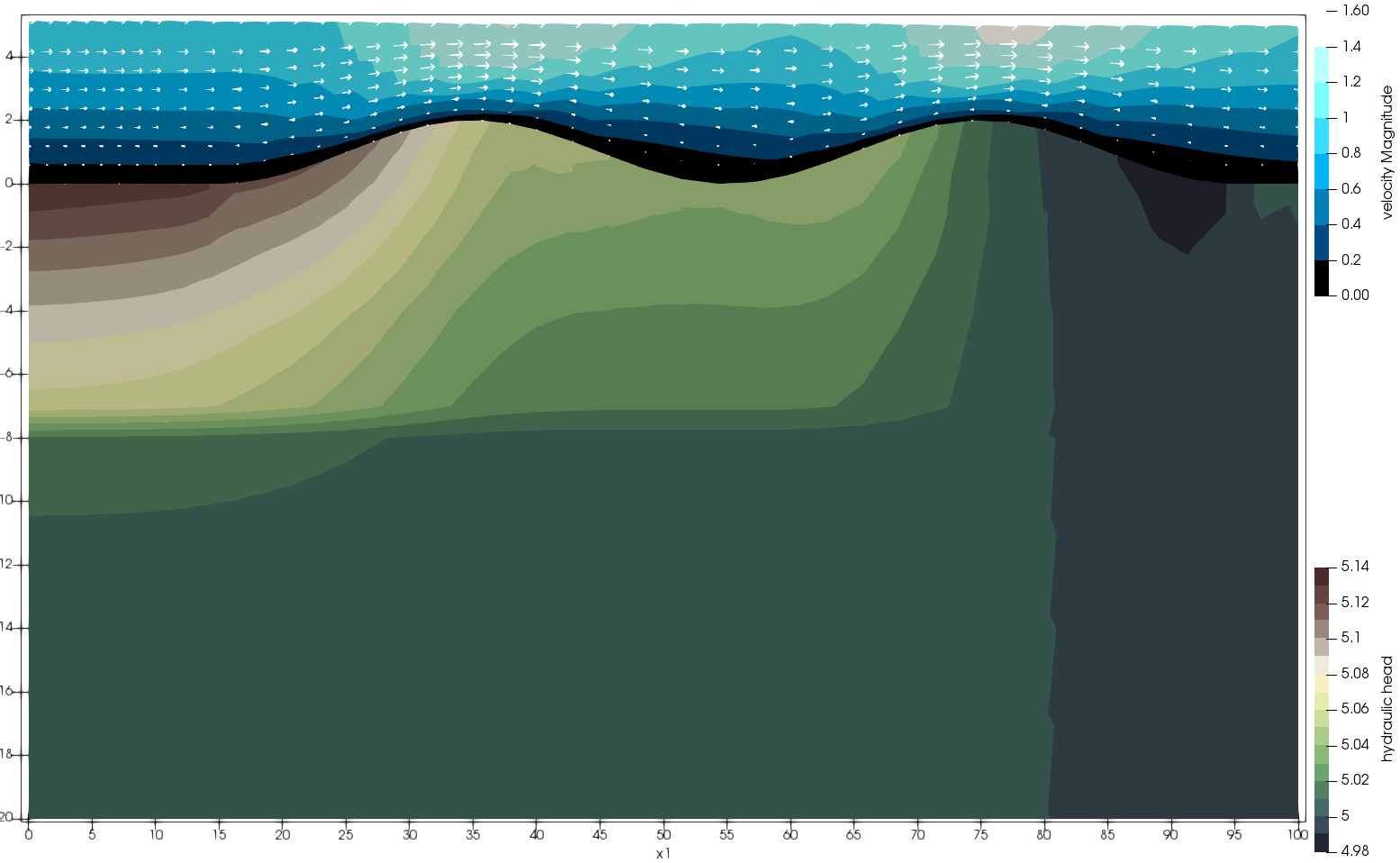}
\caption{Velocity magnitude in the free-flow domain (arrows indicating direction of the velocity field, scaled by magnitude) and hydraulic head at final state of the fully coupled example (cf.~Sec.~\ref{sec:showcase}). Illustration is compressed in the horizontal direction.}
\label{fig:showcase:vel-head}
\end{figure}

\begin{figure}[!ht]
\centering
\includegraphics[width=\textwidth]{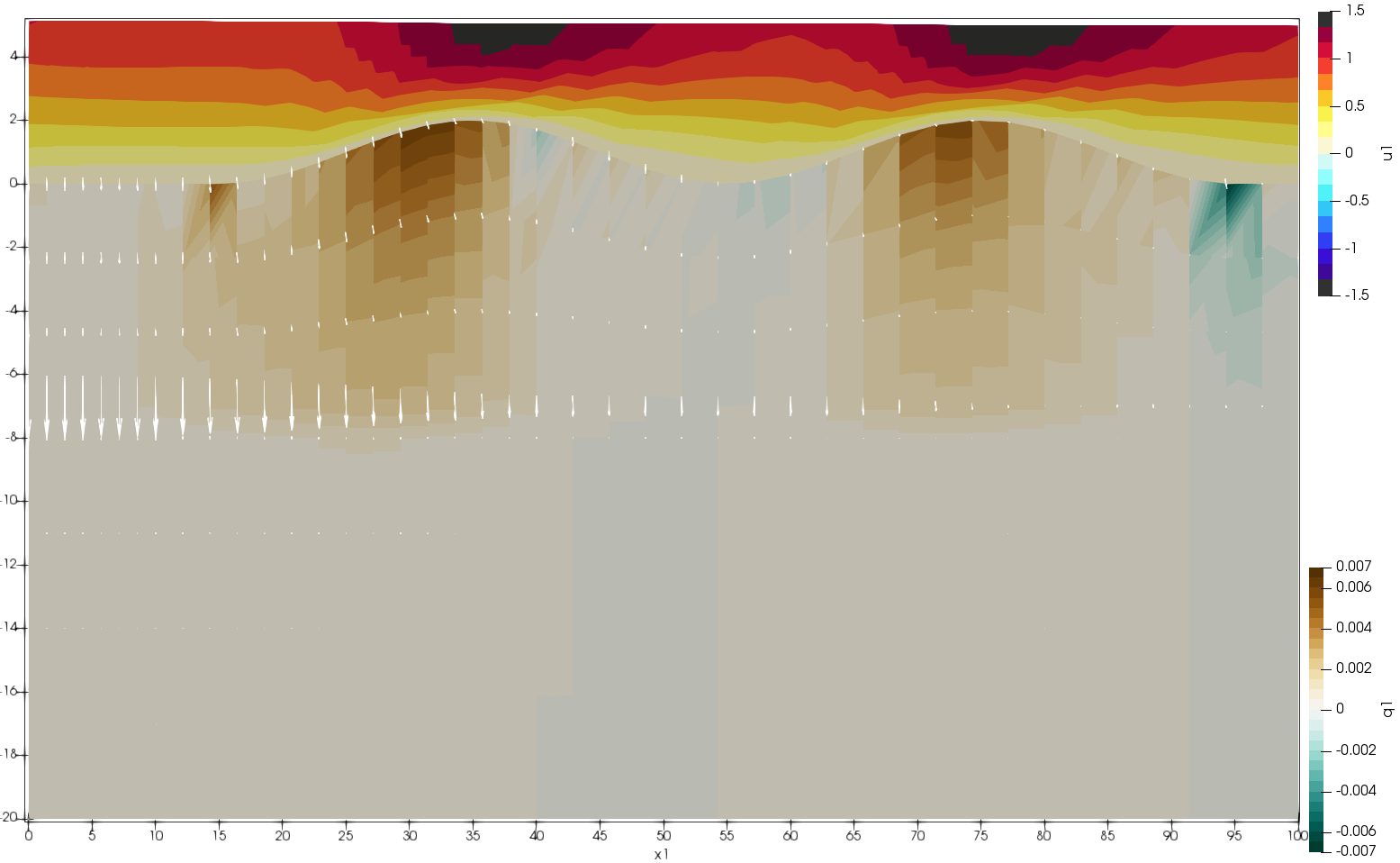}
\caption{$x^1$ components of velocity~$\vec{u}$ and flux~$\tilde{\vec{q}}$ at final state of the fully coupled example (cf.~Sec.~\ref{sec:showcase}). Arrows in the subsurface domain indicate direction and magnitude of the flux field with the magnitude~40~times amplified compared to~Fig.~\ref{fig:showcase:vel-head}. Illustration is compressed in the horizontal direction.}
\label{fig:showcase:u1-q1}
\end{figure}

\begin{figure}[!ht]
\includegraphics[width=\textwidth]{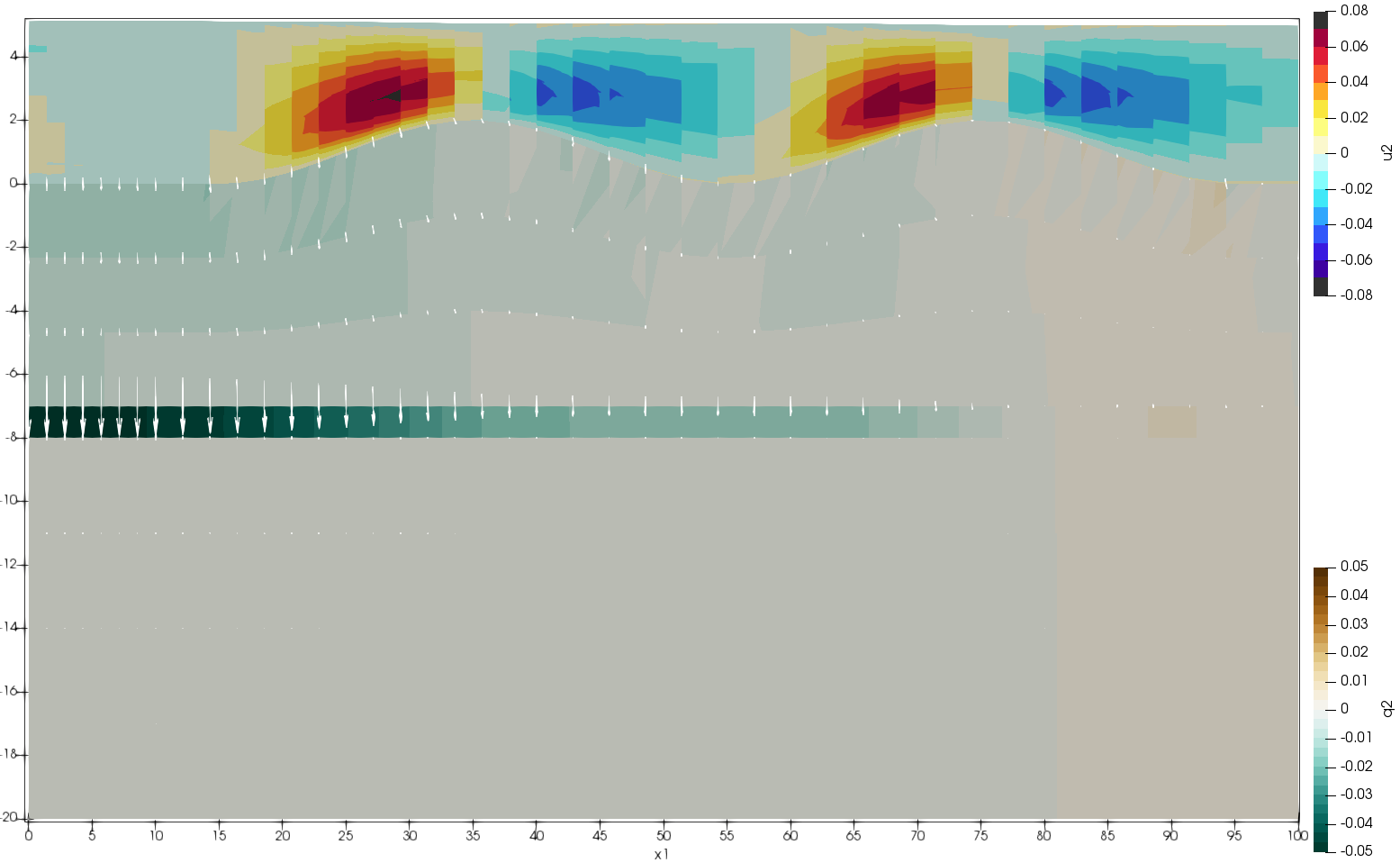}
\caption{$x^2$ components of velocity~$\vec{u}$ and flux~$\tilde{\vec{q}}$ at final state of the fully coupled example (cf.~Sec.~\ref{sec:showcase}). Arrows in the subsurface domain indicate direction and magnitude of the flux field with the magnitude~40~times amplified compared to~Fig.~\ref{fig:showcase:vel-head}. Illustration is compressed in the horizontal direction.}
\label{fig:showcase:u2-q2}
\end{figure}%

\begin{figure}[!ht]
\centering
\includegraphics[width=.95\textwidth]{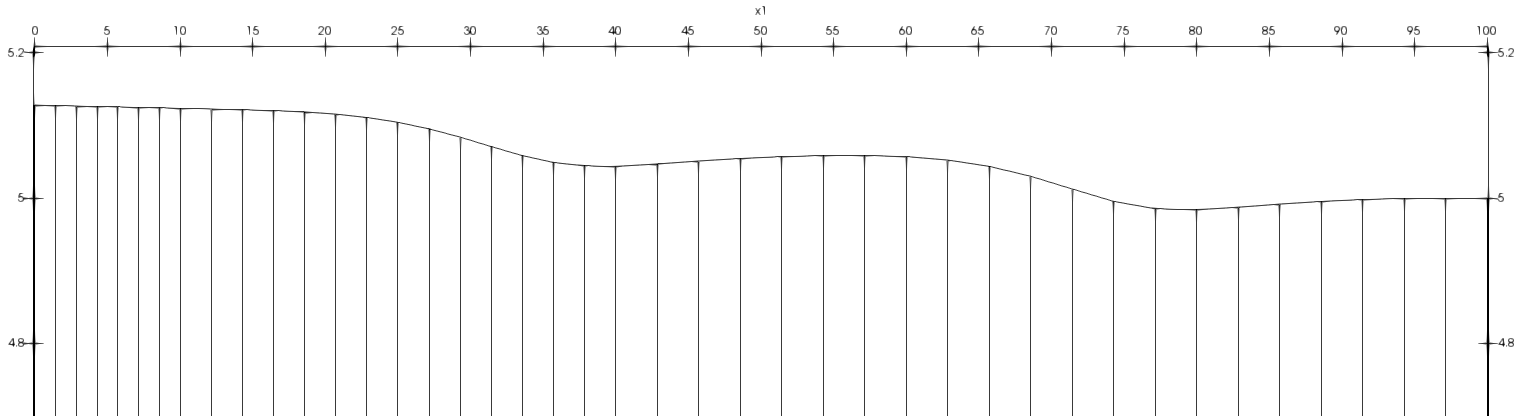}
\caption{Detail of the computational mesh at the free surface for the fully coupled example (cf.~Sec.~\ref{sec:showcase}) at final state. Illustration is compressed in the horizontal direction.}
\label{fig:showcase:grid-surface}
\end{figure}

To examine a~more realistic example, we set up a~modified channel flow with free flow domain~$\Omega(t)\coloneqq(0,100){\,[\mathrm{m}]}\times(\zb(x^1),\xi(t,x^1)) {\,[\mathrm{m}]}$ and subsurface domain~$\widetilde{\Omega}\coloneqq(0,100){\,[\mathrm{m}]}\times(-20,\zb(x^1)){\,[\mathrm{m}]}$, where the topography of the interior boundary is given as
\begin{equation*}
\zb(x^1) \coloneqq \begin{cases}
  \cos\left(\frac{x^1-35}{20}\,\pi\right) + 1\,, & \text{if } 15 \leq x^1 \leq 95 \\
  0\,, & \text{otherwise}
\end{cases}
\,[\mathrm{m}]\,.
\end{equation*}
We apply a~Neumann boundary condition~\eqref{eq:sub:N} with~$\tilde{g}_\mathrm{N}=0$ to left, right, and bottom boundary of the subsurface domain.
On the left boundary of the free flow problem, we apply a~river boundary condition (cf.~Sec.~\ref{sec:model:ff}) with~$h_\mathrm{D} = 5\,[\mathrm{m}]$, $u^1_\mathrm{D} = \ln(1 + (e - 1) x^2 / 5)\,[\frac{\mathrm{m}}{\mathrm{s}}]$ and an~open sea boundary condition with~$h_\mathrm{D}=5\,[\mathrm{m}]$ on the right. 
In the subsurface domain we choose a~homogeneous diffusion coefficient and insert a~thin layer at $x^2\in[-8,-7]$\,{\scriptsize[m]} where the diffusion coefficient is one order of magnitude smaller.
The mesh is chosen such that element boundaries coincide with the jump in the diffusion coefficient (see Fig.~\ref{fig:showcase:grid}).
Remaining parameters and initial conditions are chosen to be
\begin{align*}
  &h_0(x^1) = 5 - \zb(x^1) \, [\mathrm{m}] \,, &
  &u^1_0(\vec{x}) = 0 \,  {\textstyle\left[\frac{\mathrm{m}}{\mathrm{s}}\right]} \,, &
  &\vecc{D} = \begin{bmatrix} 0 & 0 \\ 0 & 0.08 \end{bmatrix} \, {\textstyle\left[\mathrm{\frac{m^2}{s^2}} \right]} \,, 
  \\
  &g = 10 \, {\textstyle\left[\mathrm{\frac{m}{s^2}}\right]}\,,&
  &\thh_0(\vec{x}) = 5 \,[\mathrm{m}] \,, & 
  &\vecc{\widetilde D} = \begin{cases}
    10^{-4} \cdot \vecc{I} 
    & \text{if}~~-8 \leq x^2 \leq -7\\
    10^{-3} \cdot \vecc{I} 
    & \text{otherwise}
  \end{cases} \,{\textstyle\left[\mathrm{\frac{m^2}{s}}\right]} \,,
\end{align*}
where~$\vecc{I}$ denotes the $2\times2$ identity matrix.
We run the simulation with~$42\times 8$ elements in each domain (cf.~Fig.~\ref{fig:showcase:grid}), polynomial approximation order~$p=1$, and time step size~$\Delta \tilde{t} = 5 \cdot \Delta t = 0.1 \,[\mathrm{s}]$ until~$t_\mathrm{end}=30\,000\,[\mathrm{s}]$.

Results are depicted in Figs.~\ref{fig:showcase:vel-head}--\ref{fig:showcase:grid-surface}.
The flow field in the free-flow domain exhibits the typical logarithmic flow profile in the horizontal velocity component with a~compression due to the obstacle on the bottom boundary leading to \enquote{updrafts} before and \enquote{downdrafts} after the obstacles with corresponding bumps and depressions in the free-surface elevation (cf.~Fig.~\ref{fig:showcase:grid-surface}).
The differences in elevation are primarily responsible for the induced gradient in the hydraulic head.
The flux components in the subsurface problem plotted in Fig.~\ref{fig:showcase:u1-q1} exhibit some slight oscillations since no slope limiting has been used in this test problem. 
However, the local maximum at $x^1\approx15$ and the local minimum at $x^1\approx95$ are not caused by these oscillations but rather are the effect of the reversion in the flux direction at the interface between subdomains at the corresponding locations.
We also observed similar phenomena when simulating this test case using piecewise constant approximations for all unknowns.

% ...

\section{Conclusion and Outlook}
This fourth installment in our paper series on implementing a~\MatOct{} toolbox introduces a~highly flexible problem implementation framework that is suitable for most problem classes and accommodates coupled multi-physics simulations.
We apply this framework to coupled free-surface/subsurface flow simulations in a~vertical slice and demonstrate the performance of the solver using analytical and realistic tests.
The models make use of the newly added support for quadrilateral meshes.
Plans for future work include implementing further physical models and related numerical schemes using our framework, in particular, nonlinear operators.

\small

\section*{Acknowledgments}
B.~Reuter would like to thank the German Research Foundation (DFG) for financial support under grant AI~117/1-1 \enquote{Modeling of ocean overflows using statically and dynamically adaptive vertical discretization techniques}.
A.~Rupp acknowledges financial support by the DFG~EXC~2181 \enquote{STRUCTURES: A~unifying approach to emergent phenomena in the physical world, mathematics, and complex data} and the DFG~RU~2179 \enquote{MAD~Soil -- Microaggregates: Formation and turnover of the structural building blocks of soils}.

\appendix

\section{Explicit form of matrices and vectors in the free-flow problem}\label{sec:free-flow:mat-vec}

In the interest of a~compact presentation, we left out the explicit forms of matrices and vectors in the presentation of the linear system of equations~\eqref{eq:free:matsys} for the free-flow problem.
With the corresponding matrix and vector names given for each term in system~\eqref{eq:free:sys} the explicit forms can easily be derived using the same steps as in our previous publications~\cite{FESTUNG1,FESTUNG2}.
For convenience, we provide the final forms in the following together with implementation remarks in Appendix~\ref{sec:assembly}.

\subsection[Contributions from area terms]{Contributions from area terms~$\I$--$\IV$, $\XI$--$\XIII$, $\XVI$, $\XXI$--$\XXII$}
The matrices in the remainder of this section have sparse block structure; by giving definitions for non-zero blocks we tacitly assume a~zero fill-in.
Consider system~\eqref{eq:free:sys}.
The mass matrix $\vecc{M}\in\IR^{KN\times KN}$ in terms $\I$ and $\XII$ is identical to the one of the subsurface problem and is defined component-wise as
\begin{equation*}
[\vecc{M}]_{(k-1)N+i, (k-1)N+j} \;\coloneqq\; \int_{T_k}\vphi_{ki}\,\vphi_{kj}\,\dd\vec{x}\;.
\end{equation*}
Since the basis functions~$\vphi_{ki}$, $i\in\{1,\ldots,N\}$ are supported only on~$T_k$, $\vecc{M}$ has a~block-diagonal structure 
\begin{equation}\label{eq:globMlocM}
\vecc{M} \;=\;
\begin{bmatrix}
\vecc{M}_{T_1} &          & \\
               & ~\ddots~ & \\
               &          & \vecc{M}_{T_K}
\end{bmatrix}
\qquad\text{with}\qquad
\vecc{M}_{T_k} \;\coloneqq\;
\int_{T_k}\begin{bmatrix}
\vphi_{k1}\,\vphi_{k1} & \cdots & \vphi_{k1}\,\vphi_{kN} ~\\
 \vdots                & \ddots & \vdots \\
\vphi_{kN}\,\vphi_{k1} & \cdots & \vphi_{kN}\,\vphi_{kN} 
\end{bmatrix}\,\dd\vec{x}\;,
\end{equation}
i.\,e., it consists of $K$~\emph{local mass matrices}~$\vecc{M}_{T_k}\in\IR^{N\times N}$.  
Henceforth, we write~$\vecc{M} = \diag \big(\vecc{M}_{T_1},\ldots,\vecc{M}_{T_K}\big)$.
The block-diagonal matrices~$\vecc{G}^m\in\IR^{KN\times KN}, \vecc{H}^m\in\IR^{KN\times KN}, \;m\in\{1,2\}$ in terms~$\IV$,~$\XIII$, and~$\XVI$ are the same as~$\tilde{\vecc{G}}^m,\tilde{\vecc{H}}^m$ in the sub\-surface problem (except for the name change of the diffusion matrix from~$\vecc{\widetilde D}$ to $\vecc{D}$) and are given by
\begin{align}
[\vecc{G}^m]_{(k-1)N+i, (k-1)N+j} &\;\coloneqq\; 
 \sum_{r=1}^N \sum_{l=1}^N D^{rm}_{kl}(t) \int_{T_k} \partial_{x^r} \vphi_{ki}\,\vphi_{kl} \,\vphi_{kj} \,\dd\vec{x}\,,
 \label{eq:globGm}\\
[\vecc{H}^m]_{(k-1)N+i, (k-1)N+j} &\;\coloneqq\; \int_{T_k} \partial_{x^m} \vphi_{ki}\,\vphi_{kj}\,\dd\vec{x}\,.
\label{eq:globHm}
\end{align}
For general remarks and detailed presentation regarding the assembly of these three matrices, we refer to \ref{sec:assembly:elem} and to our first paper~\cite{FESTUNG1}.
Term $\II$ defines two block-diagonal matrices~$\vecc{E}^m\in\IR^{KN\times KN}, m\in\{1,2\}$ given by
\begin{equation*}
\left[ \vecc{E}^m \right]_{(k-1)N+i,(k-1)N+j} \;\coloneqq\;
  \sum_{l=1}^N U_{kl}^1 
  \int_{T_k} \partial_{x^m} \vphi_{ki} \, \vphi_{kl} \, \vphi_{kj} \, \dd\vec{x}\,.
\end{equation*}
Their entries are very similar to those of~$\vecc{G}^m$ the only difference being a~scalar-valued coefficient function (here, $u^1_{\Delta}$) instead of a~matrix-valued one ($\vecc{D}$ in $\vecc{G}^m$).
Thus, their assembly takes the same form as given in~\eqref{eq:globGm} but without the sum over~$r$ and without~$r$ and~$m$ in the coefficient.

Term $\III$ contributes to a~rectangular matrix~$\check{\vecc{H}} \in \IR^{KN\times \overline{K}\overline{N}}$ that is applied to the representation vector of the one-dimensional water height~$\vec{H}$.
It has one non-zero ($N\times\overline{N}$-dimensional) block per block-row but possibly multiple non-zero blocks per block-column, since more than one two-dimensional element might correspond to each one-dimensional element~$\overline{T}_{\overline{k}}$. 
Its entries take a~similar form as for $\vecc{H}^m$ and are given by
\begin{equation*}
\left[\check{\vecc{H}} \right]_{(k-1)N+i,(\overline{k}-1)\overline{N}+j} \coloneqq
\int_{T_k} \partial_{x^1} \vphi_{ki} \, \phi_{\overline{k}j} \, \dd\vec{x}\,.
\end{equation*}
The right-hand side vector~$\vec{L}_u$ in term~$\XI$ is computed from the right-hand side function~$f_{\Delta}$ as in the subsurface problem.
Additionally, a~vector~$\vec{L}_{\zb}$ is assembled with the gradient of the bathymetry~$\zb$. 
For any polynomial approximation order, $\zb$ is represented by a~continuous piecewise linear function interpolating the topographic heights specified at the nodes of the one-dimensional mesh~$\Proj \setT_{\Delta}$, and thus, its $x^1$-derivative is an~element-local constant.
In free-surface equation~\eqref{eq:free:sys:surf}, term $\XXI$ contributes the one-dimensional mass matrix $\overline{\vecc{M}} \in \IR^{\overline{K}\overline{N}\times\overline{K}\overline{N}}$ defined component-wise as
\begin{equation*}
\left[\overline{\vecc{M}}\right]_{(\overline{k}-1)\overline{N}+i, (\overline{k}-1)\overline{N}+j} \;\coloneqq\; 
\int_{\overline{T}_{\overline{k}}}\phi_{\overline{k}i}\,\phi_{\overline{k}j}\,\dd x^1\;.
\end{equation*}
Term $\XXII$ contributes to a~block-diagonal matrix~$\overline{\vecc{G}} \in \IR^{\overline{K}\overline{N}\times KN}$ with entries given component-wise by
\begin{equation}
\label{eq:barG}
\left[\overline{\vecc{G}}\right]_{(\overline{k}-1)\overline{N}+i,(\overline{k}-1)\overline{N}+j} \coloneqq
\int_{\overline{T}_{\overline{k}}} \frac{1}{\Hs} \partial_{x^1} \phi_{\overline{k}i}\, \left( \sum_{l=1}^{\overline{N}} \overline{U}_{\overline{k}l} \, \phi_{\overline{k}l} \right)
\, \phi_{\overline{k}j} \, \dd x^1\,.
\end{equation}
Its assembly is detailed in~\ref{sec:assembly:elem}.

\subsection[Contributions from interior edge terms]{Contributions from interior edge terms~$\V$, $\VIII$, $\XV$, $\XVII$--$\XVIII$, $\XXIII$}
In this section, we consider a~fixed element $T_k=T_{k^-}$ with an~interior edge~$E_{k^-n^-} \in\partial T_{k^-}\cap\setE_\Omega = \partial T_{k^-} \cap \partial T_{k^+}$ shared by an adjacent element $T_{k^+}$ and associated with fixed local edge indices $n^-,n^+ \in \{1,2,3,4\}$.
Due to the numbering of our structured grid, we have a~unique mapping between $n^-$ and $n^+$, thus $n^+$ is implicitly given by $n^-$ (see Sec.~\ref{sec:gridbasis}).

Terms $\V$ and $\XV$ stem from integration by parts and are the edge integral counterparts to terms $\IV$ and $\XIII$, respectively.
Block-matrices~$\vecc{R}^m, \vecc{Q}^m \in \IR^{KN\times KN}$, $m\in\{1,2\}$ resulting from these terms have the same form as~$\tilde{\vecc{R}}^m,\tilde{\vecc{Q}}^m$ in the subsurface problem.
Diagonal blocks of~$\vecc{Q}^m$ are given by 
\begin{subequations}\label{eq:globQ}
\begin{equation}\label{eq:globQ:diag}
[\vecc{Q}^m]_{(k-1)N+i,(k-1)N+j} \;\coloneqq\;
\frac{1}{2} \sum_{E_{kn}\in\partial T_k\cap\setE_\Omega} \nu_{kn}^m \int_{E_{kn}} \vphi_{ki}\,\vphi_{kj}\,\dd\sigma \;.
\end{equation}
Entries in off-diagonal blocks in $\vecc{Q}^m$ are only non-zero for pairs of elements $T_{k^-}$, $T_{k^+}$ with $\partial T_{k^-}\cap\partial T_{k^+}\neq\emptyset$.
They consist of the terms containing basis functions from both adjacent elements and are given as
\begin{equation}\label{eq:globQ:offdiag}
[\vecc{Q}^m]_{({k^-}-1)N+i,({k^+}-1)N+j} \;\coloneqq\;
\frac{1}{2} \nu_{k^-n^-}^m \int_{E_{k^-n^-}} \vphi_{k^-i}\,\vphi_{k^+j}\,\dd\sigma\;.
\end{equation}
\end{subequations}
Entries in diagonal and off-diagonal blocks of~$\vecc{R}^m$ are given component-wise by
\begin{subequations}\label{eq:globR}
\begin{align}\label{eq:globR:diag}
&[\vecc{R}^m]_{(k-1)N+i,(k-1)N+j} \;\coloneqq\;
\frac{1}{2} \sum_{E_{kn}\in\partial T_k\cap\setE_\Omega} \sum_{r=1}^2 \nu_{kn}^r \sum_{l=1}^N D^{rm}_{kl} \int_{E_{kn}} \vphi_{ki}\,\vphi_{kl}\,\vphi_{kj}\,\dd\sigma\;,\\
&[\vecc{R}^m]_{({k^-}-1)N+i,({k^+}-1)N+j} \;\coloneqq\;
\frac{1}{2} \sum_{r=1}^2 \nu_{k^-n^-}^r \sum_{l=1}^N D^{rm}_{k^+l} \int_{E_{k^-n^-}} \vphi_{k^-i}\,\vphi_{k^+l}\,\vphi_{k^+j}\,\dd\sigma\;.
\end{align}
\end{subequations}
All off-diagonal blocks corresponding to pairs of elements not sharing an edge are zero. 
More involved are the edge integrals in terms $\VIII$, $\XVII$--$\XVIII$, and $\XXIII$ that correspond to nonlinear fluxes.

\paragraph{Interior edges $\setE_\Omega$}
We consider a~fixed element $T_k=T_{k^-}$ with an~interior edge~$E_{k^-n^-} \in\partial T_{k^-}\cap\setE_\Omega= \partial T_{k^-} \cap \partial T_{k^+}$ shared by an adjacent element $T_{k^+}$ and associated with fixed local edge indices $n^-,n^+ \in \{1,2,3,4\}$.
Due to the numbering of our structured grid we have a fixed mapping between $n^-$ and $n^+$, thus $n^+$ is implicitly given by $n^-$ (see Sec.~\ref{sec:gridbasis}).
Additionally, we have the corresponding one-dimensional elements~$\overline{T}_{\overline{k}} = \overline{T}_{\overline{k}^-} = \Proj T_k$ and~$\overline{T}_{\overline{k}^+} = \Proj T_{k^+}$.
From term~$\VIII$, we have contributions for~$\vphi_{k^-i}$ in two block-matrices~$\vecc{P}^{m} \in \IR^{KN\times KN}$ of the form
\begin{align*}
&\frac{1}{2} \sum_{m=1}^2 \nu_{k^-n^-}^m  \sum_{j=1}^N U_{k^-j}^m \sum_{l=1}^N U_{k^-l}^1 
  \int_{E_{k^-n^-}} \vphi_{k^-i} \, \vphi_{k^-l} \, \vphi_{k^-j} \, \dd \sigma \\
&\qquad+
\frac{1}{2} \sum_{m=1}^2 \nu_{k^-n^-}^m  \sum_{j=1}^N U_{k^+j}^m \sum_{l=1}^N U_{k^+l}^1 
  \int_{E_{k^-n^-}} \vphi_{k^-i} \, \vphi_{k^+l} \, \vphi_{k^+j} \, \dd \sigma\,.
\end{align*}
Entries in diagonal blocks of $\vecc{P}^{m}$ are then given component-wise by
\begin{subequations}\label{eq:free:Ph}
\begin{equation}
\left[ \vecc{P}^{m} \right]_{(k-1)N+i,(k-1)N+j} \coloneqq 
\frac{1}{2} \sum_{E_{kn}\in\partial T_k\cap\setE_\Omega} 
\nu_{kn}^m \sum_{l=1}^N U_{kl}^1 
\int_{E_{kn}} \vphi_{ki} \, \vphi_{kl} \, \vphi_{kj} \, \dd \sigma \,,
\end{equation}
and entries in off-diagonal blocks consist of
\begin{equation}
\left[ \vecc{P}^{m} \right]_{(k^--1)N+i,(k^+-1)N+j} \coloneqq 
\frac{1}{2}
\nu_{k^-n^-}^m \sum_{l=1}^N U_{k^+l}^1
\int_{E_{k^-n^-}} \vphi_{k^-i} \, \vphi_{k^+l} \, \vphi_{k^+j} \, \dd \sigma \,.
\end{equation}
\end{subequations}
Additionally, we have a~contribution for~$\vphi_{k^-i}$ from term $\VIII$ to a~rectangular matrix~$\check{\vecc{Q}}\in\IR^{KN\times \overline{K}\overline{N}}$
\begin{equation*}
\frac{1}{2} \nu_{k^-n^-}^1 \sum_{j=1}^{\overline{N}} H_{\overline{k}^-j}
  \int_{E_{k^-n^-}} \vphi_{k^-i} \, \phi_{\overline{k}^-j} \, \dd \sigma +
\frac{1}{2} \nu_{k^-n^-}^1 \sum_{j=1}^{\overline{N}} H_{\overline{k}^+j}
  \int_{E_{k^-n^-}} \vphi_{k^-i} \, \phi_{\overline{k}^+j} \, \dd \sigma \,,
\end{equation*}
which results in two blocks
\begin{equation}\begin{aligned}
\label{eq:free:tildeQh}
&\left[ \check{\vecc{Q}} \right]_{(k-1)N+i,(\overline{k}-1)\overline{N}+j} 
\coloneqq \sum_{E_{kn}\in\partial T_k\cap\setE_\Omega} \frac{\nu_{kn}^1}{2}
\int_{E_{kn}} \vphi_{ki} \, \phi_{\overline{k}j} \,\dd\sigma \,, \\
&\left[ \check{\vecc{Q}} \right]_{(k^--1)N+i,(\overline{k}^+-1)\overline{N}+j} 
\coloneqq \frac{\nu_{k^-n^-}^1}{2}
\int_{E_{k^-n^-}} \vphi_{k^-i} \, \phi_{\overline{k}^+j} \,\dd\sigma \,.
\end{aligned}\end{equation}

Term~$\XVI$ contributes to block-matrices~$\vecc{Q}_\mathrm{avg}, \vecc{Q}_\mathrm{up} \in \IR^{KN\times KN}$.
The first is equivalent to~$\vecc{Q}^1$, restricted to contributions from horizontal edges instead of all interior edges (and thus nonzero only in the case of horizontal edges that are not aligned with the $x^1$-axis).
The latter is similar to $\vecc{Q}^2$, built from its diagonal blocks for $n^- = 2$ and off-diagonal blocks for $n^- = 1$ and without the factor of~1/2.

\paragraph{Interior vertical edges $\setE_\Omega^\mathrm{v}$}
Due to the Lax--Friedrichs Riemann solver (cf.~Sec.~\ref{sec:riemann}) on vertical edges, we have some additional terms exclusively on vertical edges.
As before, we consider a~fixed element $T_k=T_{k^-}$ with an~interior edge~$E_{k^-n^-} \in \partial T_{k^-} \cap \partial T_{k^+} = \partial T_{k^-} \cap \setE_\Omega^\mathrm{v}$ shared by an adjacent element $T_{k^+}$ and associated with fixed local edge indices $n^-,n^+ \in \{3,4\}$ and corresponding one-dimensional elements~$\overline{T}_{\overline{k}} = \overline{T}_{\overline{k}^-} = \Proj T_k$ and~$\overline{T}_{\overline{k}^+} = \Proj T_{k^+}$.

From term~$\VIII$, we have the additional contribution to a~vector~$\vec{K}_u \in \IR^{KN}$ with entries of the form
\begin{equation}
\label{eq:free:Ku}
\left[\vec{K}_u\right]_{(k^--1)N+i} \coloneqq
\sum_{E_{k^-n^-}\in\partial T_k^- \cap \setE_\Omega^\mathrm{v}}
\int_{E_{k^-n^-}} \vphi_{k^-i} \, \frac{|\hat{\lambda}|}{2} \, \left(
\sum_{j=1}^N U_{k^-j}^1 \, \vphi_{k^-j} -
\sum_{j=1}^N U_{k^+j}^1 \, \vphi_{k^+j} \right)\, \dd\sigma \,.
\end{equation}
For assembly details, see~\ref{sec:assembly:jump}.
Term $\XVII$ has a~contribution for~$\vphi_{ki}$ to matrix~$\check{\vecc{P}} \in \IR^{KN\times KN}$ of the form
\begin{equation*}
\frac{1}{2} \nu_{k^-n^-}^1 \sum_{j=1}^N U_{k^-j}^1 
\sum_{l=1}^{\overline{N}} H_{\overline{k}^-l}
\int_{E_{k^-n^-}} \frac{1}{\Hs} \vphi_{k^-i} \, \phi_{\overline{k}^-l} \, \vphi_{k^-j} \, \dd\sigma
+
\frac{1}{2} \nu_{k^-n^-}^1 \sum_{j=1}^N U_{k^+j}^1 
\sum_{l=1}^{\overline{N}} H_{\overline{k}^+l}
\int_{E_{k^-n^-}} \frac{1}{\Hs} \vphi_{k^-i} \, \phi_{\overline{k}^+l} \, \vphi_{k^+j} \, \dd\sigma
\end{equation*}
with entries in diagonal blocks given component-wise by
\begin{subequations}\label{eq:free:tildeP}
\begin{equation}\label{eq:free:tildeP:diag}
\left[\check{\vecc{P}}\right]_{(k-1)N+i,(k-1)N+j} \coloneqq
\frac{1}{2} \sum_{E_{kn}\in\partial T_k \cap \setE_\Omega^\mathrm{v}} \nu_{kn}^1
\sum_{l=1}^{\overline{N}} H_{\overline{k}l}
\int_{E_{kn}} \frac{1}{\Hs} \vphi_{ki} \, \phi_{\overline{k}l} \, \vphi_{kj} \, \dd\sigma 
\end{equation}
and in off-diagonal blocks by
\begin{equation}
\left[\check{\vecc{P}}\right]_{(k^--1)N+i,(k^--1)N+j} \coloneqq
\frac{1}{2} \nu_{k^-n^-}^1
\sum_{l=1}^{\overline{N}} H_{\overline{k}^+l}
\int_{E_{k^-n^-}} \frac{1}{\Hs} \vphi_{k^-i} \, \phi_{\overline{k}^+l} \, \vphi_{k^+j} \, \dd\sigma \,.
\end{equation}
\end{subequations}
Additionally, the same term contributes to a~vector~$\vec{K}_{h} \in \IR^{KN}$ in the form
\begin{equation}
\label{eq:free:Kh}
\left[\vec{K}_{h}\right]_{(k^--1)N+i} \coloneqq
\sum_{E_{k^-n^-}\in\partial T_k^- \cap \setE_\Omega^\mathrm{v}}
\int_{E_{k^-n^-}} \frac{1}{\Hs} \vphi_{k^-i} \, \frac{|\hat{\lambda}|}{2} \, \left(
  \sum_{j=1}^{\overline{N}} H_{\overline{k}^-j} \, \phi_{\overline{k}^-j} -
  \sum_{j=1}^{\overline{N}} H_{\overline{k}^+j} \, \phi_{\overline{k}^+j} 
\right) \, \dd\sigma \,,
\end{equation}
which is assembled in the same way as~$\vec{K}_u$.
In free-surface equation~\eqref{eq:free:sys:surf}, we have contributions only for vertical and boundary edges.
Thus, we can again use the depth-averaged representation of the horizontal velocity components, and term $\XXIII$ contributes to a~matrix~$\overline{\vecc{P}} \in \IR^{\overline{K}\overline{N}\times\overline{K}\overline{N}}$ for $\phi_{\overline{k}i}$
\begin{align*}
&\frac{1}{2} \sum_{\Proj E_{k^-n^-} = \overline{E}_{\overline{k}^-\overline{n}^-}}
\nu_{k^-n^-}^1 
\int_{E_{k^-n^-}} \frac{1}{\Hs} \phi_{\overline{k}^-i} \, \left(
  \sum_{j=1}^{\overline{N}} H_{\overline{k}^-j} \, \phi_{\overline{k}^-j} \, 
  \sum_{l=1}^N U_{k^-l}^1 \, \vphi_{k^-l} +
  \sum_{j=1}^{\overline{N}} H_{\overline{k}^+j} \, \phi_{\overline{k}^+j} \, 
  \sum_{l=1}^N U_{k^+l}^1 \, \vphi_{k^+l} 
\right) \,\dd\sigma = \\
& \sum_{a^1_{\overline{k}^-n^-}}
\frac{\nu_{\overline{k}^-\overline{n}^-}}{2 \Hs} \left(
  \sum_{j=1}^{\overline{N}} H_{\overline{k}^-j} \,
  \sum_{l=1}^{\overline{N}} \overline{U}_{\overline{k}^-l} \,
  \phi_{\overline{k}^-i} \phi_{\overline{k}^-j} \phi_{\overline{k}^-l} 
\,+\,
  \sum_{j=1}^{\overline{N}} H_{\overline{k}^+j} \,
  \sum_{l=1}^{\overline{N}} \overline{U}_{\overline{k}^+l} \,
  \phi_{\overline{k}^-i}  \phi_{\overline{k}^+j} \phi_{\overline{k}^+l}
\right) \,,
\end{align*}
where basis functions, coefficients of the depth-integrated velocity, and smoothed height~$\Hs$ are evaluated in the endpoints $a^1_{\overline{k}\overline{n}}$,~$\overline{n}\in\{1,2\}$ of each one-dimensional element~$\overline{T}_{\overline{k}}$.
The resulting diagonal blocks read as
\begin{subequations}\label{eq:free:barP}
\begin{equation}\label{eq:free:barP:diag}
\left[\overline{\vecc{P}}\right]_{(\overline{k}-1)\overline{N}+i,(\overline{k}-1)\overline{N}+j} \coloneqq
\frac{1}{2} \sum_{a^1_{\overline{k}\overline{n}}\in\partial \overline{T}_{\overline{k}}\cap \overline{\setV}_\Omega}
\frac{\nu_{\overline{k}\overline{n}}^1}{\Hs\left(a^1_{\overline{k}\overline{n}}\right)}
\sum_{l=1}^{\overline{N}} \overline{U}_{\overline{k}l}\left(a^1_{\overline{k}\overline{n}}\right)\;
\phi_{\overline{k}i} \left(a^1_{\overline{k}\overline{n}}\right) \, \phi_{\overline{k}l}\left(a^1_{\overline{k}\overline{n}}\right) \, \phi_{\overline{k}j}\left(a^1_{\overline{k}\overline{n}}\right) 
\end{equation}
and off-diagonal blocks as
\begin{equation}
\left[\overline{\vecc{P}}\right]_{(\overline{k}^--1)\overline{N}+i,(\overline{k}^+-1)\overline{N}+j} \coloneqq
\frac{1}{2} \sum_{a^1_{\overline{k}^-n^-}\!\!\!\!\!\in\partial \overline{T}_{\overline{k}}\cap \overline{\setV}_\Omega}
\frac{\nu_{\overline{k}^-\overline{n}^-}}{\Hs\left(a^1_{\overline{k}\overline{n}^-}\right)}
\sum_{l=1}^{\overline{N}} \overline{U}_{\overline{k}^+l}\left(a^1_{\overline{k}\overline{n}}\right)\;
\phi_{\overline{k}^-i} \left(a^1_{\overline{k}^-\overline{n}^-}\right) \, \phi_{\overline{k}^+l}\left(a^1_{\overline{k}^-\overline{n}^-}\right) \, \phi_{\overline{k}^+j}\left(a^1_{\overline{k}^-\overline{n}^-}\right) \,.
\end{equation}
\end{subequations}
Here,~$\overline{\setV}_\Omega$ is the set of interior vertices in the one-dimensional mesh.
Additionally, the same term contributes to a~vector~$\overline{\vec{K}}_{h} \in \IR^{KN}$ in the form
\begin{equation}
\label{eq:free:barKh}
\left[\overline{\vec{K}}_{h}\right]_{(\overline{k}^--1)\overline{N}+i} \coloneqq
\sum_{\Proj E_{k^-n^-} \in \partial \overline{T}_{\overline{k}^-}}
\int_{E_{k^-n^-}} \frac{1}{\Hs}
  \phi_{\overline{k}^-\overline{i}} \, \frac{|\hat{\lambda}|}{2} \, \left(
  \sum_{j=1}^{\overline{N}} H_{\overline{k}^-j} \, \phi_{\overline{k}^-j} -
  \sum_{j=1}^{\overline{N}} H_{\overline{k}^+j} \, \phi_{\overline{k}^+j} 
\right) \, \dd\sigma \,.
\end{equation}

\subsection[Contributions from domain boundary terms]{Contributions from domain boundary terms~$\VI$--$\VII$, $\IX$--$\X$, $\XIV$, $\XIX$--$\XX$, $\XXIV$}
Contributions on domain boundary edges can belong to one of the following three types:
no given boundary data -- we rely on values from the interior and thus have only contributions to matrices that are applied to representation vectors;
for specified boundary data -- either a~vector representing these boundary values (e.\,g., boundary data for flux variable~$q_\mathrm{D}$ in term $\VII$) or both, matrix and vector contributions (e.\,g., bottom boundary data for velocity components in term~$\IX$).
To keep the presentation brief for the great number of terms introduced by different boundary conditions, we cluster together some similar contributions in the following.

Term $\VI$ represents the first of the cases above: there is no boundary data available for the diffusive flux, and we use values from the interior only, thus have a~contribution to two block-diagonal matrices~$\vecc{R}^m_\mathrm{bdr}\in\IR^{KN\times KN}$,~$m\in\{1,2\}$, that take the same form as the diagonal blocks in~$\vecc{R}^m$ without the factor~$\frac{1}{2}$ (see Eq.~\eqref{eq:globR:diag}).
The same situation is also relevant for the parts involving the water height in term~$\IX$, which integrates over horizontal boundary edges, and thus no boundary data for the water height is available, and term~$\X$ for boundary parts with no Dirichlet data for~$h$.
All together result in contributions to a~rectangular block-matrix~$\check{\vecc{Q}}_\mathrm{bdr}\in\IR^{KN\times\overline{K}\overline{N}}$ with entries given by
\begin{equation*}
\left[ \check{\vecc{Q}}_\mathrm{bdr} \right]_{(k-1)N+i,(\overline{k}-1)\overline{N}+j} 
\coloneqq \sum_{E_{kn}\in\partial T_k\cap(\setE_{\partial\Omega}\setminus\setE_\IH)} \nu_{kn}^1
\int_{E_{kn}} \vphi_{ki} \, \phi_{\overline{k}j} \,\dd\sigma \,.
\end{equation*}
On vertical edges with Dirichlet data for~$h$, term~$\X$ has a~similar contribution weighted by~$\frac{1}{2}$ (due to the Lax--Friedrichs Riemann solver), and thus the entries take the same form as the left hand blocks for~$\check{\vecc{Q}}$ in Eq.~\eqref{eq:free:tildeQh}.
For that reason, we change the sum in Eq.~\eqref{eq:free:tildeQh} to include all~$E_{kn}\in\partial T_k\cap(\setE_\Omega \cup \setE_\IH)$ and assemble matrix~$\check{\vecc{Q}}$ directly for all edges except boundary edges without Dirichlet data~$h_\mathrm{D}$, which are treated by~$\check{\vecc{Q}}_\mathrm{bdr}$ above.

The same is applied to the non-linear velocity part of terms~$\IX$ and~$\X$, which contribute to a~block-diagonal matrix~$\vecc{P}^m_\mathrm{bdr}\in\IR^{KN\times KN}$,~$m\in\{1,2\}$ for edges without boundary data~$u^1_\mathrm{D}$ as 
\begin{equation*}
\left[ \vecc{P}^m_\mathrm{bdr} \right]_{(k-1)N+i,(k-1)N+j} \coloneqq
\sum_{E_{kn}\in\partial T_k\cap(\setE_{\partial\Omega}\setminus\setE_\IU)} 
\nu_{kn}^m \sum_{l=1}^N U_{kl}^1 
\int_{E_{kn}} \vphi_{ki} \, \vphi_{kl} \, \vphi_{kj} \, \dd \sigma \,,
\end{equation*}
and, for edges with boundary data, the contribution of term~$\X$ is included into matrix~$\vecc{P}^m$ (cf.~Eq.~\eqref{eq:free:Ph}) with the sum changed to include all~$E_{kn}\in\partial T_k\cap(\setE_\Omega\cup\setE^\mathrm{v}_\IU)$.

Terms~$\XIV$ and~$\XIX$ account for the contribution of Dirichlet data~$u^1_\mathrm{D}$ to block-diagonal matrices~$\vecc{Q}^m_\mathrm{bdr}\in\IR^{KN\times KN}$,~$m\in\{1,2\}$ with entries given by
\begin{equation*}
\left[\vecc{Q}^m_\mathrm{bdr}\right]_{(k-1)N+i,(k-1)N+j} \coloneqq
\sum_{E_{kn}\in\partial T_k\cap(\setE_{\partial\Omega}\setminus\setE_\IU)} 
  \nu_{kn}^m \int_{E_{kn}} \vphi_{ki} \, \vphi_{kj} \, \dd \sigma =
\sum_{E_{kn}\in\partial T_k\cap(\setE_{\partial\Omega}\setminus\setE_\IU)} 
  \nu_{kn}^m \left[ \vecc{S}_{E_{kn}} \right]_{i,j}\,.
\end{equation*}
The boundary term for~$u^2$ in~$\XIX$ is included in the assembly of~$\vecc{Q}_\mathrm{up}$ at the free surface, and, on the bottom boundaries, the Dirichlet data~$u^2_\mathrm{D}$ is always specified.

The corresponding term~$\XX$ on vertical boundary edges with no boundary data for~$h$ and~$u^1$ contributes to a~block-diagonal matrix~$\check{\vecc{P}}_\mathrm{bdr}\in\IR^{KN\times KN}$ as
\begin{equation*}
\left[\check{\vecc{P}}_\mathrm{bdr}\right]_{(k-1)N+i,(k-1)N+j} \coloneqq
\sum_{E_{kn}\in\partial T_k \cap (\setE_{\partial\Omega}^\mathrm{v}\setminus (\setE_\IU\cup\setE_\IH))}
\nu_{kn}^1 \sum_{l=1}^{\overline{N}} H_{\overline{k}l}
\int_{E_{kn}} \frac{1}{\Hs} \vphi_{ki} \, \phi_{\overline{k}l} \, \vphi_{kj} \, \dd\sigma
\end{equation*}
and has the same contribution with factor~$\frac{1}{2}$ for edges with boundary data for~$h$,~$u^1$, or both.
This is the same as the diagonal entries of~$\check{\vecc{P}}$, and thus, these contributions are assembled together with the interior edges changing the set of relevant edges in Eq.~\eqref{eq:free:tildeP:diag} to~$E_{kn}\in\partial T_k\cap(\setE_\Omega^\mathrm{v} \cup \setE^\mathrm{v}_\IU \cup \setE^\mathrm{v}_\IH)$.

In the free-surface equation, we have similar contributions due to term~$\XXIV$ with the major difference being the one-dimensional test function.
This gives a~block-diagonal matrix~$\overline{\vecc{P}}_\mathrm{bdr}\in\IR^{\overline{K}\overline{N}\times\overline{K}\overline{N}}$ for edges without Dirichlet data, where we exploit the depth-integrated velocity (cf.~Sec.~\ref{sec:variational-freeflow}) and make use of the fact that~$\phi_{\overline{k}i}$ is constant on vertical edges resulting in entries
\begin{align*}
&\left[\overline{\vecc{P}}_\mathrm{bdr}\right]_{(\overline{k}-1)\overline{N}+i,(\overline{k}-1)\overline{N}+j} \coloneqq
\sum_{E_{kn}\in\partial T_k \cap (\setE_{\partial\Omega}^\mathrm{v}\setminus (\setE_\IU\cup\setE_\IH))}
\nu_{kn}^1 \sum_{l=1}^{N} U^1_{kl}
\int_{E_{kn}} \frac{1}{\Hs} \phi_{\overline{k}i} \, \vphi_{kl} \, \phi_{\overline{k}j} \, \dd\sigma \\
&\qquad= \sum_{a^1_{\overline{k}\overline{n}}\in\partial \overline{T}_{\overline{k}} \cap(\overline{\setV}_{\partial\Omega}\setminus (\overline{\setV}_\IU\cup\overline{\setV}_\IH))}
\frac{\nu_{\overline{k}\overline{n}}}{\Hs\left(a^1_{\overline{k}\overline{n}}\right)} \sum_{m=1}^{\overline{N}} \overline{U}_{\overline{k}m}\left(a^1_{\overline{k}\overline{n}}\right) \,
\phi_{\overline{k}i}\left(a^1_{\overline{k}\overline{n}}\right) \, \phi_{\overline{k}m}\left(a^1_{\overline{k}\overline{n}}\right) \, \phi_{\overline{k}j}\left(a^1_{\overline{k}\overline{n}}\right)\,.
\end{align*}
Here,~$\overline{\setV}$ is the set of vertices in the one-dimensional mesh with the subscript indicating the same restrictions as for sets of edges, and~$a^1_{\overline{k}\overline{n}}$ is the~$\overline{n}$-th vertex ($\overline{n}\in\{1,2\}$) of one-dimensional element~$\overline{T}_{\overline{k}}$.
For edges with Dirichlet data~$u^1_\mathrm{D}$ or~$h_\mathrm{D}$, we have the same contribution as in the diagonal blocks of~$\overline{\vecc{P}}$ for interior edges, thus we change the set of relevant vertices in Eq.~\eqref{eq:free:barP:diag} to~$a^1_{\overline{k}\overline{n}}\in\partial \overline{T}_{\overline{k}} \cap(\overline{\setV}_{\Omega}\cup \overline{\setV}_\IU\cup\overline{\setV}_\IH)$ and assemble them together.

Dirichlet data in terms~$\VII$,~$\IX$,~$\X$,~$\XIV$, and~$\XIX$ enter system~\eqref{eq:free:matsys} in vectors~$\vec{J}^m_*\in\IR^{KN}$,~$m\in\{1,2\}$ that all have the same form
\begin{equation*}
\left[\vec{J}^m_{*}\right]_{(k-1)N+i} \coloneqq \sum_{E_{kn}\in\partial T_k \cap \setE_*}
\nu_{kn}^m \, \int_{E_{kn}} \vphi_{ki} \, w_\mathrm{D} \, \dd\sigma \,,
\end{equation*}
where~$\setE_*$ corresponds to the relevant set of edges for which the contribution is to be assembled, and $w_\mathrm{D}\coloneqq w_\mathrm{D}(t,\vec{x})$ represents any of the boundary data functions (e.\,g., $u^1_\mathrm{D}$, $u^2_\mathrm{D}$, $h_\mathrm{D}$).
Note that we use compositions of functions for~$w_\mathrm{D}$ whenever necessary, e.\,g.~$w_\mathrm{D}(t,\vec{x}) = u^1_\mathrm{D}(t,\vec{x})u^1_\mathrm{D}(t,\vec{x})$ in term~$\X$.
Furthermore, instead of~$q_\mathrm{D}$ we specify~$q^1_\mathrm{D}(t,\vec{x})$ and~$q^2_\mathrm{D}(t,\vec{x})$ and substitute~$q_\mathrm{D}=q^1_\mathrm{D} \nu_{kn}^1 + q^2_\mathrm{D} \nu_{kn}^2$ for boundary condition type~\eqref{eq:free:nonmixed:bc:q} in our implementation to simplify the use of analytical test examples for non-trivial free-surface boundary shapes.
For that reason, we assemble two vectors~$\vec{J}_q^1,\vec{J}_q^2$ and compute~$\vec{J}_q=\vec{J}_q^1+\vec{J}_q^2$.

In term~$\XX$, the additional scaling by~$\Hs$ and, in case only~$u^1_\mathrm{D}$ or~$h_\mathrm{D}$ is provided, the inclusion of~$h_{\Delta}$ or~$u^1_{\Delta}$ from the interior requires the vectors with Dirichlet data to be assembled differently.
Here, we have contributions to vectors~$\vec{J}_{uh},\check{\vec{J}}_u,\check{\vec{J}}_{h}\in\IR^{KN}$ as
\begin{align*}
\left[\vec{J}_{uh}\right]_{(k-1)N+i} &\coloneqq
\sum_{E_{kn}\in\partial T_k \cap \setE_\IH \cap \setE_\IU} \nu_{kn}^1
\int_{E_{kn}} \frac{1}{\Hs} \,\vphi_{ki} \, u_\mathrm{D}\, h_\mathrm{D} \, \dd\sigma \,, \\
\left[\check{\vec{J}}_u\right]_{(k-1)N+i} &\coloneqq
\sum_{E_{kn}\in\partial T_k \cap (\setE_\IU\setminus\setE_\IH)} \nu_{kn}^1
\int_{E_{kn}} \frac{1}{\Hs} \,\vphi_{ki} \, u_\mathrm{D} \left(\sum_{j=1}^{\overline{N}} H_{\overline{k}j} \, \phi_{\overline{k}j} \right) \, \dd\sigma \,, \\
\left[\check{\vec{J}}_{h}\right]_{(k-1)N+i} &\coloneqq
\sum_{E_{kn}\in\partial T_k \cap \left(\setE_\IH \setminus \setE_\IU\right)} \nu_{kn}^1
\int_{E_{kn}} \frac{1}{\Hs} \,\vphi_{ki} \, \left(\sum_{j=1}^N U_{kj}^1 \, \vphi_{kj}\right) h_\mathrm{D} \, \dd\sigma \,.
\end{align*}
Almost identical contributions are due to term~$\XXIV$ in the free-surface equation, where the only difference is the one-dimensional test function~$\phi_{\overline{k}i}$ instead of~$\vphi_{ki}$ resulting in vectors~$\overline{\vec{J}}_{uh},\overline{\vec{J}}_{u},\overline{\vec{J}}_{h}\in\IR^{\overline{K}\overline{N}}$ with entries
\begin{align*}
\left[\overline{\vec{J}}_{uh}\right]_{(\overline{k}-1)\overline{N}+i} &\coloneqq
\sum_{\Proj E_{kn}\in\partial \overline{T}_{\overline{k}} \cap \setE_\IH \cap \setE_\IU} \nu_{kn}^1
\int_{E_{kn}} \frac{1}{\Hs} \,\phi_{\overline{k}i}  \, u_\mathrm{D}\, h_\mathrm{D} \, \dd\sigma \,, \\
\left[\overline{\vec{J}}_u\right]_{(\overline{k}-1)\overline{N}+i} &\coloneqq
\sum_{\Proj E_{kn}\in\partial \overline{T}_{\overline{k}} \cap (\setE_\IU\setminus\setE_\IH)} \nu_{kn}^1
\int_{E_{kn}} \frac{1}{\Hs} \,\phi_{\overline{k}i}  \, u_\mathrm{D} \left(\sum_{j=1}^{\overline{N}} H_{\overline{k}j} \, \phi_{\overline{k}j} \right) \, \dd\sigma \,, 
\end{align*}
\begin{align*}
\left[\overline{\vec{J}}_{h}\right]_{(\overline{k}-1)\overline{N}+i} &\coloneqq
\sum_{\Proj E_{kn}\in\partial \overline{T}_{\overline{k}} \cap \left(\setE_\IH \setminus \setE_\IU\right)} \nu_{kn}^1
\int_{E_{kn}} \frac{1}{\Hs} \,\phi_{\overline{k}i}  \, \left(\sum_{j=1}^N U_{kj}^1 \, \vphi_{kj}\right) h_\mathrm{D} \, \dd\sigma \\
&= \sum_{a^1_{\overline{k}\overline{n}}\in\partial \overline{T}_{\overline{k}}\cap\left(\overline{\setV}_\IH \setminus \overline{\setV}_\IU\right)}
\frac{\nu_{\overline{k}\overline{n}}^1}{\Hs\left(a^1_{\overline{k}\overline{n}}\right)}\,
\phi_{\overline{k}i} \left(a^1_{\overline{k}\overline{n}}\right) \, 
\left(\sum_{j=1}^{\overline{N}} \overline{U}_{\overline{k}j}\left(a^1_{\overline{k}\overline{n}}\right)\,\phi_{\overline{k}j}\left(a^1_{\overline{k}\overline{n}}\right)\right)
\, h_\mathrm{D}\left(a^1_{\overline{k}\overline{n}}\right) \,.
\end{align*}

Finally, terms~$\X$,~$\XX$, and~$\XXIV$ contribute to vectors~$\vec{K}_u$,~$\vec{K}_{h}$, and~$\overline{\vec{K}}_{h}$ for edges with Dirichlet data due to the jump term in the Lax--Friedrichs Riemann solver (see Sec.~\ref{sec:riemann}).

\section{Assembly of the free-flow problem}
\label{sec:assembly}

This section outlines the necessary steps to assemble the linear system derived in Sec.~\ref{sec:free:semi-disc}.
The terms required to build the block-matrices in system~\eqref{eq:free:matsys} are transformed to the reference square $\hat{T}=[0,1]^2$ or to the reference interval~$[0,1]$, respectively, and then evaluated using numerical quadrature rules.
The assembly of the block-matrices is then performed in vectorized operations applied to local contributions.
We presented these steps in full details in previous works in series~\cite{FESTUNG1,FESTUNG2,FESTUNG3} for similar block-matrices and restrict ourselves to some general remarks about the assembly techniques and the specific changes necessary for this model.
The main differences are the use of trapezoidal elements (instead of triangles) and some model-specifics such as the need to discretize the interplay between one- and two-dimensional functions in the free-flow problem.

\subsection{Numerical integration}
\label{sec:quadrature}
We make use of the numerical quadrature rules implemented in FESTUNG and described in previous works~\cite{FESTUNG1,FESTUNG2,FESTUNG3} to approximate edge and element integrals.
All integrals over an edge~$E_{kn}\in\setE$ are transformed to the matching edge~$\hat{E}_n$ of the reference square (cf.~Eq.~\eqref{eq:trafoRule:E}) and further transformed to the reference interval~$[0,1]$ using mappings~$\hat{\vec{\gamma}}_n: [0,1] \ni s \mapsto \hat{\vec{x}} \in \hat{E}_n$, $n\in\{1,2,3,4\}$, which are given by
\begin{align*}
\hat{\vec{\gamma}}_1(s) &\coloneqq \begin{bmatrix} s \\ 0 \end{bmatrix}, & 
\hat{\vec{\gamma}}_2(s) &\coloneqq \begin{bmatrix} s \\ 1 \end{bmatrix}, & 
\hat{\vec{\gamma}}_3(s) &\coloneqq \begin{bmatrix} 1 \\ s \end{bmatrix}, & 
\hat{\vec{\gamma}}_4(s) &\coloneqq \begin{bmatrix} 0 \\ s \end{bmatrix}.
\end{align*}
These mappings are different from the corresponding mappings in previous publications in series due to the different element shape (reference square instead of reference triangle) and are provided by a~new routine~\code{gammaMapQuadri}.
As in previous works, for interior edges~$E_{k^-n^-}=E_{k^+n^+}\in \partial T_{k^-}\cap \partial T_{k^+}$, $n^-,n^+\in\{1,2,3,4\}$, where $T_{k^-}$ and $T_{k^+}$ are two adjacent elements, we introduce the mappings
\begin{subequations}
\label{eq:map:theta}
\begin{equation}
\mapEE_{n^-} : \, \hat{E}_{n^-} \ni \hat{\vec{x}} \mapsto \mapEE_{n^-}(\hat{\vec{x}}) = \vec{F}_{k^+}\circ \vec{F}_{k^-} (\hat{\vec{x}}) \in \hat{E}_{n^+}~,
\end{equation}
which maps from one side of~$\hat{T}$ to its opposite one while keeping its orientation (cf.~Fig.~\ref{fig:setT}).
Here, the local edge index~$n^+$ is given implicitly by the numbering of mesh entities (cf.~Sec.~\ref{sec:gridbasis}) and thus can be boiled down to four cases:
\begin{align}
&\mapEE_1 : \begin{bmatrix} \hat{x}^1 \\ 0 \end{bmatrix} \mapsto \begin{bmatrix} \hat{x}^1 \\ 1 \end{bmatrix} \,, &
&\mapEE_2 : \begin{bmatrix} \hat{x}^1 \\ 1 \end{bmatrix} \mapsto \begin{bmatrix} \hat{x}^1 \\ 0 \end{bmatrix} \,, &
&\mapEE_3 : \begin{bmatrix} 1 \\ \hat{x}^2 \end{bmatrix} \mapsto \begin{bmatrix} 0 \\ \hat{x}^2 \end{bmatrix} \,, &
&\mapEE_4 : \begin{bmatrix} 0 \\ \hat{x}^2 \end{bmatrix} \mapsto \begin{bmatrix} 1 \\ \hat{x}^2 \end{bmatrix} \,.
\end{align}
\end{subequations}
With the mappings in place, it suffices to define edge quadrature rules on the reference interval~$[0,1]$:
\begin{equation}
\label{eq:quad:1d}
\int_0^1 \hat{w}(s) \, \dd s \approx \sum_{r=1}^R \omega_r \,\hat{w}(\hat{q}_r)~,
\end{equation}
where $\hat{w}:[0,1]\rightarrow \IR$ with $R$~quadrature points~$\hat{q}_r\in[0,1]$ and quadrature weights~$\omega_r\in\IR$.
We rely on standard Gauss quadrature rules implemented in function~\code{quadRule1D}.

Similarly, as all integrals over~$T_k\in\setT_\Delta$ are transformed to the reference square~$\hat{T}$, it is sufficient to define quadrature rules on~$\hat{T}$.
For that, we choose a~tensor product quadrature rule built from the one-dimensional rules~\eqref{eq:quad:1d}:
\begin{equation*}
\int_{\hat{T}} \hat{w}(\hat{\vec{x}})\,\dd\hat{\vec{x}}
=  \int_0^1\int_0^1 \hat{w}(\hat{\vec{x}}) \,\dd \hat{x}^1\dd \hat{x}^2
\approx \sum_{r=1}^R\sum_{s=1}^R \omega_r \, \omega_s \, \hat{w}\left(\begin{bmatrix}\hat{q}_r\\\hat{q}_s\end{bmatrix}\right)
\end{equation*}
for~$\hat{w}:\hat{T}\rightarrow\IR$.
This is implemented in the routine~\code{quadRuleTensorProduct} that makes use of arbitrary, user-specified one-dimensional quadrature rules.

\subsection{Computing the depth-integrated velocity}
\label{sec:depth-int-vel}
To find a~representation for the depth-integrated velocity given in Sec.~\ref{sec:free:semi-disc}, we have to compute coefficients
\begin{equation*}
\overline{U}_{\overline{k}m} = \sum_{k=1}^L U_{kj(m)} \,\left(\zeta_k(x^1) - \zeta_{k-1}(x^1)\right) \,, \qquad \overline{k}\in\{1,\ldots,\overline{K}\}, \quad \overline{m}\in\{1,\ldots,\overline{N}\}
\end{equation*}
(cf.~Eq.~\eqref{eq:averagedvelocity}), where the~$x^1$-dependent element height~$\zeta_k(x^1) - \zeta_{k-1}(x^1)>0$ is given by transformation~\eqref{eq:mappings:Fk} as
\begin{equation*}
\scalemath{0.95}{%
\zeta_k(x^1) - \zeta_{k-1}(x^1)
= F_k^2( \hat{x}^1, 1 ) - F_k^2( \hat{x}^1, 0 )
= (a_{k4}^2 - a_{k1}^2) + \left((a_{k3}^2-a_{k2}^2) - (a_{k4}^2-a_{k1}^2)\right) \hat{x}^1 
= \left[ \vecc{J}_k^1 \right]_{2,2} + \left[ \vecc{J}_k^2 \right]_{2,2} \hat{x}^1
}
\end{equation*}
with~$\vec{F}_k(\hat{\vec{x}}) \eqqcolon \transpose{\left[F_k^1(\hat{x}^1),F_k^2(\hat{\vec{x}})\right]}$ and~$\hat{x}^1$ given implicitly by back-transformation~$\hat{x}^1 = \left(F_k^1\right)^{-1}(x^1)$.
Thus, we can deduce the element height directly from the Jacobian of the mapping (cf.~Eq.~\eqref{eq:mappings:Jacobian}).
To ignore the~$\hat{x}^1$-dependency as long as possible, we pass down the splitting into constant and $\hat{x}^1$-dependent parts to the coefficients of the depth-integrated velocity, i.\,e.,
\begin{equation}
\label{eq:barU}
\overline{U}_{\overline{k}m}(\hat{x}^1)
\eqqcolon \overline{U}^1_{\overline{k}m} + \overline{U}^2_{\overline{k}m} \,\hat{x}^1
= \sum_{k=1}^L U_{kj(m)} \left[ \vecc{J}_k^1 \right]_{2,2} 
  + \sum_{k=1}^L U_{kj(m)} \left[ \vecc{J}_k^2 \right]_{2,2} \hat{x}^1\,.
\end{equation}
This allows us to compute global representation vectors~$\overline{\vec{U}}^s\in\IR^{\overline{K}\overline{N}}$, $s\in\{1,2\}$ for the depth-integrated velocity as
\begin{equation}\label{eq:depthsintegraredvelcomputation}
\left[\overline{\vec{U}}^s\right]_{:,m} \coloneqq 
\underbrace{\begin{bmatrix}
\delta_{\Proj T_1=\overline{T}_1} & \delta_{\Proj T_2=\overline{T}_1} & \cdots & \delta_{\Proj T_K=\overline{T}_1} \\
\vdots & \ddots & \ddots & \vdots \\
\delta_{\Proj T_1=\overline{T}_{\overline{K}}} & \cdots & \cdots & \delta_{\Proj T_K=\overline{T}_{\overline{K}}}
\end{bmatrix}}_{\text{\code{markT2DT.'}}} \,\left(
\left[\vec{U}^1\right]_{:,j(m)} \circ \left[ \vecc{J}^s \right]_{:,2,2}
\right)\,,
\end{equation}
where~\enquote{$\circ$} denotes the Hadamard product.
The sparse~$K\times\overline{K}$ logical matrix~\code{markT2DT} in the one-dimensional grid data structure that provides the mapping between one- and two-dimensional elements.

\subsection{Assembly of element integrals}
\label{sec:assembly:elem}
The assembly of block-matrices generated by element integrals in system~\eqref{eq:free:sys} is similar to our previous works up to the modifications due to the non-stationary determinant of the Jacobian~$\vecc{J}_k$ of the mapping~$\vec{F}_k$.
Thus, the assembly of the mass matrix~$\vecc{M}$ is well-suited to highlight the necessary changes in a~compact form.

Using the transformation rule~\eqref{eq:trafoRule:T} and the determinant of the Jacobian~\eqref{eq:mappings:detJacobian}, the following holds for the local mass matrix~$\vecc{M}_{T_k}$ as defined in~\eqref{eq:globMlocM}:
\begin{align*}
&\vecc{M}_{T_k} \;=\; \sum_{s=1}^2J_k^s\,\hat{\vecc{M}}^s \\
&\quad\text{with}\quad 
\hat{\vecc{M}}^1\;\coloneqq\;
\int_{\hat{T}}\,\begin{bmatrix}
\hat{\vphi}_{1}\,\hat{\vphi}_{1} & \cdots & \hat{\vphi}_{1}\,\hat{\vphi}_{N} ~\\
 \vdots                & \ddots & \vdots \\
\hat{\vphi}_{N}\,\hat{\vphi}_{1} & \cdots & \hat{\vphi}_{N}\,\hat{\vphi}_{N} 
\end{bmatrix}\,\dd\hat{\vec{x}}\;,\;
\hat{\vecc{M}}^2 \;\coloneqq\;
\int_{\hat{T}}\,\begin{bmatrix}
\hat{\vphi}_{1}\,\hat{\vphi}_{1} & \cdots & \hat{\vphi}_{1}\,\hat{\vphi}_{N} ~\\
 \vdots                & \ddots & \vdots \\
\hat{\vphi}_{N}\,\hat{\vphi}_{1} & \cdots & \hat{\vphi}_{N}\,\hat{\vphi}_{N} 
\end{bmatrix}\,\hat{x}^1\,\dd\hat{\vec{x}}\;,
\end{align*} 
where we split the local mass matrix on the reference element into~$\hat{\vecc{M}}^1$ and~$\hat{\vecc{M}}^2$ corresponding to the constant and~$\hat{x}^1$-dependent parts of the Jacobian, respectively.
The global mass matrix~$\vecc{M}$ can then be expressed as a~Kronecker product of a~matrix containing the contributions from the determinant of the Jacobian and the local matrices~$\hat{\vecc{M}}^1, \hat{\vecc{M}}^2$:
\begin{equation*}
\vecc{M}
\;=\;
\begin{bmatrix}
\vecc{M}_{T_1} &          & \\
               & ~\ddots~ & \\
               &          & \vecc{M}_{T_K}
\end{bmatrix}
\;=\;
\sum_{s=1}^2
\begin{bmatrix}
J_1^s &          & \\
      & ~\ddots~ & \\
      &          & J_K^s
\end{bmatrix} \otimes \hat{\vecc{M}}^s\;,
\end{equation*}
where~\enquote{$\otimes$} is the operator for the Kronecker product.
The existing assembly routine presented in~\cite{FESTUNG1} \code{assembleMatElemPhiPhi} implements the assembly of a~mass-matrix for mappings with constant Jacobians, thus the only change necessary was to add support for multiple parts in the Jacobian~$J_k^s$ and matching reference blocks~$\hat{\vecc{M}}^s$.
That way, the assembly routine is able to handle arbitrary quadrilateral elements as well as triangular elements, where the determinant of the Jacobian can be non-constant.
This generic formulation allows us to re-use this routine also to assemble the one-dimensional mass matrix~$\overline{\vecc{M}}$ in Eq.~\eqref{eq:free:sys:surf} by specifying a~different grid data structure and reference block.

The same modifications were applied to other existing assembly routines related to element integrals, for example \code{assembleMatElemDphiPhi} (to assemble~$\vecc{H}^m$) and (to assemble~$\vecc{E}^m$ and~$\vecc{G}^m$) \code{assembleMatElemDphiPhiFuncDisc}.

To assemble the entries of matrix~$\overline{\vecc{G}}$ (cf.~Eq.~\eqref{eq:barG}) in term~$\XXII$ of Eq.~\eqref{eq:free:sys:surf}, we apply the transformation rule~\eqref{eq:trafoRule:barT}, use the chain rule for the partial derivative~$\partial_{x^1} = \frac{1}{|\overline{T}_{\overline{k}}|} \partial_{\hat{x}}$, and substitute the depth-integrated velocity from~\eqref{eq:barU} to obtain
{\allowdisplaybreaks%
\begin{align*}
&\Big[\overline{\vecc{G}}\Big]_{(\overline{k}-1)\overline{N}+i,(\overline{k}-1)\overline{N}+j} 
\coloneqq
\int_{\overline{T}_{\overline{k}}} \frac{1}{\Hs} \partial_{x^1} \phi_{\overline{k}i}\, \left( \sum_{l=1}^{\overline{N}} \overline{U}_{\overline{k}l} \, \phi_{\overline{k}l} \right)
\, \phi_{\overline{k}j} \, \dd x^1
= \left|\overline{T}_{\overline{k}}\right| \sum_{l=1}^{\overline{N}} 
  \int_0^1 \frac{\overline{U}^1_{\overline{k}{l}} + \overline{U}^2_{\overline{k}{l}}\, \hat{x}}{\Hs \circ \overline{F}_{\overline{k}}(\hat{x})} \,
  \frac{\partial_{\hat{x}} \hat{\phi}_i(\hat{x})}{\left|\overline{T}_{\overline{k}}\right|}  \,
  \hat{\phi}_l(\hat{x}) \, \hat{\phi}_j(\hat{x}) \, \dd\hat{x}\\
&\qquad\approx\sum_{l=1}^{\overline{N}} \sum_{r=1}^R \frac{
  \overline{U}^1_{\overline{k}{l}} \; \overbrace{
    \omega_r\, \partial_{\hat{x}}\hat{\phi}_i(q_r) \,\hat{\phi}_l(q_r) \,\hat{\phi}_j(q_r)
  }^{\eqqcolon [\hat{\overline{\vecc{G}}}^1]_{i,j,l,r}} 
  + \overline{U}^2_{\overline{k}{l}} \; \overbrace{
    \omega_r\, \partial_{\hat{x}}\hat{\phi}_i(q_r) \,\hat{\phi}_l(q_r) \,\hat{\phi}_j(q_r) \, q_r
  }^{\eqqcolon [\hat{\overline{\vecc{G}}}^2]_{i,j,l,r}}
}{\Hs\circ\overline{F}_{\overline{k}}(q_r)}
= \sum_{s=1}^2 \sum_{l=1}^{\overline{N}}
\overline{U}^s_{\overline{k}{l}} \sum_{r=1}^{R}
  \frac{\left[ \hat{\overline{\vecc{G}}}^s \right]_{i,j,l,r}}{\Hs\circ\overline{F}_{\overline{k}}(q_r)} \,,
\end{align*}}
where we applied a~one-dimensional quadrature rule (cf.~Eq.~\eqref{eq:quad:1d}).
The smoothed height~$\Hs$ in each quadrature point is mesh-dependent (due to movements of the free surface), thus it is evaluated and stored during mesh updates and used where needed afterwards.
The assembly routine is \code{assembleMatElem1DDphiPhiFuncDiscHeight}.

\subsection{Assembly of edge integrals}
\label{sec:assembly:edge}

Edge integrals resulting in entries of block-matrices~$\vecc{P}^m$,~$\vecc{Q}^m$,~$\vecc{R}^m\in\IR^{KN\times KN}$ are assembled similarly to our previous works, the only changes concern the number of edges and the implicit mapping of local edge indices~$n^-$,~$n^+$.
These are integrated into the existing assembly routines \code{assembleMatEdgePhiPhiNu} and \code{assembleMatEdgePhiPhiFuncDiscNu}, which now support both triangular and quadrilateral meshes.

The entries of matrix~$\check{\vecc{Q}}\in\IR^{KN\times\overline{K}\overline{N}}$ (cf.~Eq.~\eqref{eq:free:tildeQh}) take a~similar form as for~$\vecc{Q}^m\in\IR^{KN\times KN}$ with a~one-dimensional basis function~$\phi_{\overline{k}j}$ replacing~$\varphi_{kj}$ in the integrand.
We split the matrix into block-diagonal and off-diagonal contributions as~$\check{\vecc{Q}}\eqqcolon \check{\vecc{Q}}^\mathrm{diag} + \check{\vecc{Q}}^\mathrm{offdiag}$ and apply transformation rule~\eqref{eq:trafoRule:E}, which allows us to assemble the diagonal blocks as
\begin{equation*} 
\check{\vecc{Q}}^\mathrm{diag} = \frac{1}{2} \sum_{n=1}^4
\begin{bmatrix}
  \delta_{E_{1n}\in\setE_\Omega} & & \\
  & \ddots & \\
  & & \delta_{E_{Kn}\in\setE_\Omega}
\end{bmatrix} \circ \begin{bmatrix}
  \nu_{1n}^1 |E_{1n}| & & \\
  & \ddots & \\
  & & \nu_{1n}^1 |E_{Kn}| 
\end{bmatrix} \otimes 
\left[ \hat{\vecc{Q}}^\mathrm{diag} \right]_{:,:,n}\,,
\end{equation*} 
where~\enquote{$\circ$} denotes the Hadamard product.
The entries of~$\hat{\vecc{Q}}^\mathrm{diag}\in\IR^{N\times\overline{N}\times4}$ are given by
\begin{equation*}
\left[ \hat{\vecc{Q}}^\mathrm{diag} \right]_{i,j,n} \coloneqq
\int_{\hat{E}_n} \hat{\vphi}_i (\hat{\vec{x}}) \, \hat{\phi_j} (\hat{x}^1) \, \dd\hat{\vec{x}}\,.
\end{equation*}
Off-diagonal blocks are assembled in the same way as
\begin{align*} 
&\check{\vecc{Q}}^\mathrm{offdiag} = \frac{1}{2} \sum_{n^-=1}^4
\begin{bmatrix}
  0 & \delta_{E_{1n^-} = E_{1n^+}} & \dots & \dots & \delta_{E_{1n^-} = E_{1n^+}} \\
  \delta_{E_{2n^-} = E_{1n^+}} & 0 & \ddots  & & \vdots \\
  \vdots & \ddots & \ddots & \ddots & \vdots \\
  \vdots & & \ddots & 0 & \delta_{E_{(K-1)n^-} = E_{Kn^+}} \\
  \delta_{E_{Kn^-} = E_{1n^+}} & \dots & \dots & \delta_{E_{Kn^-} = E_{(K-1)n^+}} & 0
\end{bmatrix} \\
& \qquad \begin{bmatrix}
  \delta_{\Proj T_1=\overline{T}_1} & \delta_{\Proj T_2=\overline{T}_1} & \cdots & \delta_{\Proj T_K=\overline{T}_1} \\
\vdots & \ddots & \ddots & \vdots \\
\delta_{\Proj T_1=\overline{T}_{\overline{K}}} & \cdots & \cdots & \delta_{\Proj T_K=\overline{T}_{\overline{K}}}
\end{bmatrix} \circ \begin{bmatrix}
  \nu_{1n^-}^1 |E_{1n^-}| & \dots & \nu_{1n^-}^1 |E_{1n^-}| \\
  \vdots & & \vdots \\
  \nu_{Kn^-}^1 |E_{Kn^-}| & \dots & \nu_{Kn^-}^1 |E_{Kn^-}| 
\end{bmatrix} \otimes 
\left[ \hat{\vecc{Q}}^\mathrm{offdiag} \right]_{:,:,n}\,,
\end{align*}
with
\begin{equation*}
\left[ \hat{\vecc{Q}}^\mathrm{offdiag} \right]_{i,j,n^-} \coloneqq
\int_{\hat{E}_n^-} \hat{\vphi}_i (\hat{\vec{x}}) \, \hat{\phi_j} \circ \hat{\vartheta}^1_{n^-} (\hat{x}^1) \, \dd\hat{\vec{x}}\,.
\end{equation*}
The first matrix is easily derived from the structured grid topology: due to the local edge numbering~$n^+$ is given implicitly by~$n^-$.
The second matrix is given by~\code{markT2DT} (cf.~Eq.~\eqref{eq:depthsintegraredvelcomputation}), and the Hadamard product is carried out using the \MatOct{} routine~\code{bsxfun} without the need to assemble the full third matrix.
Here,~$\hat{\vartheta}_n^1$ denotes the first component of the mapping defined in~\eqref{eq:map:theta}.
This is implemented in routine \code{assembleMatEdgeQuadriPhiPhi1DNu}.

Essentially, the assembly of matrix~$\check{\vecc{P}}\in\IR^{KN\times KN}$ from Eq.~\eqref{eq:free:tildeP} is the same: 
we split it into block-diagonal contributions and off-diagonal blocks and have an additional sum to account for the function~$h_\Delta$ (see~$\vecc{R}^m$ in previous works).
The facts that only vertical edges have to be considered and that the smoothed height~$\Hs$ is constant for each vertical edge allow to pull it outside of the integral.
This is done in assembly routine \code{assembleMatEdgeQuadriPhiPhiFuncDisc1DNuHeight}.

The entries of matrix~$\overline{\vecc{P}}\in\IR^{\overline{K}\overline{N}\times \overline{K}\overline{N}}$ (cf.~Eq.~\eqref{eq:free:barP}) consist of evaluating basis functions~$\phi_{\overline{k}i}$ and smoothed height~$\Hs$ in the endpoints of the one-dimensional elements and multiplying them by coefficients from the representation vector of the depth-integrated velocity.
The assembly is implemented in \code{assembleMatV0T1DPhiPhiFuncDiscNuHeight}.

Entries in matrices corresponding to boundary terms, e.\,g.,~$\vecc{P}^m_\mathrm{bdr}$, $\vecc{Q}^m_\mathrm{bdr}$, $\check{\vecc{P}}_\mathrm{bdr}$, $\check{\vecc{Q}}_\mathrm{bdr}$, or~$\overline{\vecc{P}}_\mathrm{bdr}$ are assembled in the same way as the diagonal block contributions of the corresponding matrices from the interior.

\subsection{Assembly of jump terms in flux approximations}
\label{sec:assembly:jump}

The efficient assembly of vectors~$\vec{K}_u$, $\vec{K}_h$, and~$\overline{\vec{K}}_h$ (cf.~Eqs.~\eqref{eq:free:Ku},~\eqref{eq:free:Kh},~\eqref{eq:free:barKh}) that stem from the jump term in the Lax--Friedrichs Riemann solver poses a~challenge due to the eigenvalue in the integrand.
For~$\vec{K}_u$, we apply the transformation rule~\eqref{eq:trafoRule:E} and apply a~numerical quadrature rule to obtain
{\allowdisplaybreaks%
\begin{align*}
& \int_{E_{k^-n^-}} \, \vphi_{k^-i} \, \frac{|\hat{\lambda}|}{2} \left(
  \sum_{j=1}^N U_{k^-j}^1 \, \vphi_{k^-j} -
  \sum_{j=1}^N U_{k^+j}^1 \, \vphi_{k^+j} \right)\, \dd\sigma \, = \\
&\quad \abs{E_{k^-n^-}} \int_0^1 
  \hat{\vphi_{i}} \circ \hat{\vec{\gamma}}_{n^-}(s) \, 
  \frac{|\hat{\lambda} \circ \vec{F}_{k^-} \circ \hat{\vec{\gamma}}_{n^-}(s)|}{2} \, \left(
  \sum_{j=1}^N U_{k^-j}^1 \, \hat{\vphi}_{j} \circ \hat{\vec{\gamma}}_{n^-}(s) -
  \sum_{j=1}^N U_{k^+j}^1 \, \hat{\vphi}_{j} \circ \mapEE_{n^-} \circ \hat{\vec{\gamma}}_{n^-}(s) \right) \,\dd s \, \approx \\
&\qquad \sum_{r=1}^R \underbrace{
    \abs{E_{k^-n^-}} \hat{w}_r \hat{\vphi_{i}} \circ \hat{\vec{\gamma}}_{n^-}(q_r) \, \delta_{E_{k^-n^-}=E_{k^+n^+}}
  }_{\displaystyle \eqqcolon \left[\vecc{S}\right]_{(k^--1)N+i,(k^+-1)R+r,n^--2}} \, 
  \underbrace{
    \frac{|\hat{\lambda} \circ \vec{F}_{k^-} \circ \hat{\vec{\gamma}}_{n^-}(q_r)|}{2}
  }_{\displaystyle \eqqcolon \left[ \vecc{\Lambda} \right]_{(k^--1)R+r,n^--2}} \, \Bigg(
  \underbrace{
    \sum_{j=1}^N U_{k^-j}^1 \, \hat{\vphi}_{j} \circ \hat{\vec{\gamma}}_{n^-}(q_r)
  }_{\displaystyle \eqqcolon \left[ \vecc{U}^- \right]_{(k^--1)R + r, n^--2}} \\
&\qquad\qquad\qquad-\;\underbrace{
    \sum_{j=1}^N U_{k^+j}^1 \, \hat{\vphi}_{j} \circ \mapEE_{n^-} \circ \hat{\vec{\gamma}}_{n^-}(q_r)
  }_{\displaystyle \eqqcolon \left[ \vecc{U}^+ \right]_{k^+, r, n^--2}} \Bigg) \,,
\end{align*}}
where we introduced matrix~$\vecc{\Lambda}\in\IR^{KR\times2}$ that holds eigenvalues evaluated in all quadrature points of all vertical edges, tensor~$\vecc{S}\in\IR^{KN\times KR\times2}$ with test functions and coefficients, and~matrices~$\vecc{U}^-\in\IR^{KR\times2},\vecc{U}^+\in\IR^{K\times R\times2}$ that contain the horizontal velocity evaluated in each quadrature point.
With the help of these, we can assemble the vector as
\begin{equation*}
\vec{K}_u = \sum_{n=3}^4 \left[ \vecc{S} \right]_{:,:,n-2} \, \left(
  \left[ \vecc{\Lambda}_u \right]_{:,n-2} \circ \left( 
    \left[\vecc{U}^-\right]_{:,n-2} - \left[\vecc{U}_\mathrm{lin}^\pm\right]_{:,n-2} 
  \right) \right)\,.
\end{equation*}
The matrix~$\vecc{U}_\mathrm{lin}^\pm \in \IR^{KR\times 2}$ is computed with the help of a~matrix~$\vecc{U}^\pm \in \IR^{K\times R\times 2}$ that is then linearized using \MatOct's function~\code{reshape} to obtain the following form:
{\scriptsize
\begin{equation*}
\vecc{U}^\pm_\mathrm{lin} \coloneqq \begin{bmatrix}
  \left[\vecc{U}^\pm\right]_{1,1,1} & \left[\vecc{U}^\pm\right]_{1,1,2} \\
  \vdots & \vdots \\
  \left[\vecc{U}^\pm\right]_{1,R,1} & \left[\vecc{U}^\pm\right]_{1,R,2} \\
  \left[\vecc{U}^\pm\right]_{2,1,1} & \left[\vecc{U}^\pm\right]_{2,1,2} \\
  \vdots & \vdots \\
  \left[\vecc{U}^\pm\right]_{K,R,1} & \left[\vecc{U}^\pm\right]_{K,R,2} \\
\end{bmatrix}\\
\;\text{ with }\;
\left[\vecc{U}^\pm\right]_{:,:,n^--2} \coloneqq\setlength{\arraycolsep}{2pt}
\begin{bmatrix}
    0 & \delta_{E_{1n^-} = E_{1n^+}} & \dots & \dots & \delta_{E_{1n^-} = E_{1n^+}} \\
    \delta_{E_{2n^-} = E_{1n^+}} & 0 & \ddots  & & \vdots \\
    \vdots & \ddots & \ddots & \ddots & \vdots \\
    \vdots & & \ddots & 0 & \delta_{E_{(K-1)n^-} = E_{Kn^+}} \\
    \delta_{E_{Kn^-} = E_{1n^+}} & \dots & \dots & \delta_{E_{Kn^-} = E_{(K-1)n^+}} & 0
  \end{bmatrix} \left[\vecc{U}^+\right]_{:,:,n^--2}.
\end{equation*}
}

\section*{Index of notation}
\noindent
\begin{footnotesize}
\begin{tabularx}{\linewidth}{@{}lX@{}}\toprule
\textbf{Symbol}    & \textbf{Definition}\\\midrule
$\uplus$            & {Disjoint union of two sets.} \\
$[\vecc{A}]_{i,j}$  & entry in $i$th row and $j$th column of matrix~$\vecc{A}$ (analog for multidimensional arrays).\\ 
$\diag(\vecc{A},\vecc{B})$ & $\setlength\arraycolsep{2pt} \coloneqq \begin{bmatrix}\vecc{A}&\quad\\[-4pt] \quad&\vecc{B}\end{bmatrix}$, block-diagonal matrix with blocks~$\vecc{A}$, $\vecc{B}$.\\
$\card{\mathcal{M}}$& Cardinality of set~$\mathcal{M}$.\\
$\vec{a}\cdot\vec{b}$& $\coloneqq a_1 b_1 + a_2 b_2$, Euclidean scalar product in~$\IR^2$.\\
$\grad$             & $\coloneqq \transpose{[\partial_{x^1}, \partial_{x^2}]}$, spatial gradient in the physical domain~$\Omega\in\{\tO,\Omega(t)\}$.\\
$\circ$             & Composition of functions or Hadamard product.\\
$\otimes$           & Kronecker product.\\
$\vec{a}_{ki}$      & $i$th vertex of the physical element~$T_k$.\\
$\widetilde{\vecc{D}}, \vecc{D}$  & {Hydraulic conductivity, divided by specific storativity~$S_0$~$(\mathrm{m^2\,s^{-1}})$, diffusion coefficient in the free flow domain~$(\mathrm{m^2\,s^{-1}})$.}\\
$\delta_\mathrm{[condition]}$ & $\coloneqq \{1~\text{if condition is true, 0~otherwise}\}$, Kronecker delta.\\
$\Delta$            & Mesh fineness of $\setT_{\Delta}$.\\
$\eta$              & Penalty parameter.\\
$E_{kn}$, $\hat{E}_n$ & $n$th edge of the physical element~$T_k$, $n$th edge of the reference square~$\hat{T}$.\\
$\setE_{\Omega}$,\; $\setE_{\partial\Omega}$ & Set of interior edges, set of boundary edges.\\
$\setE^\mathrm{h}$,\; $\setE^\mathrm{v}$ & {Set of horizontal edges, set of vertical edges.}\\
$f$                 & {Source\,/\,sink in the free-flow domain (scalar-valued coefficient function) $(\mathrm{m\,s^{-2}})$}.\\ 
$\tilde{f}$         & {Source\,/\,sink in the subsurface domain (scalar-valued coefficient function) $(\mathrm{m\,s^{-1}})$}.\\ 
$\vec{F}_k$, $\overline{F}_{\overline{k}}$         & Affine mapping from $\hat{T}$ to $T_k$, affine mapping from $[0,1]$ to $\overline{T}_{\overline{k}}$.\\
$g$                 & {Acceleration due to gravity $(\mathrm{m\,s^{-2}})$.}\\
$\hat{\vec{\gamma}}_n$ & {Mapping from $[0,1]$ to $\hat{E}_n$.}\\
$h$                 & {Water height ($x^2$~direction), $h=\xi - \zb$ $(\mathrm{m})$.}\\
$\tilde{h}$         & {Hydraulic head ($x^2$~direction) $(\mathrm{m})$.}\\
$\widetilde{\vec{H}}$ & {$\in\IR^{KN}$, representation vector of~$\thh_{\Delta}\in\IQ_p(\setT_{\Delta})$ with respect to $\{\varphi_{kj}\}$.}\\
$\vec{H}$           & {$\in\IR^{\overline{K}\overline{N}}$, representation vector of~$h_{\Delta}\in\IQ_p(\Proj\setT_{\Delta})$ with respect to $\{\phi_{\overline{k}j}\}$}.\\
$J$                 & $\coloneqq (0,t_\mathrm{end})$, open time interval.\\
$K$                 & $\coloneqq \card{\setT_{\Delta}}$, number of elements.\\
$\vecc{K}$          & {Hydraulic conductivity $(\mathrm{m\,s^{-1}})$.}\\
$\vec{\nu}$, $\vec{\nu}_{T}$    & Unit normal on~$\partial \Omega$ pointing outward of~$\Omega$, unit normal on~$\partial T$ pointing outward of~$T$.\\
$\vec{\nu}_k$       & $\coloneqq\vec{\nu}_{T_k}$.\\
$N=N_p$             & $\coloneqq (p+1)(p+2)/2$, number of local degrees of freedom of~$\IQ_p(T)$.\\
$\omega_r$          & Quadrature weight associated with~$\hat{\vec{q}}_r$.\\
$\Omega(t)$,\; $\partial\Omega(t)$  & Free flow subdomain in two dimensions, boundary of $\Omega(t)$.\\
$\tO$,\; $\partial\tO$  & Subsurface subdomain in two dimensions, boundary of $\tO$.\\
$p$                 &  Polynomial degree.\\ 
$\phi_{\overline{k}i}$,\; $\hat{\phi}_{i}$ & $i$th basis function on~$\overline{T}_{\overline{k}}$, $i$th basis function on~$[0,1]$.\\
$\vphi_{ki}$,\;  $\hat{\vphi}_i$ & $i$th basis function on~$T_k$, $i$th basis function on~$\hat{T}$.\\
$\IQ_p(\setT_{\Delta})$    & $\coloneqq \{ w_{\Delta}:\Omega\rightarrow \IR\,;\forall T\in\setT_{\Delta},\,   {w_{\Delta}}|_T\in\IQ_p(T)\}$.\\
$\IQ_p(T)$          & Space of polynomials on~$T\in\setT_{\Delta}$ of degree at most~$p$.\\
$\widetilde{\vec{Q}}^m$,\; $\vec{Q}^m$ & {$\in\IR^{KN}$, representation vector of~$q_\Delta^m,\tq^m_{\Delta}\in\IQ_p(\setT_{\Delta})$ with respect to $\{\varphi_{kj}\}$.}\\
$\hat{\vec{q}}_r$   & $r$th quadrature point in~$\hat{T}$.\\
$\tq$               & {Flux in subsurface domain}.\\
$R$                 & Number of quadrature points.\\
\end{tabularx}

\begin{tabularx}{\linewidth}{@{}lX@{}}
$\specS$            & {Specific storativity coefficient $(\mathrm{m^{-1}})$.}\\
$t$, $t^n$, $t_\mathrm{end}$ & Time variable $(\mathrm{s})$, $n$th time level, end time.\\
$\mapEE_{n^-}$      & Mapping from $\hat{E}_{n^-}$ to $\hat{E}_{n^+}$.\\
$\Delta t^n$        & $\coloneqq t^{n+1} - t^n$, time step size.\\
$T_k$,\; $\partial T_k$  & $k$th physical (trapezoidal) element, boundary of~$T_k$.\\
$\overline{T}_{\overline{k}}$ & $\overline{k}$th physical one-dimensional element.\\
$\hat{T}$           & Reference square.\\
$\vec{u}$           & $=\transpose{[u^1,u^2]}$, {water velocity $(\mathrm{m\,s^{-1}})$.}\\
$\overline{u}^1$    & {Depth integrated velocity.}\\
$\vec{U}^m$         & {$\in\IR^{KN}$, representation vector of~$u_\Delta^m\in\IQ_p(\setT_{\Delta})$ with respect to $\{\varphi_{kj}\}$.}\\
$\vec{x}$           & $=\transpose{[x^1,x^2]}$, space variable in the physical domain~$\Omega$.\\
$\hat{\vec{x}}$     & $=\transpose{[\hat{x}^1, \hat{x}^2]}$, space variable in the reference square~$\hat{T}$.\\
$\xi$               & {Free-surface elevation~$(\mathrm{m})$.}\\
$\zb$               & {Bathymetry $(\mathrm{m})$.}\\
\bottomrule
\end{tabularx}
\end{footnotesize}


\bibliography{FESTUNG}
\bibliographystyle{elsarticle-num}

\end{document}